\RequirePackage{fix-cm}
\documentclass[smallextended]{svjour3}       
\usepackage{amssymb}
\usepackage{amsmath} \usepackage{bm} \usepackage{color}
\usepackage{multirow} \usepackage{lineno} \usepackage{aliases}
\smartqed
\usepackage{graphicx}

\begin{document}

\title{Single-step Arbitrary Lagrangian-Eulerian discontinuous Galerkin method
for 1-D Euler equations}
\author{Jayesh Badwaik \and Praveen Chandrashekar \and Christian Klingenberg}

\institute{Jayesh Badwaik\at
Dept. of Mathematics, Univ. of
W\"urzburg, W\"urzburg, Germany
              \email{badwaik.jayesh@gmail.com}           
           \and
           Praveen Chandrashekar\at
TIFR Center for Applicable
Mathematics, Bangalore, India
           \email{praveen@math.tifrbng.res.in}
           \and
          Christian Klingenberg \at
Dept. of Mathematics, Univ. of
W\"urzburg, W\"urzburg, Germany
          \email{klingenberg@mathematik.uni-wuerzburg.de}
}

\maketitle

\begin{abstract} We propose an explicit, single step discontinuous Galerkin (DG)
  method on moving grids using the arbitrary Lagrangian-Eulerian (ALE) approach
  for one dimensional Euler equations. The grid is moved with  the local fluid
  velocity modified by some smoothing, which is found to considerably reduce the
  numerical dissipation introduced by Riemann solvers.  The scheme preserves
  constant states for any mesh motion and we also study its positivity
  preservation property. Local grid refinement and coarsening are performed to
  maintain the mesh quality and avoid the appearance of very small or large
  cells. Second, higher order methods are developed and several test cases are
  provided to demonstrate the accuracy of the proposed scheme.
\end{abstract}

\section{Introduction} Finite volume schemes based on exact or approximate
Riemann solvers are used for solving hyperbolic conservation laws like the Euler
equations governing compressible flows. These schemes are able to compute
discontinuous solutions in a stable manner since they have implicit dissipation
built into them due to the upwind nature of the schemes. Higher order schemes
are constructed following a reconstruction approach combined with a high order
time integration scheme. Discontinuous Galerkin methods can be considered as
higher order generalizations of finite volume methods which also make use of
Riemann solver technology but do not need a reconstruction step since they
evolve a polynomial solution inside each cell. While these methods are formally
high order accurate on smooth solutions, they can still introduce too much
numerical dissipation in some situations. Springel~\cite{springel2010pur} gives
the example of a Kelvin-Helmholtz instability in which adding a large  constant
velocity to both states leads to suppression of the instability due to excessive
numerical dissipation. This behaviour is attributed to the fact that fixed grid
methods based on upwind schemes are not Galilean invariant. Upwind schemes, even
when they are formally high order accurate, are found to be too dissipative when
applied to turbulent flows~\cite{Johnsen:2010:AHM:1670833.1670927} since the
numerical viscosity can overwhelm the physical viscosity.

For the linear convection equation $u_t + a u_x=0$, the first order upwind
scheme has the modified partial differential equation
\[
  \df{u}{t} + a \df{u}{x}
= \half|a|h (1 - \nu) \df{^2 u}{x^2} + O(h^2), \qquad \nu = \frac{|a| \Delta
t}{h}
\]

which shows that the numerical dissipation is proportional to $|a|$
which is the wave speed. In case of Euler equations simulated with a Riemann
solver, e.g., the Roe scheme, the wave speeds are related to the eigenvalues of
the flux Jacobian and the numerical dissipation would be proportional to the
absolute values of the eigenvalues, e.g., $|v-c|, |v|, |v+c|$ where $v$ is the
fluid velocity and $c$ is the sound speed. This type of numerical viscosity is
not Galilean invariant since the fluid velocity depends on the coordinate frame
adopted for the description of the flow. Adding a large translational velocity
to the coordinate frame will increase the numerical viscosity and reduce the
accuracy of the numerical solution. Such high numerical viscosity can be
eliminated or minimized if the grid moves along with the flow as in Lagrangian
methods~\cite{doi:10.1137/0731002,Carre:2009:CLH:1552584.1553086,Maire20092391}.
However pure Lagrangian methods encounter the issue of large grid deformations
that occur in highly sheared flows as in the Kelvin-Helmholtz problem requiring
some form of re-meshing. A related approach is to use arbitrary
Lagrangian-Eulerian approach~\cite{HIRT1974227,doneaale} where the mesh velocity
can be chosen to be close to the local fluid velocity but may be regularized to
maintain the mesh quality. Even in the ALE approach, it may be necessary to
perform some local remeshing to prevent the grid quality from degrading.
In~\cite{springel2010pur}, the mesh is regenerated after every time step based
on a Delaunay triangulation, which allows it to maintain good mesh quality even
when the fluid undergoes large shear deformation. However these methods have
been restricted to second order accuracy as they rely on unstructured finite
volume schemes on general polygonal/polyhedral cells, where achieving higher
order accuracy is much more difficult compared to structured grids.

Traditionally, ALE methods have been used for problems involving moving
boundaries as in wing flutter, store separation and other problems involving
fluid structure
interaction~\cite{LOMTEV1999128},\cite{Venkatasubban19951743},\cite{wang2015aa},
\cite{Boscheri2017}, \cite{Gaburro2018}.
In these applications, the main reason to use ALE is not to minimize the
dissipation in upwind schemes but to account for the moving boundaries, and
hence the grid velocities are chosen based on boundary motion and with a view to
maintain good mesh quality. Another class of methods solve the PDE on moving
meshes where the mesh motion is determined based on a monitor function which is
designed to detect regions of large gradients in the solution,
see~\cite{ammmbook} and the references therein. These methods achieve automatic
clustering of grid points in regions of large gradients. ALE schemes have been
used to compute multi-material flows as in~\cite{Luo2004304}, since they are
useful to accurately track the material interface. The mesh velocity was chosen
to be equal to the contact speed but away from the material contact, the
velocity was chosen by linear interpolation and was not close to Lagrangian.
There are other methods for choosing the mesh velocity which have been studied
in~\cite{Boscheri2013,Lohner2004}.
Lax-Wendroff type ALE schemes for compressible flows have been developed
in~\cite{Liu20098872}. Finite volume schemes based on ADER approach have been
developed on unstructured grids~\cite{Boscheri_Dumbser_2013},
\cite{Boscheri201648}, \cite{boscheriIJNMF2014}. The theoretical analysis of
ALE-DG schemes in the framework of Runge-Kutta time stepping for conservation
laws has been done in~\cite{gero}.

In the present work, we consider only the one dimensional problem
in order to set down the fundamental principles with which in an upcoming work,
we shall solve the multidimensional problem. The numerical method developed here
will be usable in the multiple dimensions, but additional work is required in
multiple dimensions to maintain a good mesh quality under fluid flow deformation.
We develop an explicit discontinuous Galerkin scheme that is conservative on
moving meshes and automatically satisfies the geometric conservation law. The
scheme is a single step method which is achieved by using a predictor computed
from a Runge-Kutta scheme that is local to each cell in the sense that it does
not require any data from neighbouring cells and belongs to the class of schemes called ADER method. Due to the single step nature of the scheme, the TVD limiter has to be applied only once in each time step unlike in multi-stage Runge-Kutta schemes where the limiter is applied after each stage update. This nature of the ADER scheme can reduce its computational expense especially in multi-dimensional problems and while performing parallel computations. The mesh velocity is specified at
each cell face as  the local velocity with some smoothing. We
analyze the positivity of the first order scheme using Rusanov flux and derive a
CFL condition. The scheme is shown to be exact for steady moving contact waves
and the solutions are invariant to the motion of the coordinate frame. Due to
Lagrangian nature, the Roe scheme does not require any entropy fix. However, we
identify the possibility of spurious contact waves arising in some situations.
This is due to the vanishing of the eigenvalue corresponding to the contact
wave. While the cell averages are well predicted, the higher moments of the
solution can be inaccurate. This behaviour of Lagrangian DG schemes does not
seem to have been reported in the literature. We propose a fix for the
eigenvalue in the spirit of the entropy fix of Harten \cite{Harten1983} that
prevents the spurious contact waves from occurring in the solution. The
methodology developed here
will be extended to multi-dimensional flows in a future work with a view towards
handling complex sheared flows.

The rest of the paper is organized as follows. Section~(\ref{sec:euler})
introduces the Euler equation model that is used in the rest of the paper. In
section~(\ref{sec:dgfem}), we explain the derivation of the scheme on a moving
mesh together with the quadrature approximations and computation of mesh
velocity. The computation of the predicted solution is detailed in
section~(\ref{sec:pred}). The TVD type limiter is presented in
section~(\ref{sec:lim}) for a non-uniform mesh, section~(\ref{sec:pos}) shows
the positivity of the first order scheme and section~(\ref{sec:const}) shows the
preservation of constant states. The grid coarsening and refinement strategy is
explained in section~(\ref{sec:adap}) while section~(\ref{sec:res}) presents a
series of numerical results.

\section{Euler equations}
\label{sec:euler}
The Euler equations model the conservation of mass, momentum and energy, and can
be written as a system of coupled partial differential equations laws of the
form
\begin{equation}
  \df{\con}{t} + \df{\fl(\con)}{x} = 0 \label{eq:claw}
\end{equation}
where $\con$ is called the vector of {\em conserved variables} and $\fl(\con)$
are the corresponding fluxes given by
\[
  \con
  = \begin{bmatrix}
    \rho \\
    \rho v\\
    E
    \end{bmatrix},
    \qquad \fl(\con)
    = \begin{bmatrix}
        \rho v \\
        p + \rho v^2 \\
        \rho H v
      \end{bmatrix}
\]
In the above expressions, $\rho$ is the density, $v$ is the fluid velocity, $p$
is the pressure and $E$ is the total energy per unit volume, which for an ideal
gas is given by $E = p/(\gamma-1) + \rho v^2/2$, with $\gamma > 1$ being the
ratio of specific heats at constant pressure and volume, and $H=(E+p)/\rho$ is
the enthalpy. The Euler equations form a hyperbolic system; the flux Jacobian
$A(\con) = \fl'(\con)$ has real eigenvalues and linearly independent
eigenvectors. The eigenvalues are $v-c, \ v, \ v+c$ where $c=\sqrt{\gamma
p/\rho}$ is the speed of sound and the corresponding right eigenvectors are
given by
\begin{equation}
  r_1 =
  \begin{bmatrix} 1 \\ v-c \\ H - vc \end{bmatrix}, \qquad r_2 =
  \begin{bmatrix} 1 \\ v \\ \half v^2 \end{bmatrix}, \qquad r_3 =
  \begin{bmatrix} 1 \\ v + c \\ H + vc \end{bmatrix} \label{eq:eigvec}
\end{equation}
The hyperbolic property implies that $A$ can be diagonalized as
$A = R \Lambda R^{-1}$ where $R$ is the matrix formed by the right
eigenvectors as the columns and $\Lambda$ is the diagonal matrix of
eigenvalues.

\section{Discontinuous Galerkin method} \label{sec:dgfem}

\subsection{Mesh and solution space}

\begin{figure} \begin{center} \includegraphics[width=0.7\textwidth]{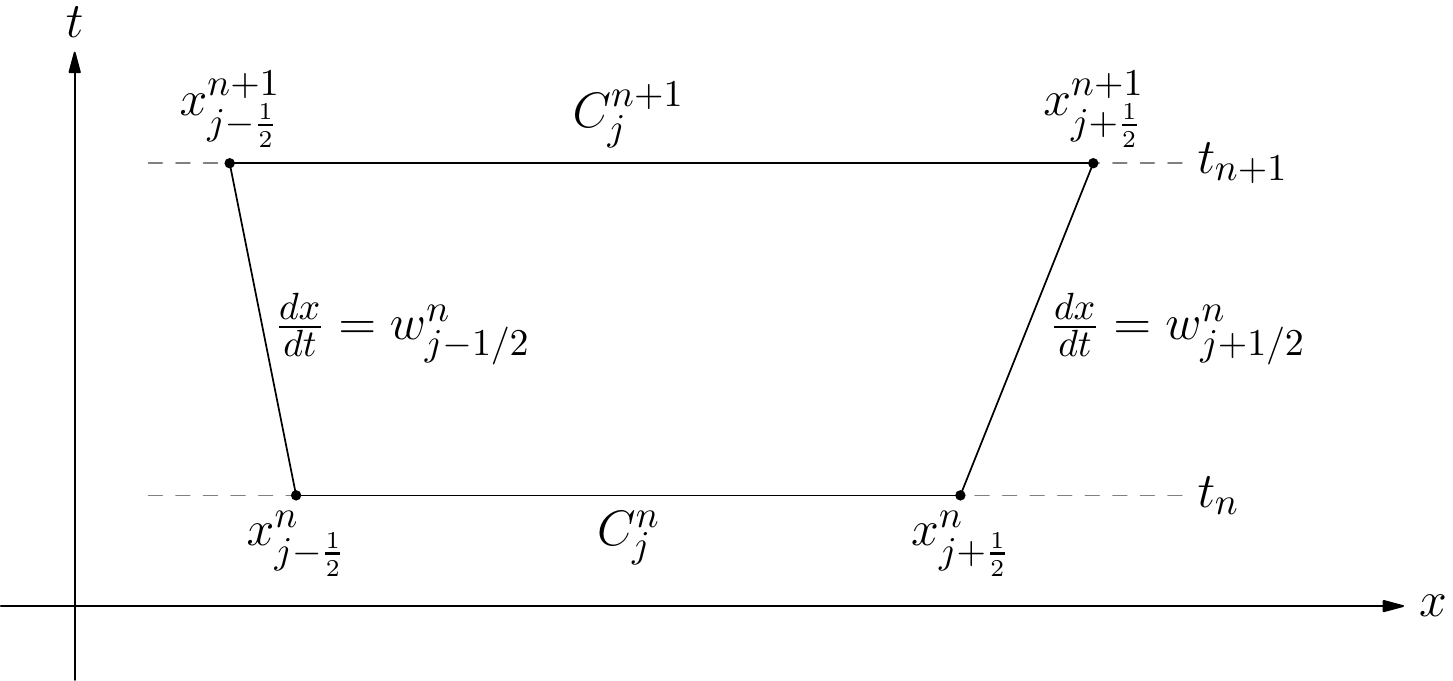}
\end{center} \caption{Example of a space-time cell in the time interval
$(t_n,t_{n+1})$} \label{fig:cell} \end{figure}

Consider a partition of the domain into disjoint cells with the $j$'th cell
being denoted by  $C_j(t) = [x_\jmh(t), x_\jph(t)]$. As the notation shows, the
cell boundaries are time dependent which means that the cell is moving in some
specified manner. The time levels are denoted by $t_n$ with the time step
$\Delta t_n = t_{n+1} - t_n$. The boundaries of the cells move with a constant
velocity in the time interval $(t_n, t_{n+1})$ given by \[ w_\jph(t) = w_\jph^n,
\qquad t_n < t < t_{n+1} \] which defines a cell in space-time as shown in
figure~(\ref{fig:cell}). The algorithm to choose the mesh velocity $w_\jph^n$ is
explained in a later section. The location of the cell boundaries is given by \[
x_\jph(t) = x_\jph^n +  (t - t_n) w_\jph^n, \qquad t_n \le t \le t_{n+1} \] Let
$x_j(t)$ and $h_j(t)$ denote the center of the cell $C_j(t)$ and its length,
i.e., \[ x_j(t) = \half( x_\jmh(t) + x_\jph(t)), \qquad h_j(t) = x_\jph(t) -
  x_\jmh(t) \] Let $w(x,t)$ be the continuous linear interpolation of the mesh
  velocity which is given by \[ w(x,t) = \frac{x_\jph(t) - x}{h_j(t)} w_\jmh^n +
    \frac{x - x_\jmh(t)}{h_j(t)} w_\jph^n, \qquad x \in C_j(t), \quad t \in
    (t_n, t_{n+1}) \] \begin{figure} \begin{center}
    \includegraphics[width=0.4\textwidth]{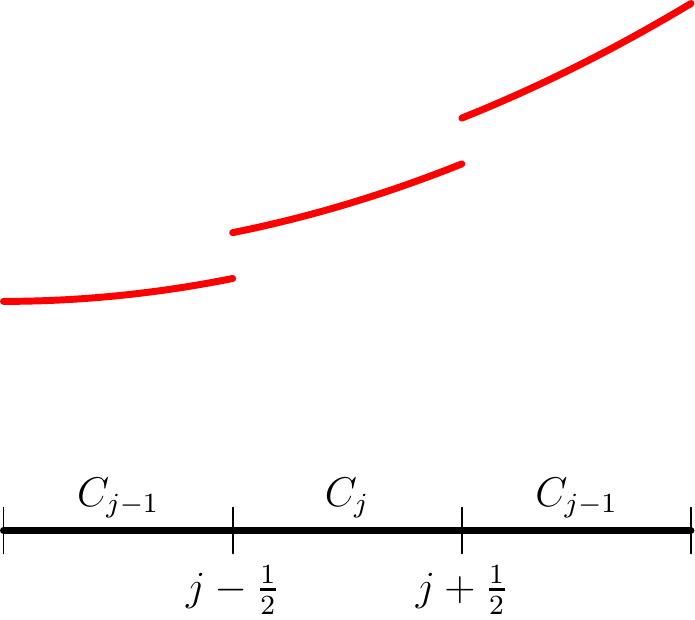} \caption{Example of a
  discontinuous piecewise polynomial solution} \label{fig:dgsol} \end{center}
\end{figure} We will approximate the solution of the conservation law by
piecewise polynomials which are allowed to be discontinuous across the cell
boundaries as shown in figure~(\ref{fig:dgsol}). For a given degree $k \ge 0$,
the solution in the $j$'th cell is given by \[ \con_h(x,t) = \sum_{m=0}^k
\con_{j,m}(t) \basis_m(x,t), \qquad x \in C_j(t) \] where $\{\con_{j,m} \in
\re^3, 0 \le m \le k\}$ are the {\em degrees of freedom} associated with the
$j$'th cell. The basis functions $\basis_m$ are defined in terms of Legendre
polynomials \[ \basis_m(x,t) = \rbasis_m(\xi) = \sqrt{2m+1} P_m(\xi), \qquad \xi
= \frac{x - x_j(t)}{\half h_j(t)} \] where $P_m :[-1,+1] \to \re$ is the
Legendre polynomial of degree $m$. The above definition of the basis functions
implies the following orthogonality property \begin{equation}
\int_{x_\jmh(t)}^{x_\jph(t)} \basis_l(x,t) \basis_m(x,t) \ud x = h_j(t)
\delta_{lm}, \qquad 0 \le l,m \le k \label{eq:ortho} \end{equation} We will
sometimes also write the solution in the $j$'th cell in terms of the reference
coordinates $\xi$ as \[ \con_h(\xi,t) = \sum_{m=0}^k \con_{j,m}(t)
\rbasis_m(\xi) \] and we will use the same notation $\con$ to denote both
functions.

\subsection{Derivation of the scheme} In order to derive the DG scheme on a
moving mesh, let us introduce the change of variable $(x,t) \to (\xi,\tau)$
given by \begin{equation} \tau = t, \qquad \xi = \frac{x - x_j(t)}{\half h_j(t)}
  \label{eq:vars} \end{equation} For any $0 \le l \le k$, we now calculate the
  rate of change of the $l$'th moment of the solution starting from
  \begin{eqnarray*} \dd{}{t} \int_{x_\jmh(t)}^{x_\jph(t)} \con(x,t)
    \basis_l(x,t) \ud x &=& \dd{}{\tau} \int_{-1}^{+1} \con(\xi,\tau)
    \rbasis_l(\xi) \half h_j(\tau) \ud\xi \\ &=& \half \int_{-1}^{+1} \left[
  h_j(\tau) \df{\con}{\tau} + \con \dd{h_j}{\tau} \right] \rbasis(\xi) \ud\xi
    \end{eqnarray*} wherein we used the change of variables given
    by~(\ref{eq:vars}). But we also have the inverse transform \[ t = \tau,
      \qquad x = x_j(\tau) + \frac{\xi}{2} h_j(\tau) \] and hence \[
      \df{t}{\tau} = 1, \qquad \df{x}{\tau} = \dd{x_j}{\tau} + \frac{\xi}{2}
      \dd{h_j}{\tau} = \frac{1}{2}(w_\jmh + w_\jph) + \frac{\xi}{2}(w_\jph -
      w_\jmh) = w(x,t) \] Using the above relations, we can easily show that \[
    \df{\con}{\tau}(\xi,\tau) = \df{\con}{t}(x,t) + w(x,t) \df{\con}{x}(x,t) \]
    Moreover \[ \dd{h_j}{\tau} = w_\jph - w_\jmh = h_j \df{w}{x} \] since
    $w(x,t)$ is linear in $x$ and hence $\partial w/\partial x$ is constant
    inside each cell. Hence the $l$'th moment evolves according to
    \begin{eqnarray*} \dd{}{t} \int_{x_\jmh(t)}^{x_\jph(t)} \con(x,t)
      \basis_l(x,t) \ud x &=& \int_{-1}^{+1} \left[ \df{\con}{t} + w
      \df{\con}{x} + \con \df{w}{x} \right] \rbasis_l(\xi) \half h_j \ud \xi \\
                          &=& \int_{x_\jmh(t)}^{x_\jph(t)} \left[ -
                        \df{\fl(\con)}{x} + \df{}{x}(w \con) \right]
                      \basis_l(x,t) \ud x \end{eqnarray*} where we have
                      transformed back to the physical coordinates and made use
                      of the conservation law~(\ref{eq:claw}) to replace the
                      time derivative of the solution with the flux derivative.
                      Define the {\em ALE flux} \begin{equation} \mfl(\con,w) =
                        \fl(\con) - w \con \label{eq:aleflux} \end{equation}
                        Performing an integration by parts in the $x$ variable,
                        we obtain \begin{eqnarray*} \dd{}{t}
                          \int_{x_\jmh(t)}^{x_\jph(t)} \con_h(x,t) \basis_l(x,t)
                          \ud x &=& \int_{x_\jmh(t)}^{x_\jph(t)} \mfl(\con_h,w)
                        \df{}{x}\basis_l(x,t) \ud x \\ &&  +
                      \nfl_\jmh(\con_h(t)) \basis_l(x_\jmh^+,t) -
                    \nfl_\jph(\con_h(t)) \basis_l(x_\jph^-,t) \end{eqnarray*}
                    where we have introduced the numerical flux \[
                    \nfl_\jph(\con_h(t)) = \nfl( \con_\jph^-(t), \con_\jph^+(t),
                  w_\jph(t)) \] which provides an approximation to the ALE flux,
                  see Appendix. Integrating over the time interval $(t_n,
                  t_{n+1})$ and using (\ref{eq:ortho}), we obtain
                  \begin{eqnarray*} h_j^{n+1} \con_{j,l}^{n+1} &=& h_j^n
                    \con_{j,l}^n + \int_{t_n}^{t_{n+1}}
                    \int_{x_\jmh(t)}^{x_\jph(t)} \mfl(\con_h,w)
                    \df{}{x}\basis_l(x,t) \ud x \ud t \\ &&  +
                    \int_{t_n}^{t_{n+1}} [\nfl_\jmh(t) \basis_l(x_\jmh^+,t) -
                    \nfl_\jph(t) \basis_l(x_\jph^-,t)] \ud t \end{eqnarray*} The
                    above scheme has an implicit nature since the unknown
                    solution $\con_h$ appears on the right hand side integrals
                    whereas we only know the solution at time $t_n$. In order to
                    obtain an explicit scheme, we assume that we have available
                    with us a {\em predicted solution} $\conp_h$ in the time
                    interval $(t_n,t_{n+1})$, which is used in the time
                    integrals to obtain an explicit scheme. Moreover, the
                    integrals are computed using quadrature in space and time
                    leading to the fully discrete scheme \begin{eqnarray}
                      \nonumber h_j^{n+1} \con_{j,l}^{n+1} &=& h_j^n
                      \con_{j,l}^n \\ &+& \Delta t_n \sum_r \wt_r h_j(\tau_r)
                      \sum_q \wx_q\mfl(\conp_h(x_q,\tau_r),w(x_q,\tau_r))
                      \df{}{x}\basis_l(x_q,\tau_r) \\ \nonumber &+& \Delta t_n
                    \sum_r \wt_r [\nfl_\jmh(\conp_h(\tau_r))
                  \basis_l(x_\jmh^+,\tau_r) - \nfl_\jph(\conp_h(\tau_r))
                \basis_l(x_\jph^-,\tau_r)] \label{eq:fulld} \end{eqnarray} where
                $\wt_r$ are weights for time quadrature and $\wx_q$ are weights
                for spatial quadrature. For the spatial integral, we will use
                $q=k+1$ point Gauss quadrature. For the time integral we will
                use mid-point rule for $k=1$ and two point Gauss quadrature for
                $k=2,3$. Since the mesh is moving,  the spatial quadrature
                points $x_q$ depend on the quadrature time $\tau_r$ though this
                is not clear from the notation. In practice, the integrals are
                computed by mapping the cell to the reference cell, and the
                basis functions and its derivatives are also evaluated on the
                reference cell. The quadrature points in the reference cell are
                independent of time due to the linear mesh evolution.

\subsection{Mesh velocity} The mesh velocity must be close to the local fluid
velocity in order to have a Lagrangian character to the scheme. Since the
solution is discontinuous, there is no unique fluid velocity at the mesh
boundaries. Some researchers, especially in the context of Lagrangian methods,
solve a Riemann problem at the cell face to determine the face velocity. Since
we use an ALE formulation, we do not require the exact fluid velocity which is
anyway not available to use since we only have a predicted solution. Following
the exact trajectory of the fluid would also lead to curved trajectories for the
grid point, which is an unnecessary complication. In our work, we make two
different choices for the mesh velocities:
\begin{enumerate}
  \item The
  first choice is to take an average of the two velocities at every face
    In the numerical results, we refer to this as ADG.
  \[
    \tilde{w}_\jph^n = \half[ v(x_\jph^-,t_n) + v(x_\jph^+, t_n)] \]

  \item
    The second choice is to solve a linearized Riemann problem
      at the face at time $t_n$.
      In the numerical results, we refer to this as RDG.
      For simplicity of notation, let the solution
      to the left of the face $x_{j+\frac{1}{2}}$ be represented as
      $u_{j+\frac{1}{2}}^-$ and the solution to the right be represented as
      $u_{j+\frac{1}{2}}^+$, then \[
      \tilde{w}_\jph^n = \frac{\rho_j^nc_j^nv_j^n +
        \rho_{j+1}^nc_{j+1}^nv_{j+1}^n}{\rho_j^nc_j^n +
        \rho_{j+1}^nc_{j+1}^n} + \frac{p_j^n - p_{j+1}^n}{\rho_j^nc_j^n + \rho_{j+1}^nc_{j+1}^n} \]

\end{enumerate} We will also perform some
  smoothing of the mesh velocity, e.g., the actual face velocity is computed
  from \[ w_\jph^n = \frac{1}{3}( \tilde{w}_{j-\half}^n + \tilde{w}_\jph^n +
    \tilde{w}_{j+\frac{3}{2}}^n) \] Note that our algorithm to choose the mesh
    velocity is very local and hence easy and efficient to implement as it does
    not require the solution of any global problems. In
    Springel~\cite{springel2010pur}, the mesh velocity is adjusted so that the
    cells remain nearly isotropic which leads to smoothly varying cell sizes.
    Such an approach leads to many parameters that need to be selected and we
    did not find a good way to make this choice that works well for a range of
    problems. Instead, we will make use of mesh refinement and coarsening to
    maintain the quality of cells, i.e., to prevent very small or large cells
    from occurring in the grid. The use of a DG scheme makes it easy to perform
    such local mesh adaptation without loss of accuracy.
\begin{remark} Consider the application of the proposed ALE-DG scheme to the
  linear advection equation $u_t + a u_x=0$. In this case the mesh velocity is
  equal to the advection velocity $w_\jph = a$, i.e., the cells move along the
  characteristics. This implies that the ALE flux $g(u,w)=au-wu=0$ and also the
  numerical flux $\hat{g}_\jph = 0$. Thus the DG scheme reduces to \[
  \int_{x_\jmh^{n+1}}^{x_\jph^{n+1}} u_h(x,t_{n+1}) \basis_l(x,t_{n+1}) \ud x =
\int_{x_\jmh^n}^{x_\jph^n} u_h(x,t_n) \basis_l(x,t_{n}) \ud x, \qquad
l=0,1,\ldots,k \] so that the solution at time $t_n$ has been advected exactly
to the solution at time $t_{n+1}$. Note that there is no time step restriction
involved in this case and the accuracy of the predicted solution is also not
relevant. If the initial condition has a discontinuity coinciding with a cell
face, then the scheme advects the discontinuity exactly without any diffusion.
\end{remark}
\section{Computing the predictor} \label{sec:pred} The predicted solution is
used to approximate the flux integrals over the time interval $(t_n,t_{n+1})$
and the method to compute this must be local, i.e., it must not require solution
from neighbouring cells. Several methods for computing the predictor have been
reviewed in~\cite{Gassner20114232}. The simplest approach is to use a Taylor
expansion in space and time. Since the cells are moving, the Taylor expansion
has to be performed along the trajectory of the mesh motion. For a second order
scheme, an expansion retaining only linear terms in $t$ and $x$ is sufficient.
Consider a quadrature point $(x_q,\tau_r)$; the Taylor expansion of the solution
around the cell center $x_j^n$ and time level $t_n$ yields \begin{eqnarray*}
  \con_h(x_q,\tau_r) &=& \con_h(x_j^n,t_n) + (\tau_r - t_n)
  \df{\con_h}{t}(x_j^n,t_n) + (x_q - x_j^n) \df{\con_h}{x}(x_j^n,t_n) \\ && +
  O(\tau_r-t_n)^2 + O(x_q-x_j^n)^2 \end{eqnarray*} and the predicted solution is
  given by truncating the Taylor expansion at linear terms, leading to \[
    \conp(x_q,\tau_r) = \con_h(x_j^n,t_n) + (\tau_r - t_n)
    \df{\con_h}{t}(x_j^n,t_n) + (x_q - x_j^n) \df{\con_h}{x}(x_j^n,t_n) \] Using
    the conservation law, the time derivative is written as $\df{\con}{t} =
    -\df{\fl}{x} = - A \df{\con}{x}$ so that the predictor is given by
    \begin{equation} \conp_h(x_q,\tau_r) = \con_h^n(x_j^n) - (\tau_r - t_n)
      \left[ A(\con_h^n(x_j^n)) - w_q I \right] \df{\con_h^n}{x}(x_j^n)
    \label{eq:pred1} \end{equation} The above predictor is used for the case of
    polynomial degree $k=1$. This procedure can be extended to higher orders by
    including more terms in the Taylor expansion but the algebra becomes
    complicated. Instead we will adopt the approach of continuous explicit
    Runge-Kutta (CERK) schemes~\cite{Owren:1992:DEC:141072.141096} to
    approximate the predictor.

Let us choose a set of $(k+1)$ distinct nodes, e.g., Gauss-Legendre or
Gauss-Lobatto nodes, which uniquely define the polynomial of degree $k$. These
nodes are moving with velocity $w(x,t)$, so that the time evolution of the
solution at node $x_m$ is governed by \begin{eqnarray*} \dd{\conp_m}{t} &=&
  \df{}{t}\conp_h(x_m,t) + w(x_m,t) \df{}{x}\conp_h(x_m,t) \\ &=&
  -\df{}{x}\fl(\conp_h(x_m,t)) + w(x_m,t) \df{}{x} \conp_h(x_m,t) \\ &=& -[
  A(\conp_m(t)) - w_m(t) I] \df{}{x}\conp_h(x_m,t) =: \rhs_m(t) \end{eqnarray*}
  wherein we have made use of the PDE to write the time derivative in terms of
  spatial derivative of the flux. This equation is solved with initial condition
  \[ \conp_m(t_n) = \con_h(x_m, t_n) \] Using a Runge-Kutta scheme of sufficient
  order (see Appendix), we will approximate the solution at these nodes as \[
    \conp_m(t) = \con_h(x_m, t_n) + \sum_{s=1}^{n_s} b_s((t-t_n)/\Delta t_n)
    \rhs_{m,s}, \quad t \in [t_n,t_{n+1}), \quad m=0,1,\ldots,k \] where
    $\rhs_{m,s} = \rhs_m(t_n + \theta_s\Delta t_n)$, $\theta_s\Delta t_n$ is the
    stage time and $b_s$ are certain polynomials related to the CERK scheme and
    given in the Appendix. Note that we are evolving the nodal values but the
    computation of $\rhs_{m,s}$ requires the modal representation of the
    solution in order to calculate spatial derivative of the solution.

\begin{figure} \begin{center} \includegraphics[width=0.7\textwidth]{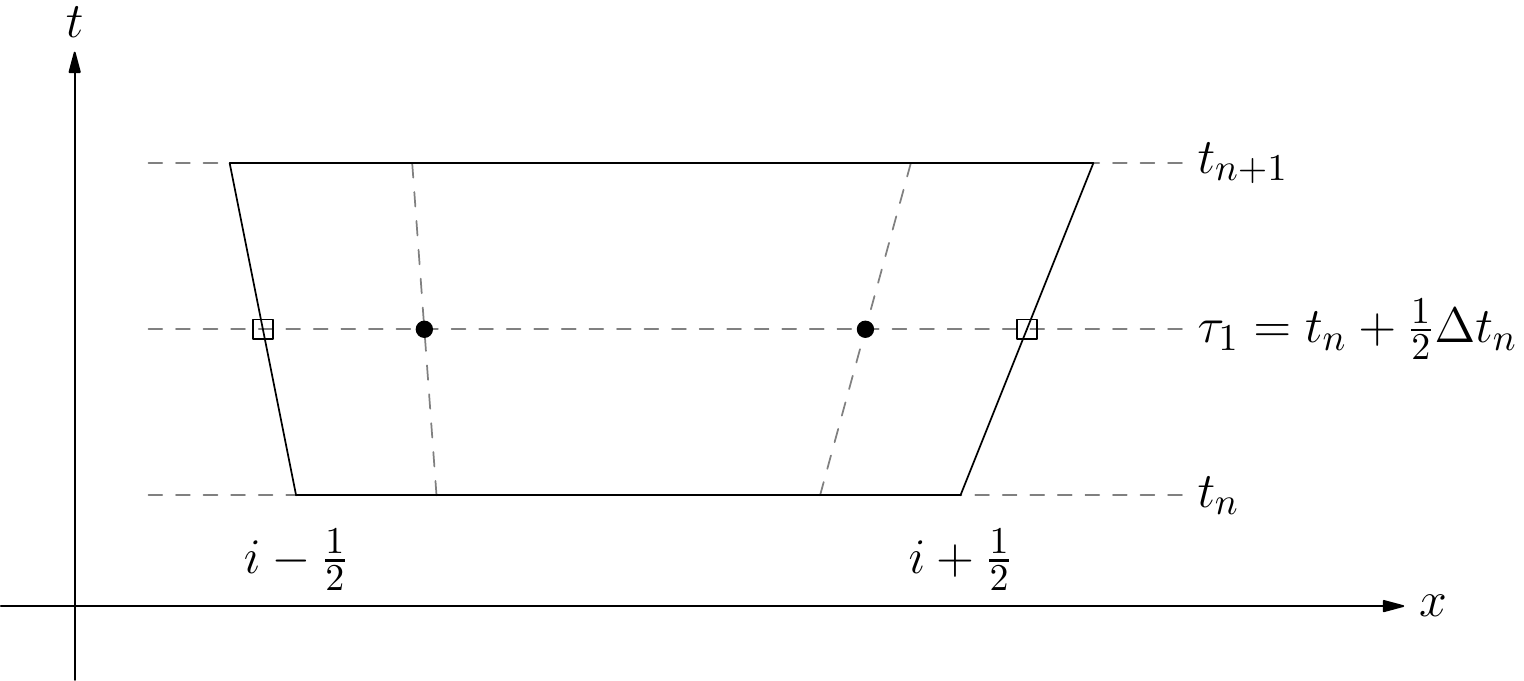}
\end{center} \caption{Quadrature points for second order scheme}
\label{fig:cerk1} \end{figure}

\begin{figure} \begin{center} \includegraphics[width=0.7\textwidth]{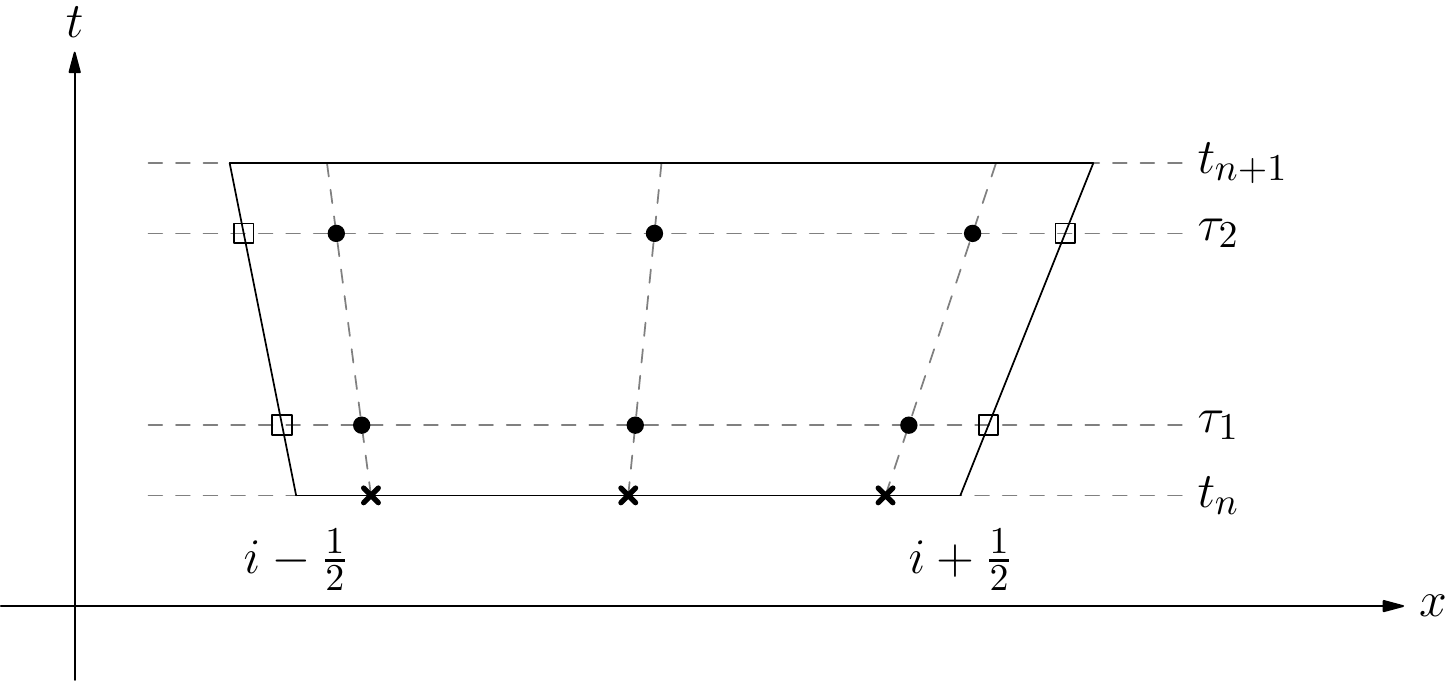}
\end{center} \caption{Quadrature points for third order scheme}
\label{fig:cerk2} \end{figure}

Once the predictor is computed as above, it must be evaluated at the quadrature
point $(x_q, \tau_r)$ as follows. For each time quadrature point $\tau_r \in
(t_n, t_{n+1})$, \begin{enumerate} \item Compute nodal values
  $\conp_m(\tau_r)$, $m=0,1,\ldots,k$ \item For each $r$, convert the nodal
  values to modal coefficients $\con_{m,r}$, $m=0,1,\ldots,k$ \item Evaluate
  predictor $\conp_h(x_q,\tau_r) = \sum_{m=0}^k \con_{m,r} \basis_m(x_q,
  \tau_r)$ \end{enumerate} The conversion from nodal to modal values is
  accomplished through a Vandermonde matrix of size $(k+1) \times (k+1)$ which
  is the same for every cell and can be inverted once before the iterations
  start. The predictor is also computed at the cell boundaries using the above
  procedure. Figure~(\ref{fig:cerk1}) and (\ref{fig:cerk2}) show the quadrature
  points used in the second, third and fourth order scheme. For the second order
  scheme, the values at $\bullet$ and $\Box$ points are obtained from the
  predictor based on Taylor expansion as given in equation~(\ref{eq:pred1}). For
  third and fourth order schemes, the nodal values $\times$ are evolved forward
  in time by the CERK scheme and evaluated at the $\bullet$ points. The
  $\bullet$ point values are converted to modal coefficients using which the
  solution at the $\Box$ points is computed.
\begin{remark} By looping over each cell in the mesh, the predicted solution is
  computed in each cell and the cell integral in~(\ref{eq:fulld}) is evaluated.
  The trace values $\conp_\jph^-(\tau_r)$ and $\conp_\jmh^+(\tau_r)$ at the
  $\Box$ points needed for quadrature in time are computed and stored. These are
  later used in a loop over the cell faces where the numerical flux is
  evaluated. Thus the algorithm is easily parallelizable on multiple core
  machines and/or using threads.  \end{remark}
\section{Control of Oscillation by Limiting} \label{sec:lim} High order schemes for hyperbolic equations
suffer from spurious numerical oscillations when discontinuities or large
gradients are present in the solution which cannot be accurately resolved on the
mesh. In the case of scalar problems, this is a manifestation of loss of TVD
property and hence limiters are used to satisfy some form of TVD condition. In
the case of DG schemes, the limiter is used as a post processor which is applied
on the solution after the time update has been performed. If the limiter detects
that the solution is oscillatory, then the solution polynomial is reduced to
at most a linear polynomial with a limited slope. In the present scheme, the
limiter is applied after the solution is updated from time $t_n$ to time
$t_{n+1}$, i.e., the solution $\con_h^{n+1}$ obtained from~(\ref{eq:fulld}) is
post-processed by the limiter. Since the mesh is inherently non-uniform due to
it being moved with the flow, we modify the standard TVD limiter to account for
this non-uniformity. Also, since we are solving a system of conservation laws,
the limiter is applied on the local characteristic variables which gives better
results than applying it directly on the conserved
variables~\cite{Cockburn:1989:TRL:69978.69982}.

The solution in cell $j$ can be written as \footnote{We suppress the time
variable for clarity of notation.}
\[
  \con_h(x) = \acon_j + \frac{x - x_j}{\half h_j} \slope_j + \textrm{higher order terms}
\]

where $\acon_j$ is the
cell average value and $\slope_j$ is proportional to the derivative of the
solution at the cell center. Let $R_j$, $L_j$ denote the matrix of right and
left eigenvectors evaluated at the cell average value $\acon_j$ which satisfy
$L_j = R_j^{-1}$, and the right eigenvectors are given in
equation~(\ref{eq:eigvec}). The local characteristic variables are defined by
$\acon^* = L_j \acon$ and $\slope^* = L_j \slope$. We first compute the limited
slope of the characteristic variables from \[ \slope_j^{**} = h_j \ \minmod{
\frac{\slope_j^*}{h_j}, \frac{\acon^*_j - \acon^*_{j-1}}{\half(h_{j-1} + h_j)},
\frac{\acon^*_{j+1} - \acon^*_{j}}{\half(h_{j} + h_{j+1})} } \] where the minmod
function is defined by \[ \minmod{a,b,c} = \begin{cases} s \min(|a|,|b|,|c|) &
\textrm{if } s = \sign{a} = \sign{b} = \sign{c} \\ 0 & \textrm{otherwise}
\end{cases} \] If $\slope_j^{**} = \slope_j^*$ then we retain the solution as it
is, and otherwise, the solution is modified to \[ \con_h(x) = \acon_j + \frac{x
- x_j}{\half h_j} R_j \slope_j^{**} \] This corresponds to a TVD limiter which
is known to lose accuracy at smooth extrema~\cite{cockburn1989} since the minmod
function returns zero slope at local extrema. The TVB limiter corresponds to
replacing the minmod limiter function with the following function
\begin{equation} \minmodb{a,b,c} = \begin{cases} a & \textrm{if } |a| \le M h^2
\\ \minmod{a,b,c} & \textrm{otherwise} \end{cases} \label{eq:tvblim}
\end{equation} where the parameter $M$ is an estimate of the second derivative
of the solution at smooth extrema~\cite{cockburn1989}, and has to be chosen by
the user.
\section{Positivity property} \label{sec:pos} The solutions of Euler equations
are well defined only if the density and pressure are positive quantities. This
is not a priori guaranteed by the DG scheme even when the TVD limiter is
applied. In the case of Runge-Kutta DG schemes, a positivity limiter has been
developed in~\cite{Zhang:2010:PHO:1864819.1865066} which preserves accuracy in
smooth regions. This scheme is built on a positive first order finite volume
scheme. Consider the first order version of the ALE-DG scheme which is a finite
volume scheme given by \begin{equation} h_j^{n+1} \acon_j^{n+1} = h_j^n
\acon_j^n - \Delta t_n [ \nfl_\jph^n - \nfl_\jmh^n] \label{eq:fvm}
\end{equation} The only degree of freedom is the cell average value and the
solution is piecewise constant. We will analyze the positivity of this scheme
for the case of Rusanov flux which is given in~\ref{sec:rusanov}. The update
equation can be rewritten as \begin{eqnarray*} h_j^{n+1} \acon_j^{n+1} &=&
  \left[ h_j^n - \frac{\Delta t_n}{2}( \lambda_\jmh^n + w_\jmh^n +
  \lambda_\jph^n - w_\jph^n) \right] \acon_j^n \\ && + \frac{\Delta t_n}{2}
  \left[ (\lambda_\jmh^n - w_\jmh^n) \acon_{j-1}^n + \fl_{j-1}^n  \right] \\ &&
  + \frac{\Delta t_n}{2} \left[ (\lambda_\jph^n + w_\jph^n) \acon_{j+1}^n -
\fl_{j+1}^n \right] \\ &=& a_j^n \acon_j^n + \frac{\Delta t_n}{2} B_j^n +
\frac{\Delta t_n}{2} C_j^n \end{eqnarray*} From the definition of the Rusanov
flux formula, we can easily see that\footnote{We drop the superscript $^n$ in
  some of these expressions.} \[ (\lambda_\jmh - w_\jmh) + v_{j-1} \ge c_{j-1} >
  0, \qquad (\lambda_\jph + w_\jph) - v_{j+1} \ge c_{j+1} > 0 \] Consider the
  first component of $B_j^n$ \begin{eqnarray*} (\lambda_\jmh - w_\jmh)
    \rho_{j-1} + \rho_{j-1} v_{j-1} &\ge& (|v_{j-1} - w_\jmh| + c_{j-1} - w_\jmh
    + v_{j-1}) \rho_{j-1} \\ &\ge& c_{j-1} \rho_{j-1} > 0 \end{eqnarray*}
    Consider the first component of $C_j^n$ \begin{eqnarray*} (\lambda_\jph +
      w_\jph) \rho_{j+1} - \rho_{j+1} v_{j+1} &\ge& (|v_j - w_\jph| + c_{j+1} +
      w_{j-1} - v_{j+1}) \rho_{j+1} \\ &\ge& c_{j+1} \rho_{j+1} > 0
    \end{eqnarray*} The pressure corresponding to $B_j^n$ is \[ p_{j-1} \left(
      -p_{j-1} + \frac{2\rho_{j-1} (v_{j-1} + \lambda_\jmh -
    w_\jmh)^2}{\gamma-1} \right) \ge p_{j-1} \rho_{j-1} c_{j-1}^2
    \frac{\gamma+1}{\gamma-1} \ge 0 \] and similarly the pressure corresponding
    to $C_j^n$ is non-negative. Hence if the coefficient in term $a_j^n$ is
    positive, then the scheme is positive. This requires the CFL condition \[
      \Delta t_n \le \frac{2h_j^n}{(\lambda_\jmh^n + w_\jmh^n + \lambda_\jph^n -
      w_\jph^n)} \] The time step will also be restricted to ensure that the
      cell size does not change too much in one time step. If we demand that the
      cell size does not change by more than a fraction $\beta \in (0,1)$, then
      we need to ensure that the time step satisfies \[ \Delta t_n \le
      \frac{\beta h_j^n}{|w_\jph^n - w_\jmh^n|} \] Combining the previous two
      conditions, we obtain the following condition on the time step
      \begin{equation} \Delta t_n \le \Delta t_n^{(1)} := \min_j\left\{
        \frac{(1-\thalf\beta) h_j^n}{\thalf(\lambda_\jmh^n + \lambda_\jph^n)},
      \frac{\beta h_j^n}{|w_\jph^n - w_\jmh^n|} \right\} \label{eq:dtpos}
    \end{equation} We can now state the following result on the positivity of
    the first order finite volume scheme on moving meshes.
\begin{theorem} The scheme~(\ref{eq:fvm}) with Rusanov flux is positivity
preserving if the time step condition~(\ref{eq:dtpos}) is satisfied.
\end{theorem}

\begin{remark} In this work we have not attempted to prove the positivity of the
  scheme for other numerical fluxes. We also do not have a proof of positivity
  for higher order version of the scheme. In the computations, we use the
  positivity preserving limiter of~\cite{Zhang:2010:PHO:1864819.1865066} which
  leads to robust schemes which preserve the positivity of the cell average
  value in all the test cases.  \end{remark}
\section{Preservation of constant states} \label{sec:const} An important
property of schemes on moving meshes is their ability to preserve constant
states for any mesh motion. This is related to the conservation of cell volumes
in relation to the mesh motion. In our scheme, if we start with a constant state
$\con^n_h = \const$, then the predictor is also constant in the space-time
interval, i.e., $\conp_h = \const$. The space-time terms in (\ref{eq:fulld}) are
polynomials with degree $k+1$ in space and degree one in time and these are
exactly integrated by the chosen quadrature rule. The flux terms at cell
boundaries in (\ref{eq:fulld}) are of degree one in time and these are also
exactly integrated. Hence the scheme (\ref{eq:fulld}) can be written as
\begin{eqnarray*} h_j^{n+1} \con_{j,l}^{n+1} &=& h_j^n \con_{j,l}^n +
  \int_{t_n}^{t_{n+1}} \int_{x_\jmh(t)}^{x_\jph(t)} \mfl(\const,w)
  \df{}{x}\basis_l(x,t) \ud x \ud t \\ &&  + \int_{t_n}^{t_{n+1}} [\nfl_\jmh(t)
\basis_l(x_\jmh^+,t) - \nfl_\jph(t) \basis_l(x_\jph^-,t)] \ud t \end{eqnarray*}
where $\con_{j,0}^n = \const$ and $\con_{j,l}^n = 0$ for $l=1,2,\ldots,k$. Due
to the constant predictor and by consistency of the numerical flux \[
\nfl_\jph(t) = \fl(\const) - w_\jph^n \const \] Moreover, for $l=1,2,\ldots,k$
\begin{eqnarray*} && \int_{t_n}^{t_{n+1}} \int_{x_\jmh(t)}^{x_\jph(t)}
  \mfl(\const,w) \df{}{x}\basis_l(x,t) \ud x \ud t   + \int_{t_n}^{t_{n+1}}
  [\nfl_\jmh(t) \basis_l(x_\jmh^+,t) - \nfl_\jph(t) \basis_l(x_\jph^-,t)] \ud t
  \\ &=& \int_{t_n}^{t_{n+1}} \int_{x_\jmh(t)}^{x_\jph(t)} \df{}{x}
  \mfl(\const,w) \basis_l(x,t) \ud x \ud t   + \int_{t_n}^{t_{n+1}}
  [\nfl_\jmh(t) \basis_l(x_\jmh^+,t) - \nfl_\jph(t) \basis_l(x_\jph^-,t)] \ud t
  \\ && - \int_{t_n}^{t_{n+1}} \int_{x_\jmh(t)}^{x_\jph(t)} \basis_l(x,t)
  \df{}{x} \mfl(\const,w)  \ud x \ud t \\ &=& - \int_{t_n}^{t_{n+1}} \df{}{x}
\mfl(\const,w) \left( \int_{x_\jmh(t)}^{x_\jph(t)} \basis_l(x,t) \ud x   \right)
\ud t =0 \end{eqnarray*} where we have used the property that $w$ is an affine
function of $x$ and $\basis_l$ are orthogonal. This implies that
$\con_{j,l}^{n+1}=0$ for $l=1,2,\ldots,k$. For $l=0$, we get \begin{eqnarray*}
h_j^{n+1} \con_{j,0}^{n+1} &=& h_j^n \const  + \int_{t_n}^{t_{n+1}}
[\nfl_\jmh(t) - \nfl_\jph(t) ] \ud t = [ h_j^n + (w_\jmh^n - w_\jph^n) \Delta t]
\const \end{eqnarray*} and since $h_j^{n+1} = h_j^n + (w_\jmh^n - w_\jph^n)
\Delta t$, we obtain $\con_{j,0} = \const$ which implies that $\con_h^{n+1} =
\const$.

\section{Grid coarsening and refinement} \label{sec:adap} The size of the cells
can change considerably during the time evolution         process   due to the
near Lagrangian movement of the cell boundaries. Near      shocks, the cells
will be compressed to smaller sizes which will reduce the     allowable time
step since a CFL condition has to be satisfied. In some regions, e.g., inside
expansion fans, the cell size can increase considerably which may  lead to loss
of accuracy. In order to avoid too small or too large cells from   occurring in
the grid, we implement cell merging and refinement into our scheme.  If a cell
becomes smaller than some specified size $h_{min}$, then it is merged with one
of its neighbouring cells and the solution is transferred from the two cells to
the new cell by performing an $L^2$ projection.  If a cell size becomes larger
than some specified size $h_{max}$, then this cell is refined into two cells by
division and the solution is again transferred by $L^2$ projection. The use of
$L^2$ projection for solution transfer ensures the conservation of mass,
momentum and energy and preserves the accuracy in smooth regions. We also ensure
that the cell sizes do not change drastically   between neighbouring cells. To
keep a track of refinement of cells, each cell is assigned an initial level
equal to $0$. The daughter cells created during       refinement are assigned a
level incremented from the parent cell, while the     coarsened cells are
assigned a level decremented from the parent cell.

The algorithm for refinement and coarsening is carried out in three sweeps
over all the active cells. In the first sweep, we mark the cells for
refinement or coarsening based on their size and the level of neighboring
cells.  Cells are marked for coarsening if the size is less than a pre-specified
minimum size. They are marked for refinement if either the size of    the cell
is larger than the maximum size or if the level of the cell is less    than the
level of the neighboring cells. If none of the conditions are
satisfied, the cells are marked for no change.  In the second sweep, a cell is
marked for refinement if both the neighboring    cells are marked for
refinement. A cell is also marked for refinement if the    size of the cell is
larger than twice the size of either of the neighboring     cells, and is also
larger than twice the minimum size. The last condition is    inserted in order
to prevent a cell being alternately marked for refinement and coarsening in
consecutive adaptation cycles.  In the third and final sweep, we again mark
cells for refinement if both the    neighboring cells are marked for refinement.
Further, we ensure that a cell     marked for refinement does not have a
neighboring cell marked for a coarsening, since this can lead to an inconsistent
mesh.
\section{Numerical results} \label{sec:res} The numerical tests are performed
with polynomials of degree one, two and three, together with the linear Taylor
expansion, two stage CERK and four stage CERK, respectively, for the computation
of the predictor. For the quadrature in time, we use the mid-point rule, two and
three point Gauss-Legendre quadrature, respectively. The time step is chosen
using the CFL condition, \[ \Delta t_n = \frac{\cfl}{2k+1} \Delta t_n^{(1)} \]
where $\Delta t_n^{(1)}$ is given by equation~(\ref{eq:dtpos}), and the factor
$(2k+1)$ comes from linear stability analysis~\cite{cockburn1989}. In most of
the computations we use $\cfl = 0.9$ unless stated otherwise. We observe that
the results using average or linearized Riemann velocity are quite similar.
We use the average velocity for most of the results and show the comparison
between the two velocities for some results.  The main steps in the algorithm
within one time step $t_n \to t_{n+1}$ are as follows.  \begin{enumerate} \item
  Choose mesh velocity $w_\jph$ \item Choose time step $\Delta t_n$ \item
  Compute the predictor $\conp_h$ \item Update solution $\con_h^n$ to the next
  time level $\con_h^{n+1}$ \item Apply TVD/TVB limiter on $\con_h^{n+1}$ \item
  Apply positivity limiter on $\con_h^{n+1}$ from
  \cite{Zhang:2010:PHO:1864819.1865066}. \item Perform grid
  refinement/coarsening \end{enumerate} In all the solution plots given below,
  symbols denote the cell average value.
\subsection{Order of accuracy} We study the convergence rate of the schemes by
applying them to a problem with a known smooth solution. The initial condition
is taken as \[ \rho(x,0) = 1 + \exp(-10 x^2), \qquad u(x,0) = 1, \qquad p(x,0) =
1 \] whose exact solution is $\rho(x,t) = \rho(x-t,0)$, $u(x,t)=1$, $p(x,t)=1$.
The initial domain is $[-5,+5]$ and the final time is $t=1$ units. The results
are presented using Rusanov and HLLC numerical fluxes. The $L^2$ norm of the
error in density are shown in table~(\ref{tab:ord0}),~(\ref{tab:hllcord0}) for
the static mesh and in table~(\ref{tab:ord1}),~(\ref{tab:hllcord1}) for the
moving mesh. In each case, we see that the error behaves as $O(h^{k+1})$ which
is the optimal rate we can expect for smooth solutions.
\textcolor{blue}{In table~(\ref{tab:limitedord0}), we show that  the ALE DG
methods preserves its higher order in presence of a limiter.}

\begin{figure} \begin{center}
\includegraphics[width=0.5\textwidth]{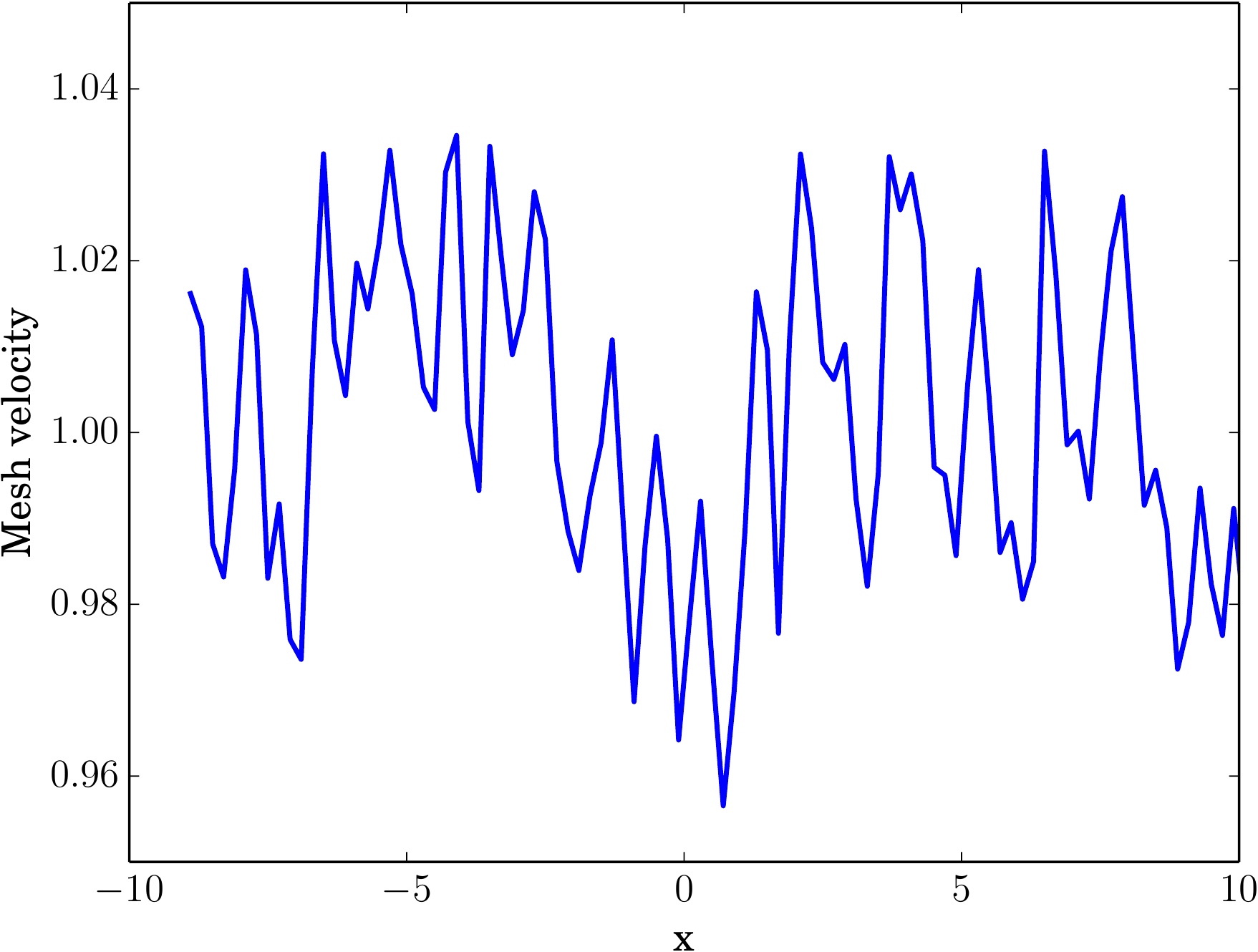} \caption{Example
of randomized velocity distribution for smooth test case} \label{fig:randvel}
\end{center} \end{figure}

The mesh velocity is constant since the fluid velocity is constant. In order to
study the effect of perturbations in mesh velocity, we add a random perturbation
to each mesh velocity, $w_\jph \leftarrow (1 + \alpha r_\jph)w_\jph$ where
$r_\jph$ is a uniform random variable in $[-1,+1]$ and $\alpha=0.05$ and a
sample velocity distribution is shown in figure~(\ref{fig:randvel}). Note that
this randomization is performed in each time step with different random
variables drawn for each face. For moving mesh, there is no unique cell size,
and the convergence rate is computed based on initial mesh spacing which is
inversely proportional to the number of cells. From table~(\ref{tab:hllcord2})
which shows results using HLLC flux, we again observe that the error reduces at
the optional rate of $k+1$ even when the mesh velocity is not very smooth.

\begin{table} \begin{center} \begin{tabular}{|c|c|c|c|c|c|c|} \hline
  \multirow{2}{*}{$N$}  & \multicolumn{2}{|c|}{$k=1$} &
  \multicolumn{2}{|c|}{$k=2$} & \multicolumn{2}{|c|}{$k=3$} \\ \cline{2-7} &
  Error & Rate & Error & Rate & Error & Rate \\ \hline 100 &4.370E-02& -
        &3.498E-03 & - &3.883E-04& - \\ 200 &6.611E-03 &2.725 &4.766E-04 &2.876
        &1.620E-05 &4.583 \\ 400 &1.332E-03 &2.518 &6.415E-05 &2.885 &9.376E-07
        &4.347\\ 800 &3.151E-04 &2.372 &8.246E-06 &2.910 &5.763E-08 &4.239 \\
1600 &7.846E-05 &2.280 &1.031E-06 &2.932 &3.595E-09 &4.180\\ \hline
\end{tabular} \caption{Order of accuracy study on static mesh using Rusanov
flux} \label{tab:ord0} \end{center} \end{table}

\begin{table}
\begin{center} \begin{tabular}{|c|c|c|c|c|c|c|} \hline
  \multirow{2}{*}{$N$}  & \multicolumn{2}{|c|}{$k=1$} &
  \multicolumn{2}{|c|}{$k=2$} & \multicolumn{2}{|c|}{$k=3$} \\ \cline{2-7} &
  Error & Rate & Error & Rate & Error & Rate \\ \hline 100&2.331E-02& -
        &3.979E-03& - &8.633E-04& - \\
  200&6.139E-03&1.9250&4.058E-04&3.294&1.185E-05& 6.186 \\
  400&1.406E-03&2.0258&5.250E-05&3.122&7.079E-07& 5.126 \\
  800&3.375E-04&2.0366&6.626E-06&3.077&4.340E-08& 4.760 \\
1600&8.278E-05&2.0344&8.304E-07&3.057&2.689E-09& 4.573\\ \hline \end{tabular}
\caption{Order of accuracy study on  moving mesh using Rusanov flux}
\label{tab:ord1} \end{center}
\end{table}

\begin{table} \begin{center} \begin{tabular}{|c|c|c|c|c|c|c|} \hline
  \multirow{2}{*}{$N$}  & \multicolumn{2}{|c|}{$k=1$} &
  \multicolumn{2}{|c|}{$k=2$} & \multicolumn{2}{|c|}{$k=3$} \\ \cline{2-7} &
  Error & Rate & Error & Rate & Error & Rate \\ \hline 100&4.582E-02
        &&3.952E-03&&3.464E-04& \\
  200&9.611E-03&2.253&4.048E-04&3.287&2.058E-05& 4.073\\
  400&2.052E-03&2.240&4.640E-05&3.206&1.287E-06& 4.036\\ 800 &4.803E-04  &2.192
     &5.623E-06  &3.152  &8.061E-08  &4.023\\ 1600  &1.184E-04  &2.149
     &6.929E-07  &3.119  &5.050E-09  &4.016\\ \hline \end{tabular}
   \caption{Order of accuracy study on static mesh using HLLC flux}
 \label{tab:hllcord0} \end{center} \end{table}

\begin{table} \begin{center} \begin{tabular}{|c|c|c|c|c|c|c|} \hline
  \multirow{2}{*}{$N$}  & \multicolumn{2}{|c|}{$k=1$} &
  \multicolumn{2}{|c|}{$k=2$} & \multicolumn{2}{|c|}{$k=3$} \\ \cline{2-7} &
  Error & Rate & Error & Rate & Error & Rate \\ \hline 100 &1.590E-02  &&
  1.626E-03 &&  1.962E-04&  \\ 200 &4.042E-03  &1.977  &2.072E-04  &2.972
            &1.269E-05  &3.950\\ 400 &1.014E-03  &1.985  &2.605E-05  &2.982
            &7.983E-07  &3.971\\ 800 &2.538E-04  &1.990  &3.261E-06  &2.988
            &4.997E-08  &3.980\\ 1600  &6.349E-05  &1.992  &4.077E-07  &2.991
            &3.124E-09  &3.985\\ \hline \end{tabular} \caption{Order of accuracy
          study on  moving mesh using HLLC flux} \label{tab:hllcord1}
\end{center} \end{table}

\begin{table}
\begin{center}
\begin{tabular}{|c|c|c|c|c|} \hline
  \multirow{2}{*}{$N$}  & \multicolumn{2}{|c|}{$k=1$} &
  \multicolumn{2}{|c|}{$k=2$}  \\ \cline{2-5} &
  Error & Rate & Error & Rate  \\
  \hline
  100 &2.053E-02& - &2.277E-03 & - \\
  200 &4.312E-03 &2.251 &3.425E-04 &  2.732   \\
  400 &1.031E-03 &2.064 &4.565E-05 &  2.907   \\
  800 &2.550E-04 &2.015 &5.812E-06 &  2.973   \\
  1600 &6.356E-05 &2.004 &7.315E-07&  2.990   \\
  \hline
\end{tabular}
\caption{Order of accuracy study on moving mesh using Rusanov flux using
Higher Order Limiter \cite{Zhong:2013:SWE:2397205.2397446}}
\label{tab:limitedord0}
\end{center}
\end{table}

\begin{table} \begin{center} \begin{tabular}{|c|c|c|c|c|c|c|} \hline
  \multirow{2}{*}{$N$}  & \multicolumn{2}{|c|}{$k=1$} &
  \multicolumn{2}{|c|}{$k=2$} & \multicolumn{2}{|c|}{$k=3$} \\ \cline{2-7} &
  Error & Rate & Error & Rate & Error & Rate \\ \hline 100&  1.735E-02  &   -
        &1.798E-03 &    -  &2.351E-04 &  -    \\ 200&  4.179E-03  &2.051
        &2.848E-04 &2.676  &1.416E-05 &4.069  \\ 400&  1.054E-03  &2.035
        &4.301E-05 &2.703  &8.578E-07 &4.041  \\ 800&  2.615E-04  &1.943
        &6.012E-06 &2.838  &5.476E-08 &3.958  \\ 1600& 7.279E-05  &1.852
        &8.000E-07 &2.909  &3.505E-09 &3.966  \\ \hline \end{tabular}
\end{center} \caption{Order of accuracy study on  moving mesh using HLLC flux
with randomly perturbed mesh velocity} \label{tab:hllcord2} \end{table}
\color{blue}
\subsection{Smooth Test Case with Non-Constant Velocity}
We also test the accuracy of our schemes on a isentropic problem with smooth
solutions. The test case The initial conditions are given by
\begin{align}
  \rho(x,0) &= 1 + 0.9999995 \sin(\pi x) & u(x,0) &= 0 & p(x,0) &=
  \rho^\gamma(x,0)
\end{align}
with $\gamma =3$ and periodic boundary conditions. For this kind of special
isentropic problem, the Euler equations are equivalent to the two Burgers
equations in terms of their two Riemann invariants which can then be used to
derive the analytical solution. The errors are then computed with respect to the
given analytical solution. In contrast to the previous test case, the velocity
and pressure are not constant which makes this a more challenging test case. We
run the simulation with a WENO-type limiter from
\cite{Zhang:2010:PHO:1864819.1865066} and positivity limiter enabled.  As we can
see from Tables~\ref{tab:errorsmoothfixed}, \ref{tab:errorsmooth},
the rate of convergence is maintained for the moving mesh method with the moving
mesh methods exhibiting much lower errors.

\begin{table} \begin{center} \begin{tabular}{|c|c|c|c|c|} \hline
  \multirow{2}{*}{$N$}  & \multicolumn{2}{|c|}{$k=1$} &
  \multicolumn{2}{|c|}{$k=2$} \\ \cline{2-5} &
  Error & Rate & Error & Rate \\
  \hline 100&  8.535E-03  &   - &1.033E-03 &    -  \\
  200& 1.958E-03  &2.124 &1.221E-04 &3.08  \\
  400& 4.721E-04  &2.052 &1.581E-05 &2.95  \\
  800& 1.238E-04  &1.931 &2.14E-06 &2.89  \\
  1600& 3.563E-05  &1.796 &2.63E-07 &3.02 \\ \hline
\end{tabular}
\end{center} \caption{Order of accuracy study on  fixed mesh using Roe flux
with Non-Constant Velocity Smooth Test Case} \label{tab:errorsmoothfixed} \end{table}

\begin{table} \begin{center} \begin{tabular}{|c|c|c|c|c|} \hline
  \multirow{2}{*}{$N$}  & \multicolumn{2}{|c|}{$k=1$} &
  \multicolumn{2}{|c|}{$k=2$} \\ \cline{2-5} &
  Error & Rate & Error & Rate \\
  \hline 100&  4.235E-03  &   - &2.238E-04 &    -  \\
  200& 1.058E-03  &2.001 &3.255E-05 &2.87  \\
  400& 2.586E-04  &2.035 &4.301E-05 &3.133  \\
  800& 5.804E-05  &2.155 &5.762E-06 &2.901  \\
  1600& 1.271E-05  &2.192 &7.401E-07 &2.96 \\ \hline
\end{tabular}
\end{center} \caption{Order of accuracy study on  moving mesh using Roe flux
with Non-Constant Velocity Smooth Test Case} \label{tab:errorsmooth} \end{table}

\color{black}
\subsection{Single contact wave} 
\begin{figure} \begin{center}
    \begin{tabular}{cc}
      \includegraphics[width=0.45\textwidth]{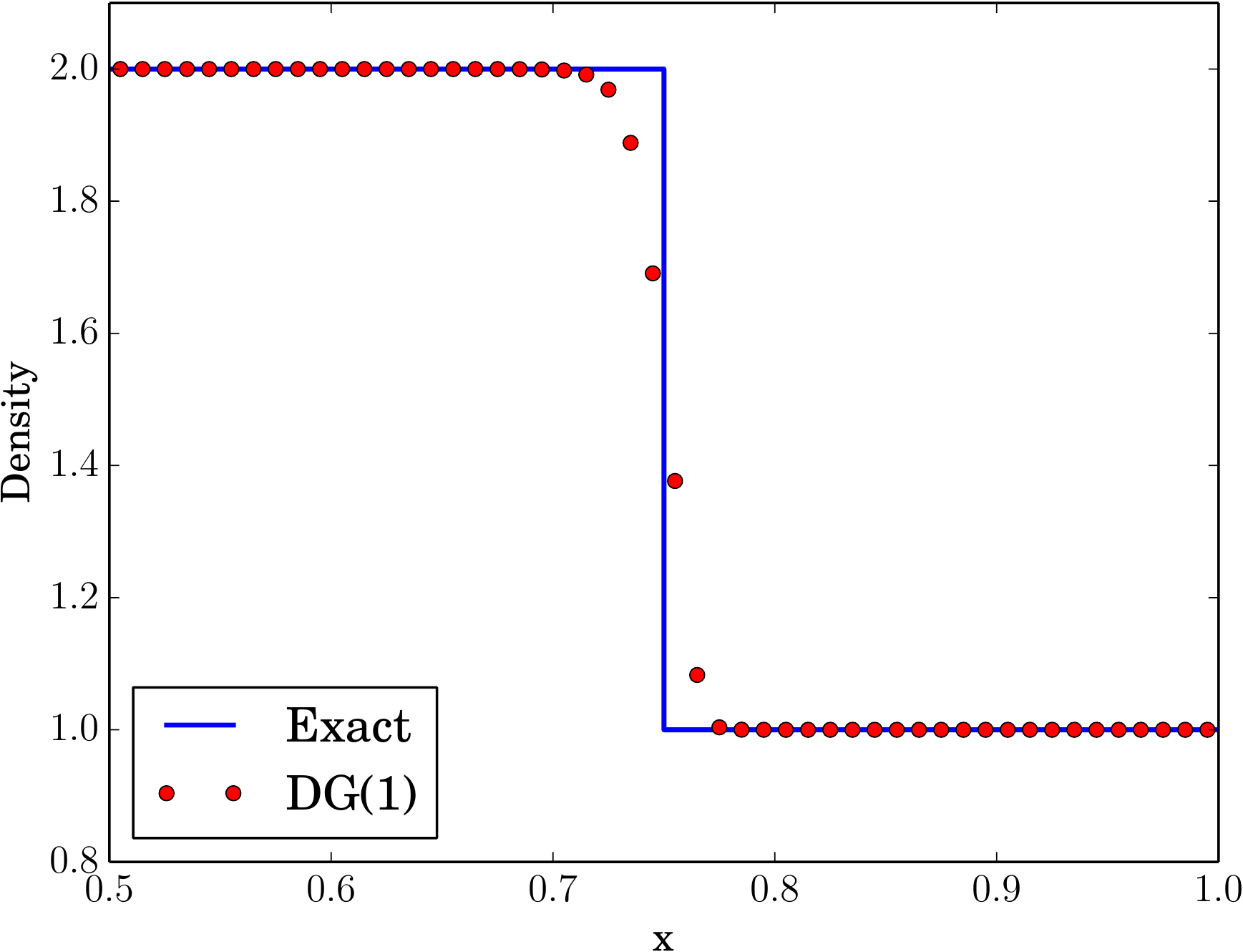} &
      \includegraphics[width=0.45\textwidth]{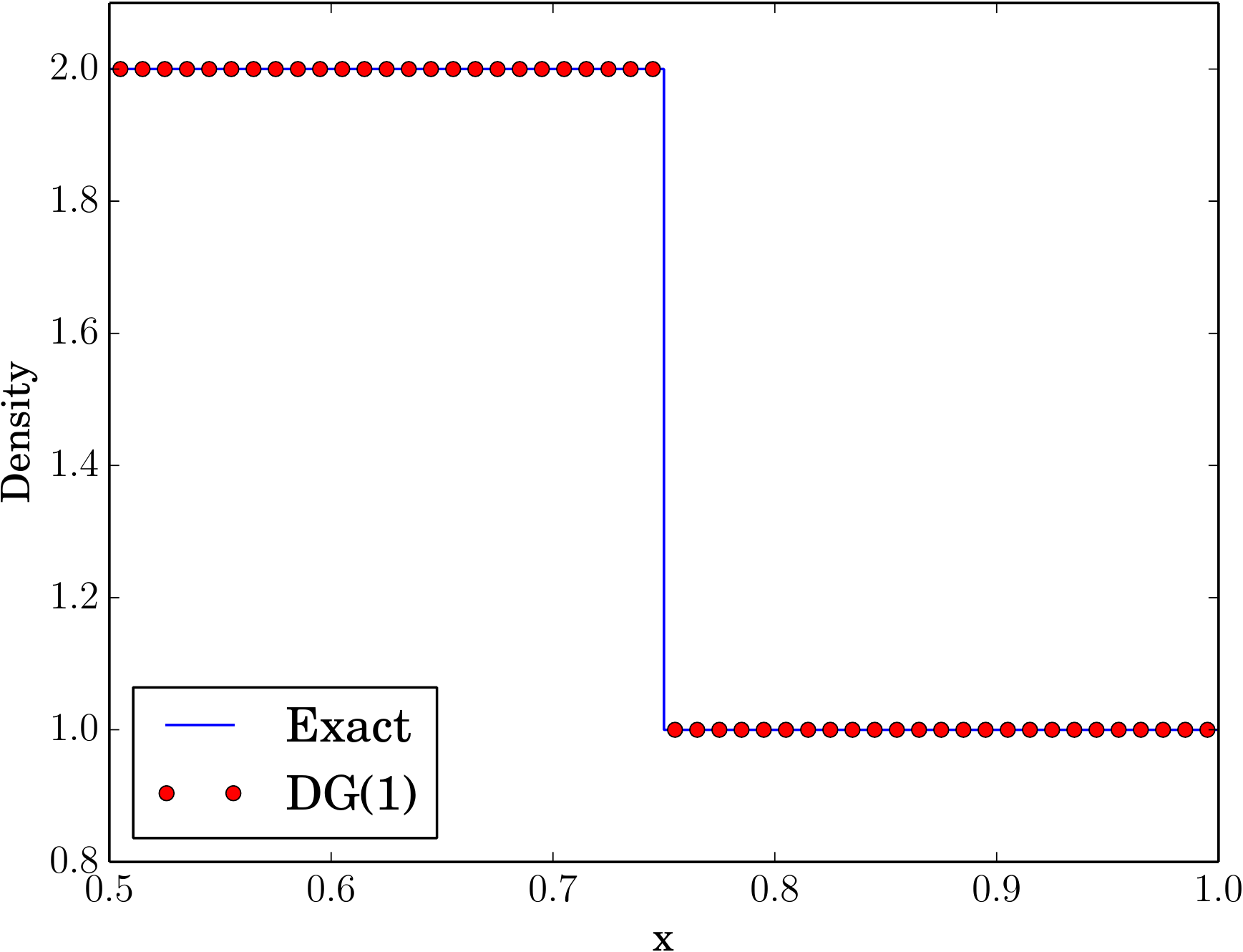} \\ (a) &
      (b) \end{tabular} \end{center} \caption{Single contact wave using Roe flux
    and 100 cells: (a) static mesh, (b) moving mesh} \label{fig:contact}
    \end{figure}

In this example, we choose a Riemann problem
  which gives rise to a single contact wave in the solution that propagates with
  a constant speed. The initial condition is given by \[ (\rho,v,p) =
    \begin{cases} (2.0, 1.0, 1.0) & \textrm{ if } x < 0.5 \\ (1.0, 1.0, 1.0) &
  \textrm{ if } x > 0.5 \end{cases} \] and the contact wave moves with a
  constant speed of 1.0. The solution on static and moving meshes are shown in
  figure~(\ref{fig:contact}) at time $t=0.5$ using Roe flux. The moving mesh is
  able to exactly resolve the contact wave while the static mesh scheme adds
  considerable numerical dissipation that smears the discontinuity over many
  cells. The accurate resolution of contact waves is a key advantage of such
  moving mesh methods, which are capable of giving very good resolution of the
  contact discontinuity even on coarse meshes.  
\subsection{Sod problem} 
\begin{figure}
  \begin{center} \begin{tabular}{cc}
    \includegraphics[width=0.48\textwidth]{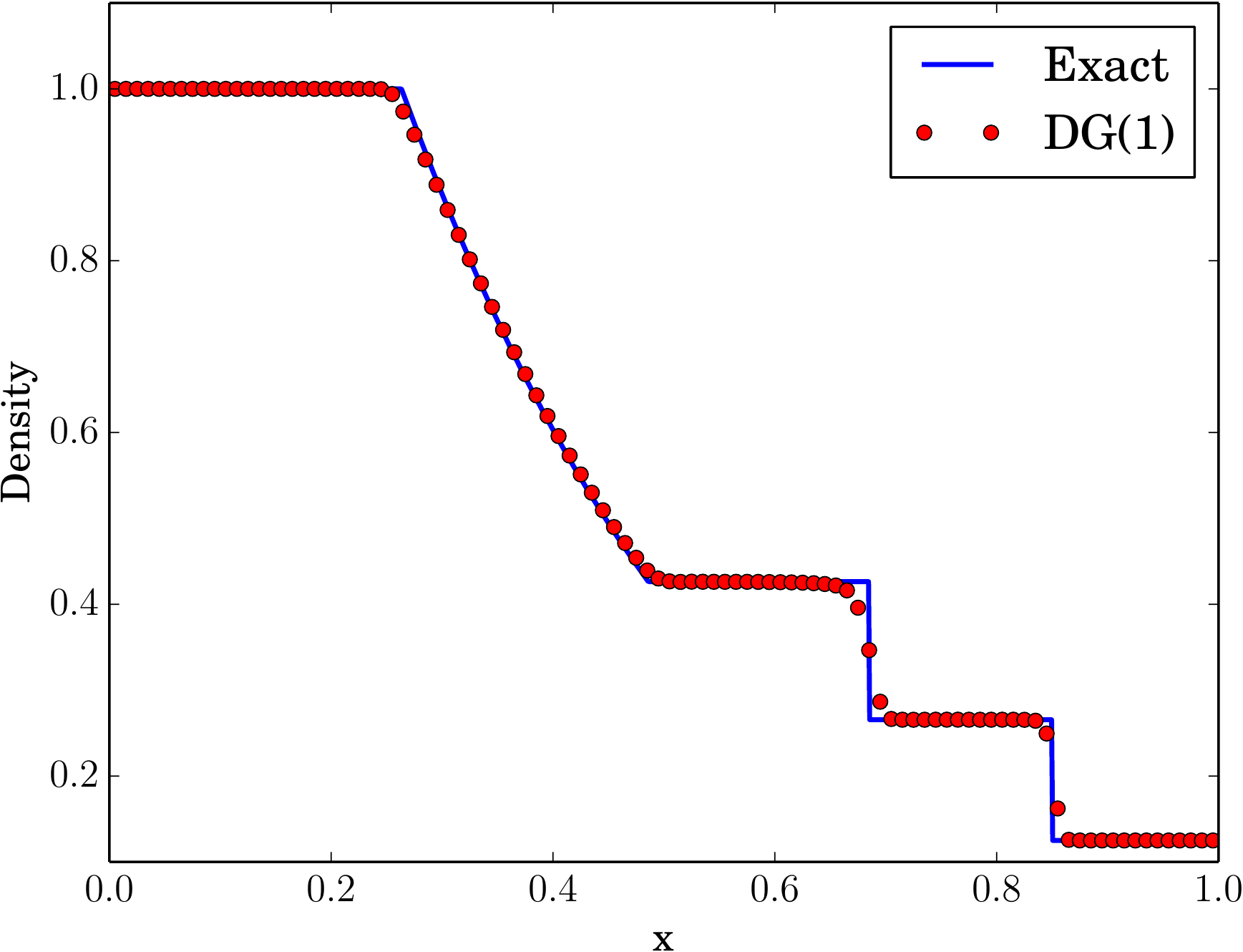} &
    \includegraphics[width=0.48\textwidth]{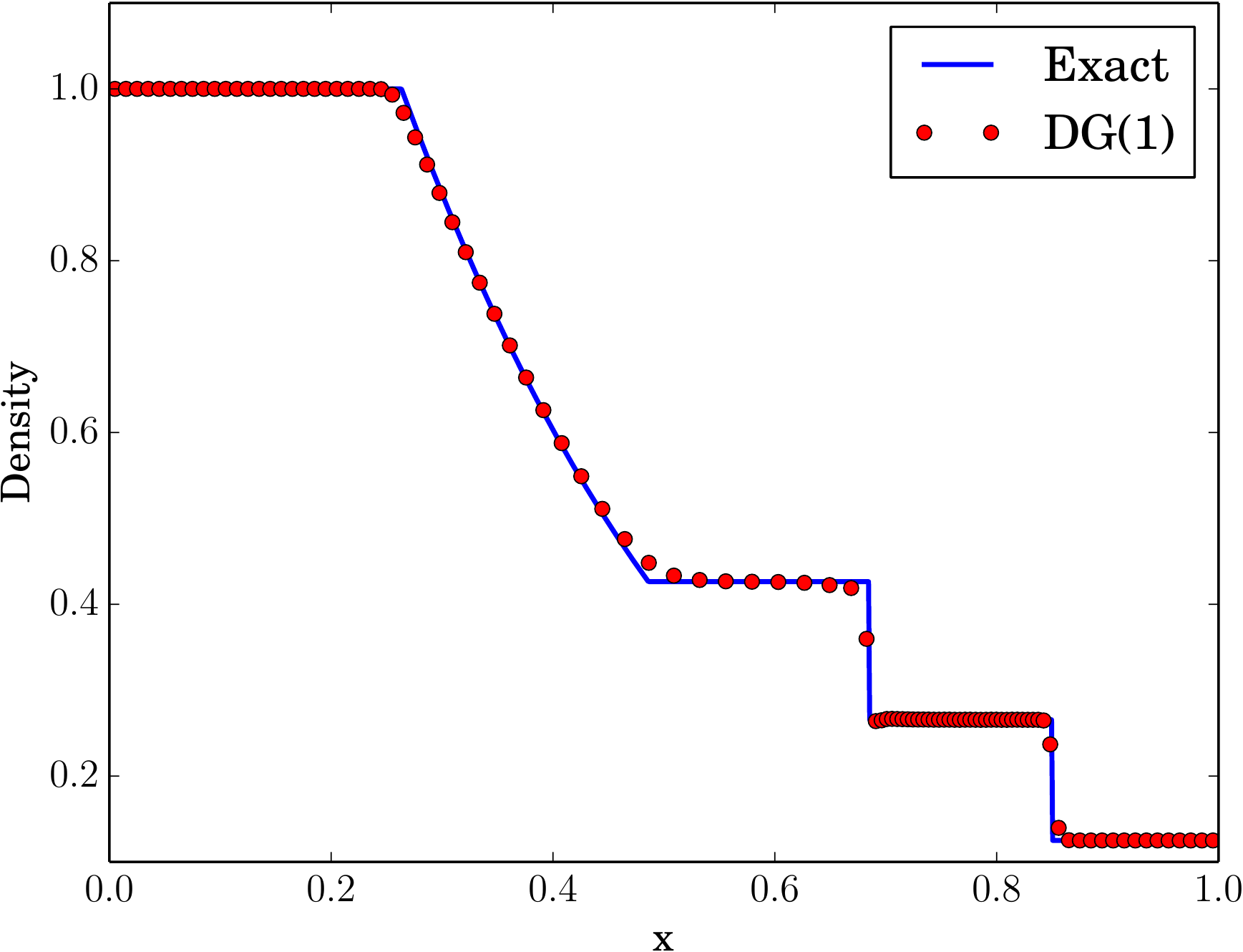} \\ (a) & (b)
  \end{tabular} \end{center} \caption{Sod problem using Roe flux, 100 cells and
TVD limiter: (a) static mesh (b) moving mesh} \label{fig:sod1} \end{figure}

The initial condition for the Sod test case is given
by~\cite{sod1978} \[ (\rho,v,p) = \begin{cases} (1.0, 0.0, 1.0) & \textrm{ if }
x < 0.5 \\ (0.125, 0.0, 0.1) & \textrm{ if } x > 0.5 \end{cases} \] and the
solution is computed upto a final time of $T=0.2$ with the computational domain
being $[0,1]$.  Since the fluid velocity is zero at the boundary, the
computational domain does not change with time for the chosen final time. The
exact solution consists of a rarefaction fan, a contact wave and a shock wave.
In figure~(\ref{fig:sod1}), we show the results obtained using Roe flux with 100
cells and TVD limiter on static and moving mesh. The contact wave is
considerably well resolved on the moving mesh as compared to the static mesh due
to reduced numerical dissipation on moving meshes.  

\begin{figure} \begin{center} \begin{tabular}{cc}
\includegraphics[width=0.48\textwidth]{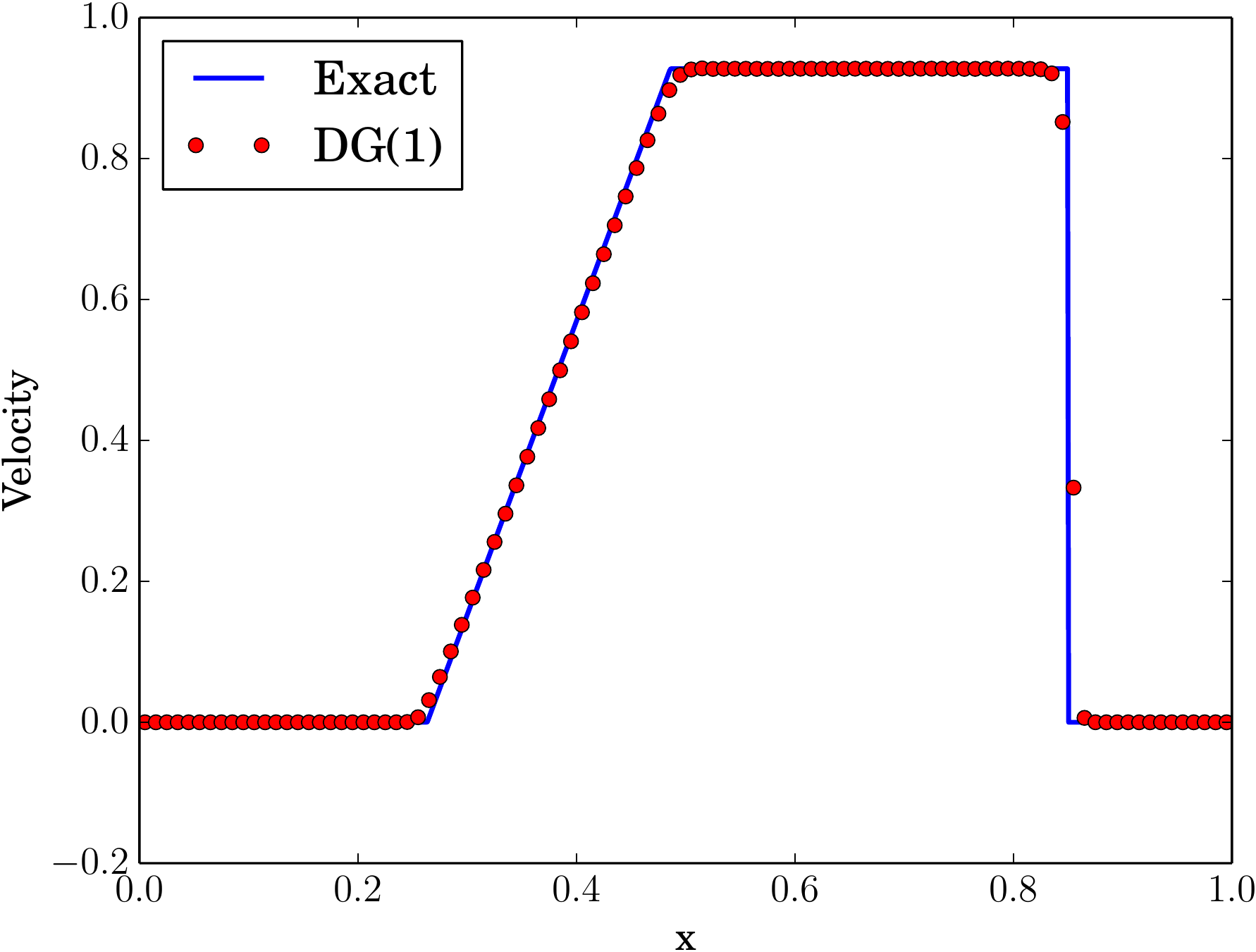} &
\includegraphics[width=0.48\textwidth]{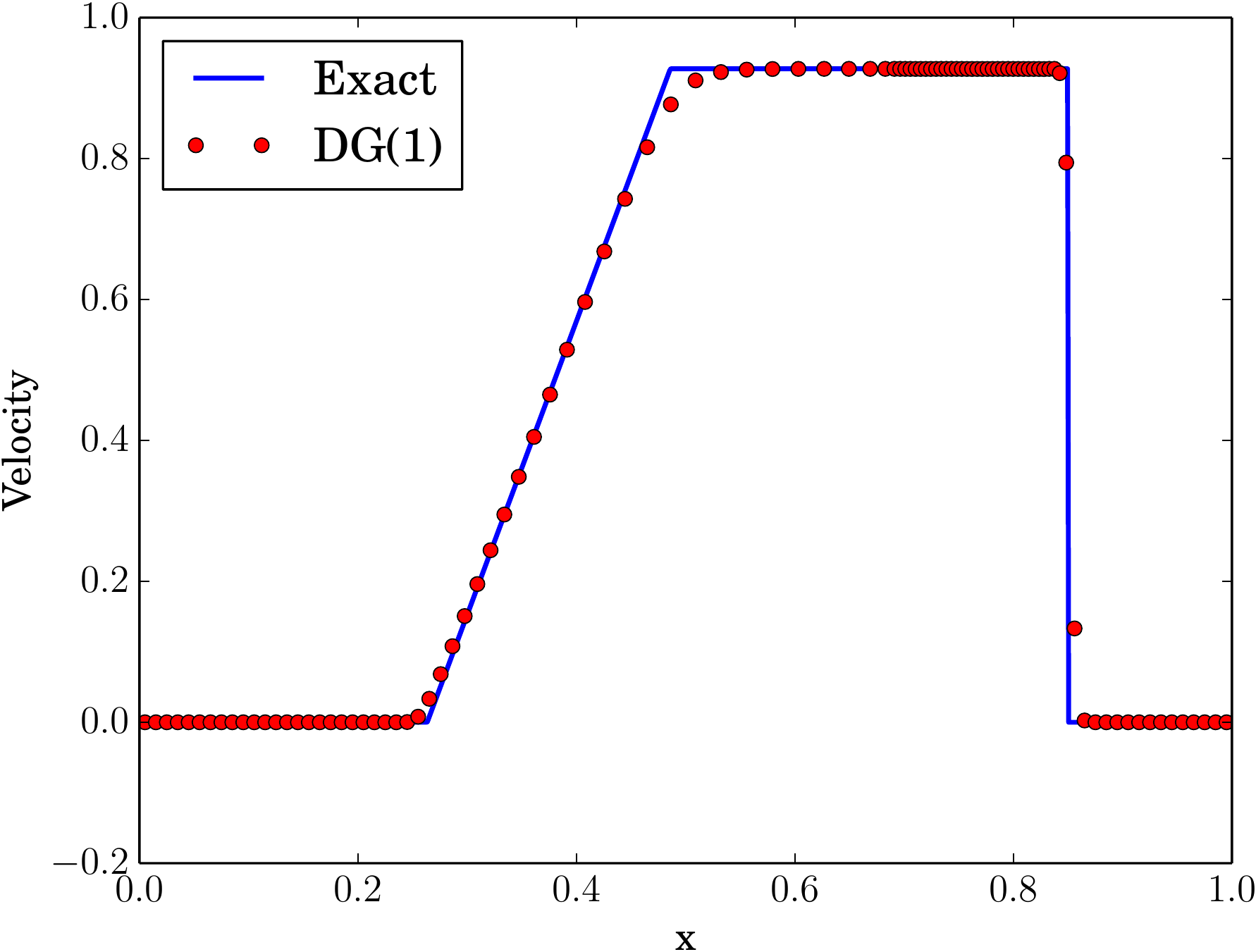} \\ (a) & (b)
\end{tabular} \end{center} \caption{Sod problem using Roe flux, 100 cells and
TVD limiter: (a) static mesh (b) moving mesh} \label{fig:sod2} \end{figure} To
study the Galilean invariance or the dependence of the solution on the choice of
coordinate frame, we add a boost velocity of $V=10$ or $V=100$ to the coordinate
frame, while implies the initial fluid velocity is $v(x,0)=V$ and the other
quantities remain as before. Figure~(\ref{fig:sod3}a) shows that the accuracy of
the static mesh results degrades with increase in velocity of the coordinate
frame, particularly the contact discontinuity is highly smeared. The results
given in figure~(\ref{fig:sod3}b) clearly show the independence of the results
on the moving mesh with respect to the coordinate frame velocity.  The allowed
time step from CFL condition decreases with increase in coordinate frame speed
for the static mesh case, while in case of moving mesh, it remains invariant.
This means that in case of static mesh, we have to perform more time steps to
reach the same final time as shown in table~(\ref{tab:sodboost}), which
increases the computational time. Thus the moving mesh scheme has the additional
advantage of allowing a larger time step compared to the fixed mesh scheme.

\begin{figure} \centering
\includegraphics[width=0.6\textwidth]{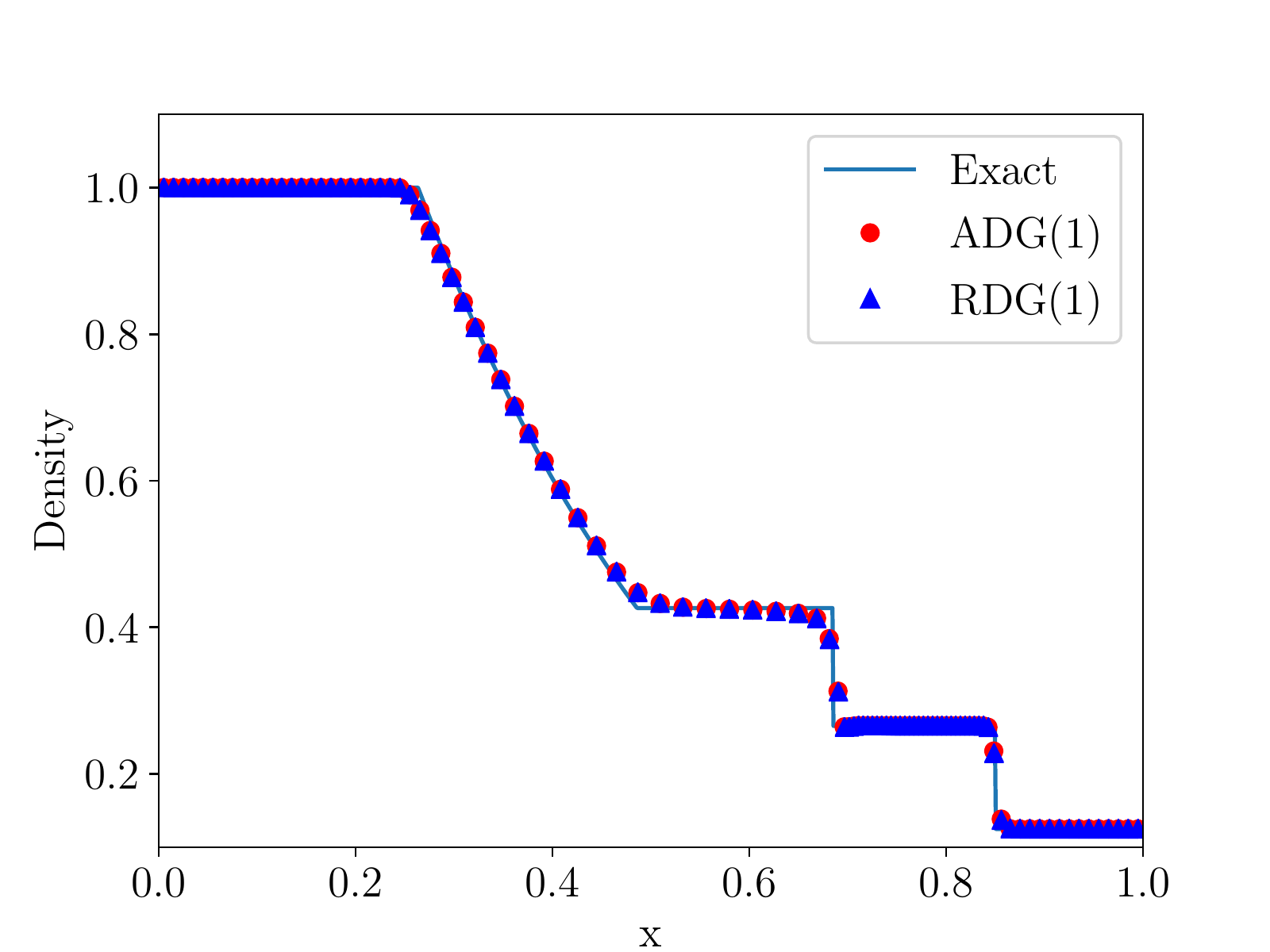}
\caption{Sod problem using Roe flux, 100 cells and TVD limiter. ADG : Average
Velocity, RDG : Linearized Riemann Velocity} \end{figure}

\begin{table} \begin{center} \begin{tabular}{|c|c|c|c|} \hline $V$ & $0$ & $10$
& $100$ \\ \hline static mesh & 144 & 810  & 6807 \\ moving mesh & 176 & 176 &
176 \\ \hline \end{tabular} \end{center} \caption{Number of iterations required
to reach time $t=0.2$ for Sod test for different boost velocity of the
coordinate frame} \label{tab:sodboost} \end{table}

\begin{figure} \begin{center} \begin{tabular}{cc}
\includegraphics[width=0.48\textwidth]{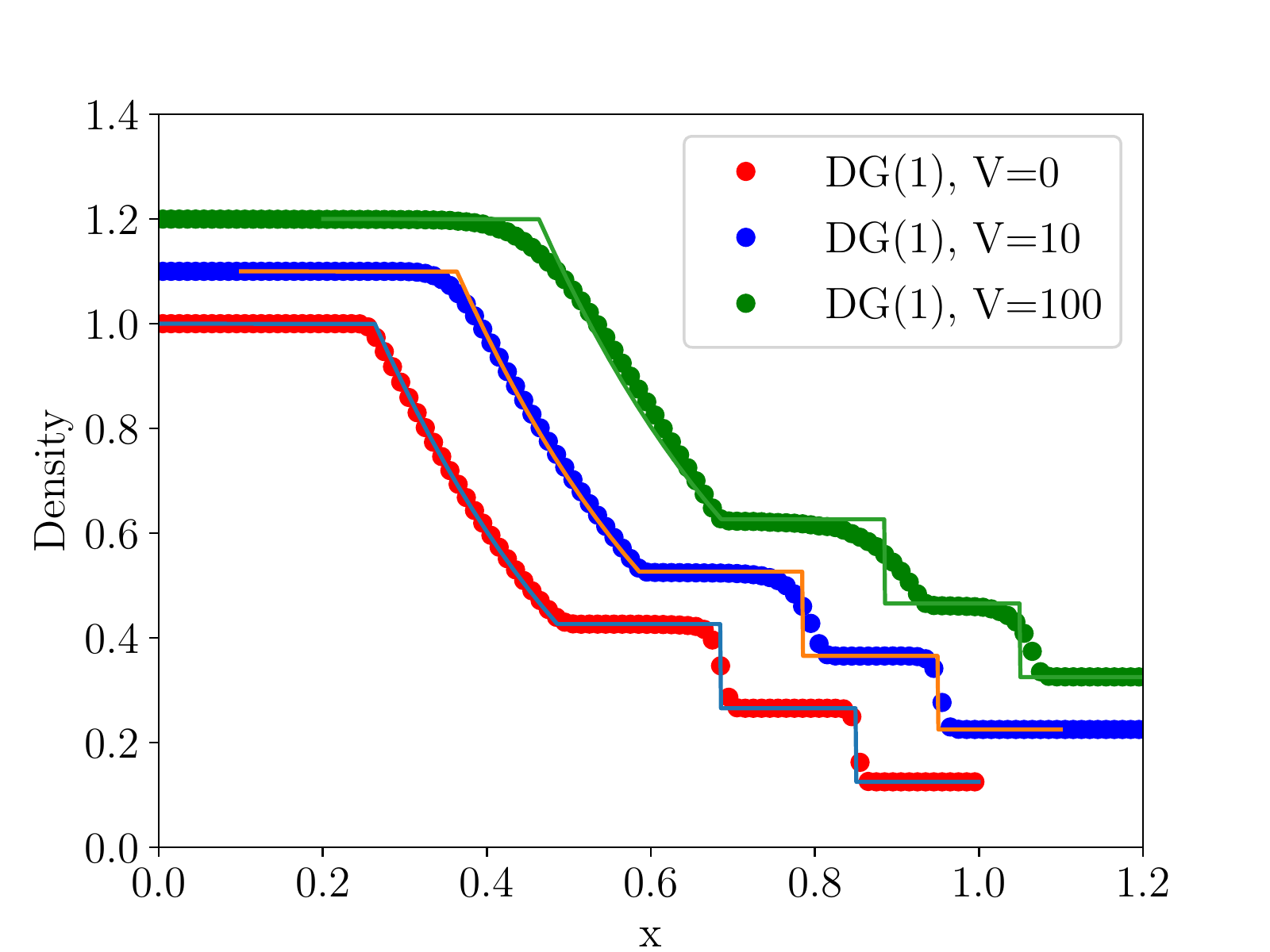} &
\includegraphics[width=0.48\textwidth]{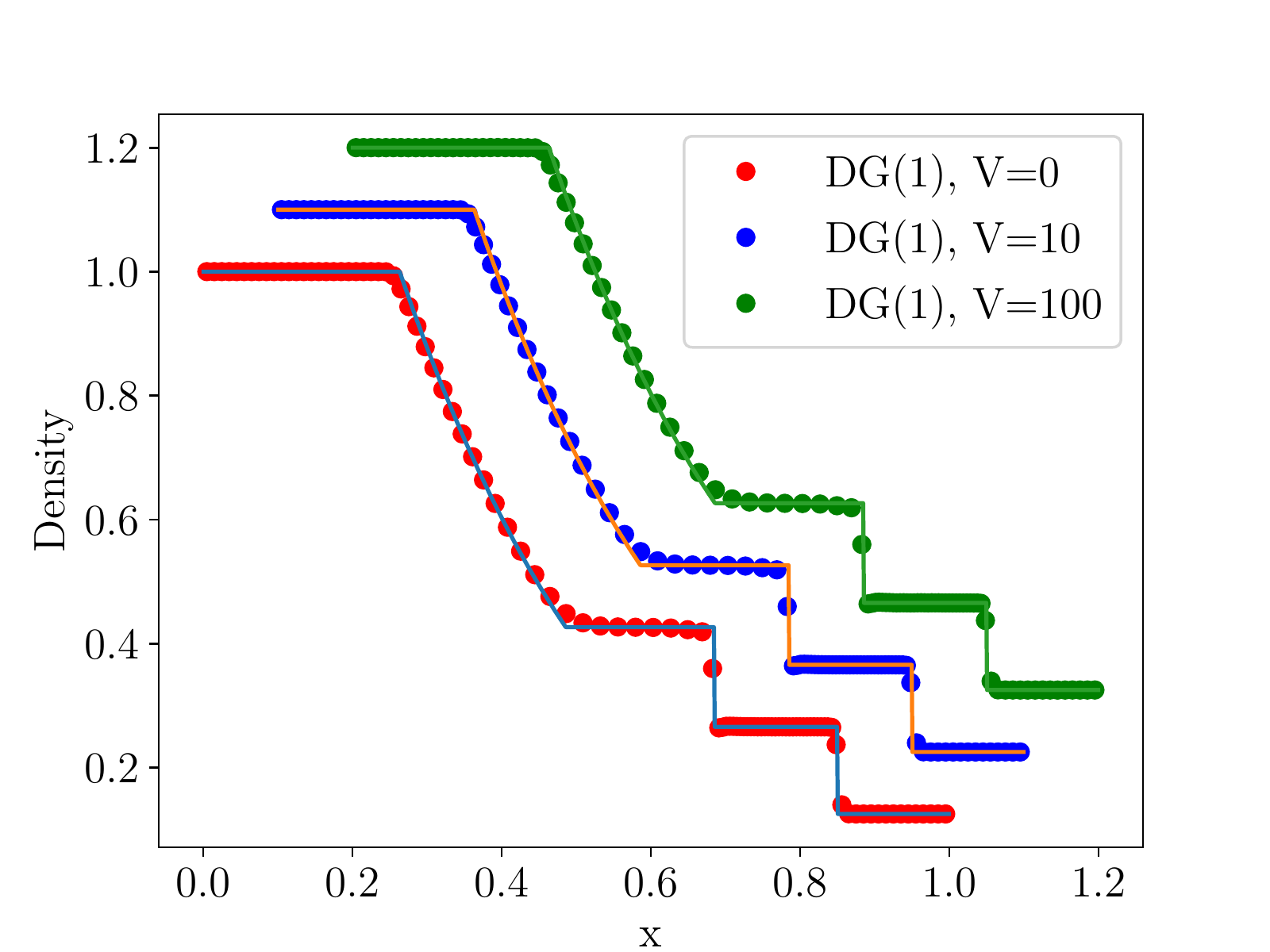} \\ (a) & (b)
\end{tabular} \end{center} \caption{Effect of coordinate frame motion on Sod
problem using Roe flux, 100 cells and TVD limiter: (a) static mesh (b) moving
mesh} \label{fig:sod3} \end{figure}
\begin{figure} \begin{center} \begin{tabular}{cc}
    \includegraphics[width=0.48\textwidth]{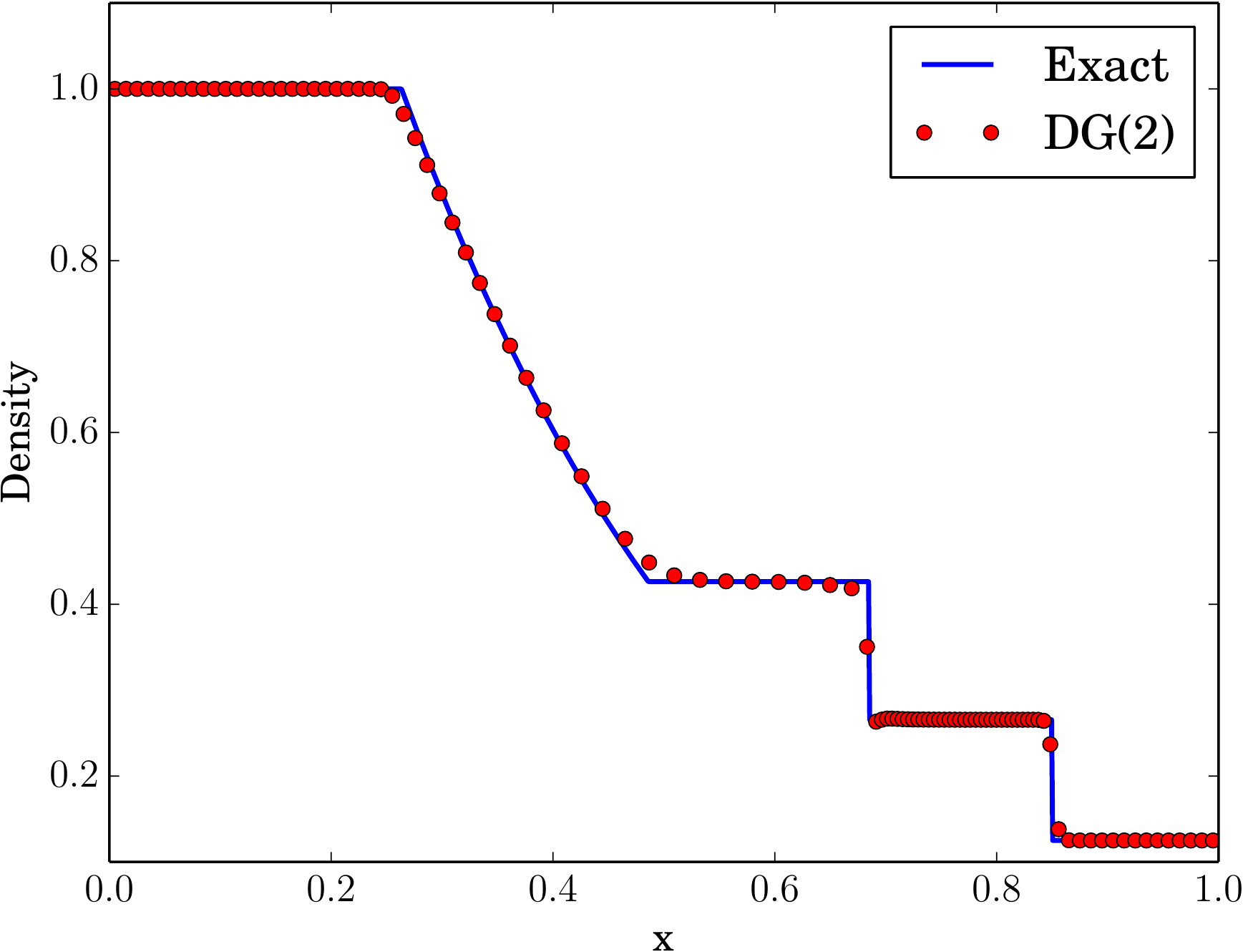} &
    \includegraphics[width=0.48\textwidth]{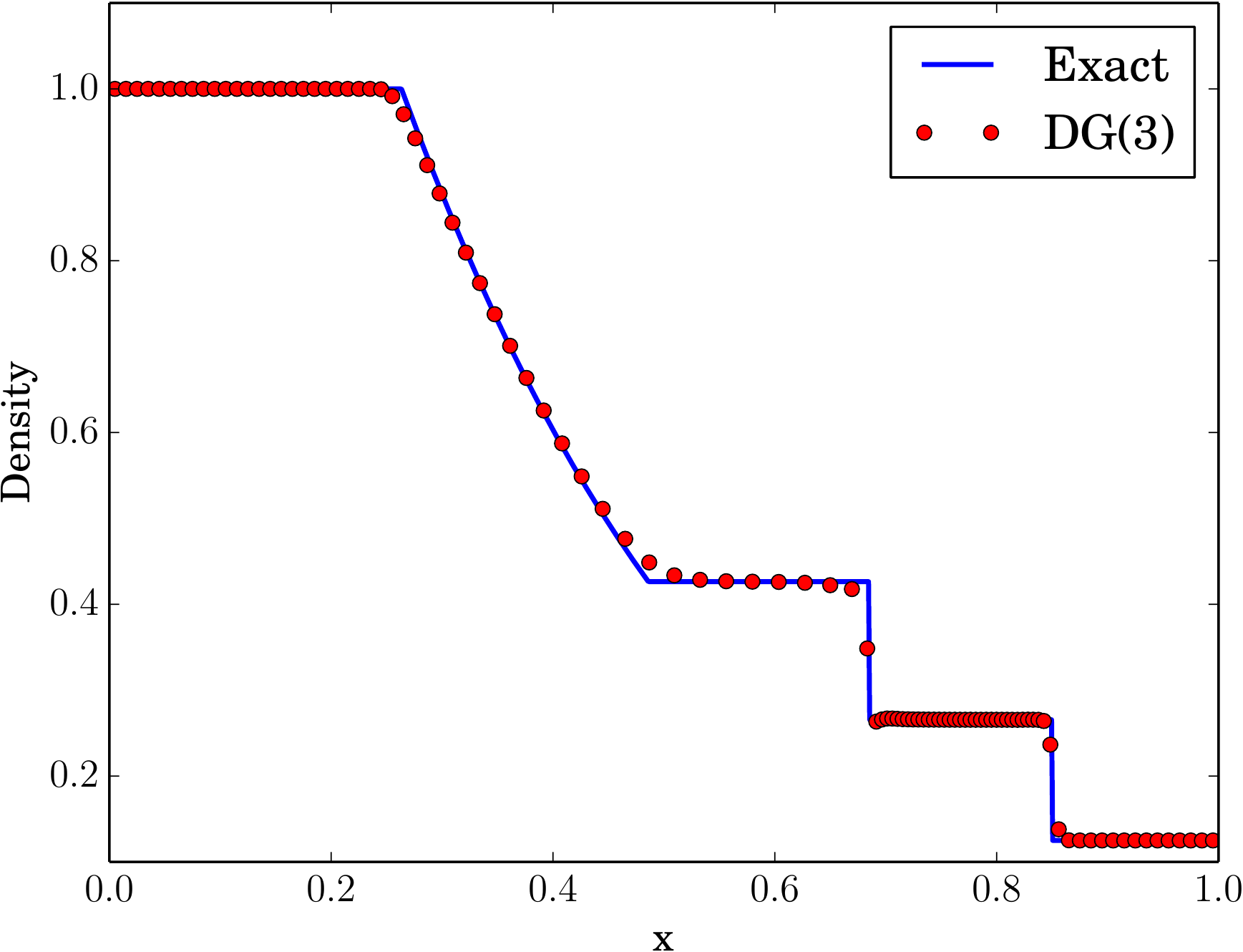} \\ (a)
                                                                         & (b)
  \end{tabular} \end{center} \caption{Sod problem on moving mesh using Roe flux,
100 cells and TVD limiter: (a) Degree = 2 (b) Degree = 3} \label{fig:sod4}
\end{figure}

 Finally, we compute the solutions using
quadratic and cubic polynomials and the results are shown in
figure~(\ref{fig:sod4}). The solutions look similar to the case of linear
polynomials and have the same sharp resolution of discontinuities.
\subsection{Lax problem} 
\begin{figure} \begin{center} \begin{tabular}{cc}
\includegraphics[width=0.48\textwidth]{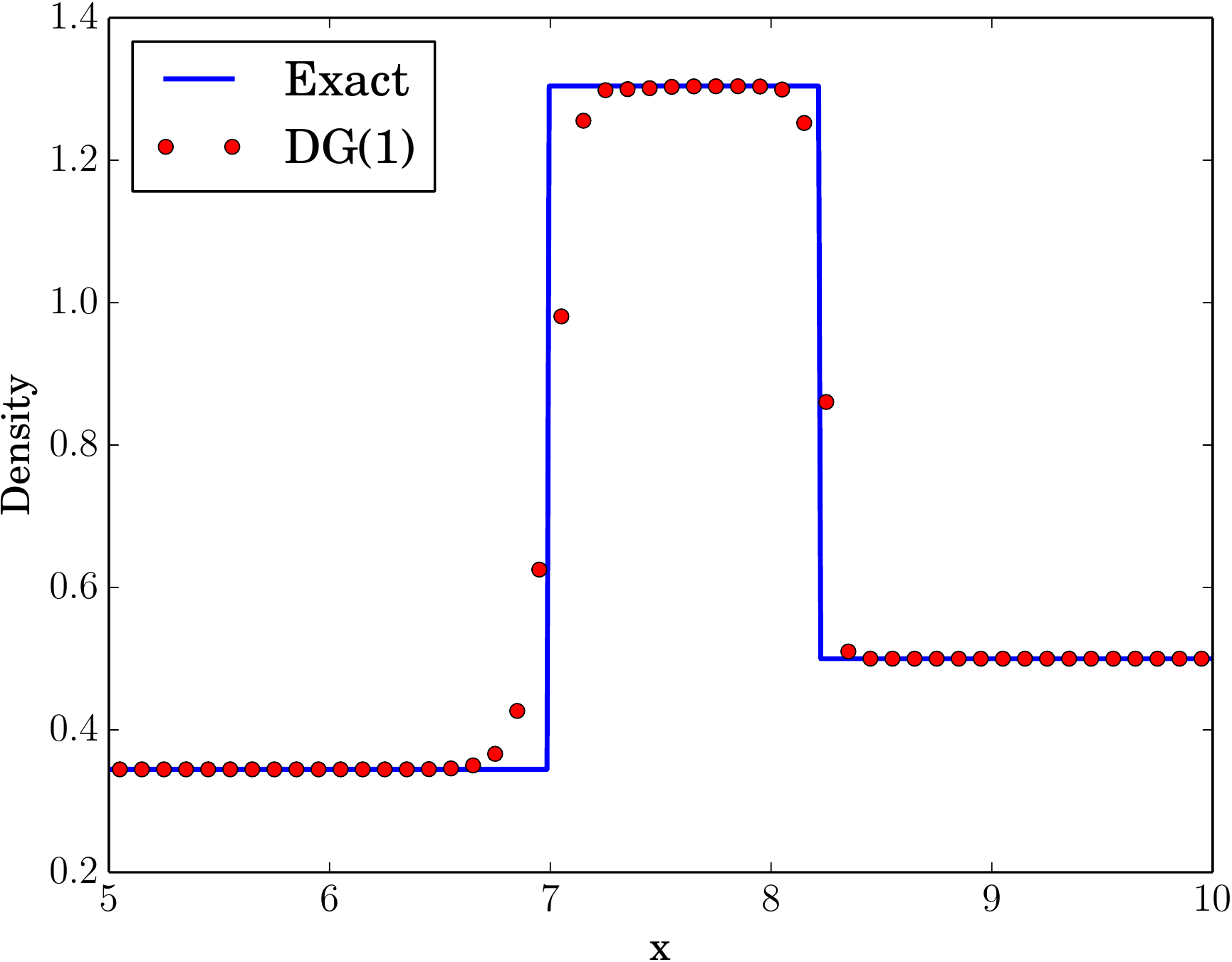} &
\includegraphics[width=0.48\textwidth]{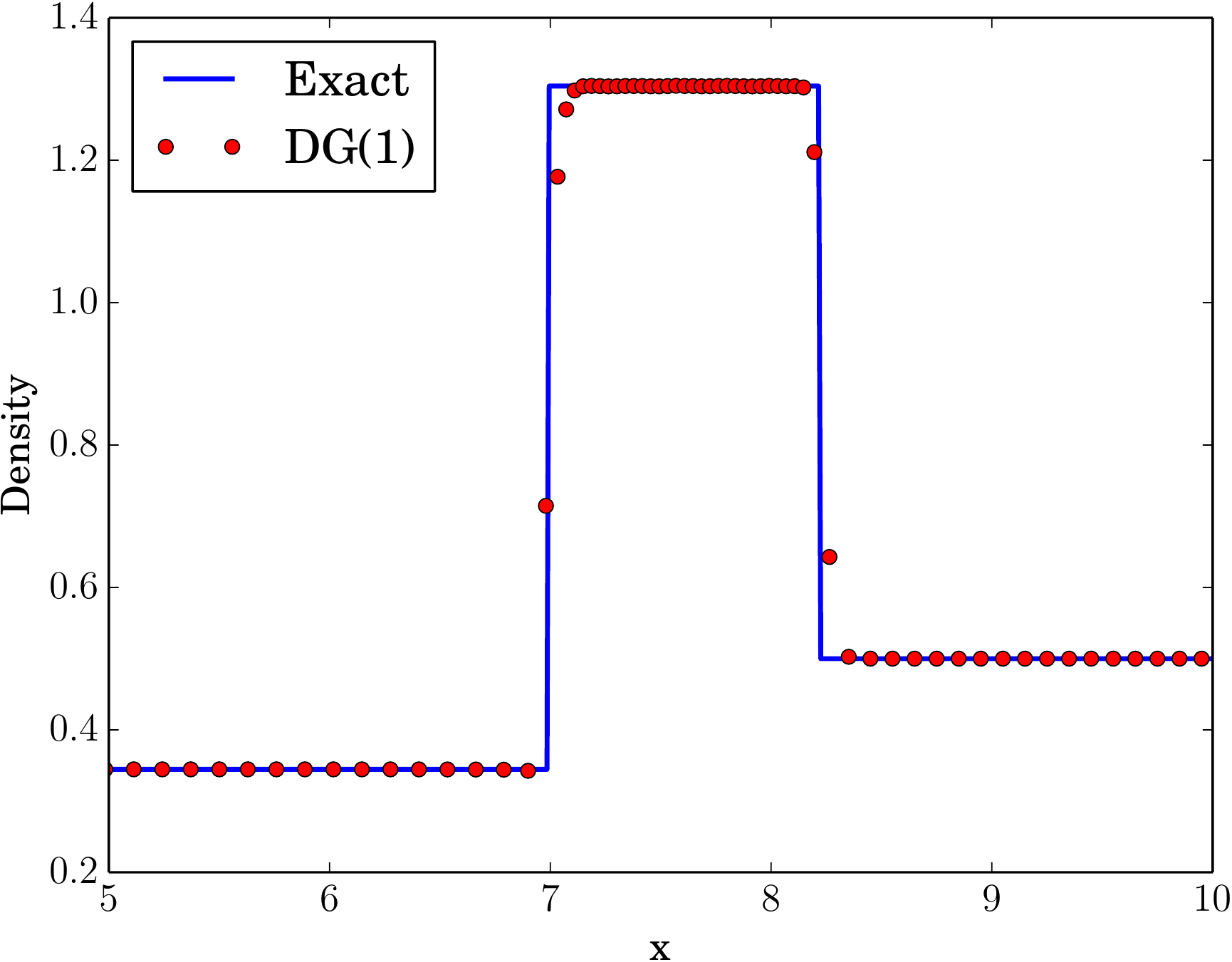} \\ (a) & (b)
\end{tabular} \end{center} \caption{Lax problem using HLLC flux, 100 cells and
TVD limiter: (a) static mesh (b) moving mesh} \label{fig:lax1} \end{figure}

The initial condition is given by \[ (\rho,v,p) =
\begin{cases} (0.445, 0.698, 3.528) & \textrm{if } x < 0 \\ (0.5, 0.0, 0.571) &
\textrm{if } x > 0 \end{cases} \] The computational domain is $[-10,+10]$ and we
compute the solution up to a final time of $T=1.3$. This problem has a strong
shock and a contact wave that is difficult to resolve accurately. The zoomed
view of density is shown at the final time in figure~(\ref{fig:lax1}), and we
observe the moving mesh results are more accurate for the contact wave, which is
the first discontinuity in the figure. The second discontinuity is a shock which
is equally well resolved in both cases. We can observe that the grid is
automatically clustered in the region between the contact and shock wave, but no
explicit grid adaptation was used in this simulation.

\begin{figure} \centering
\includegraphics[width=0.6\textwidth]{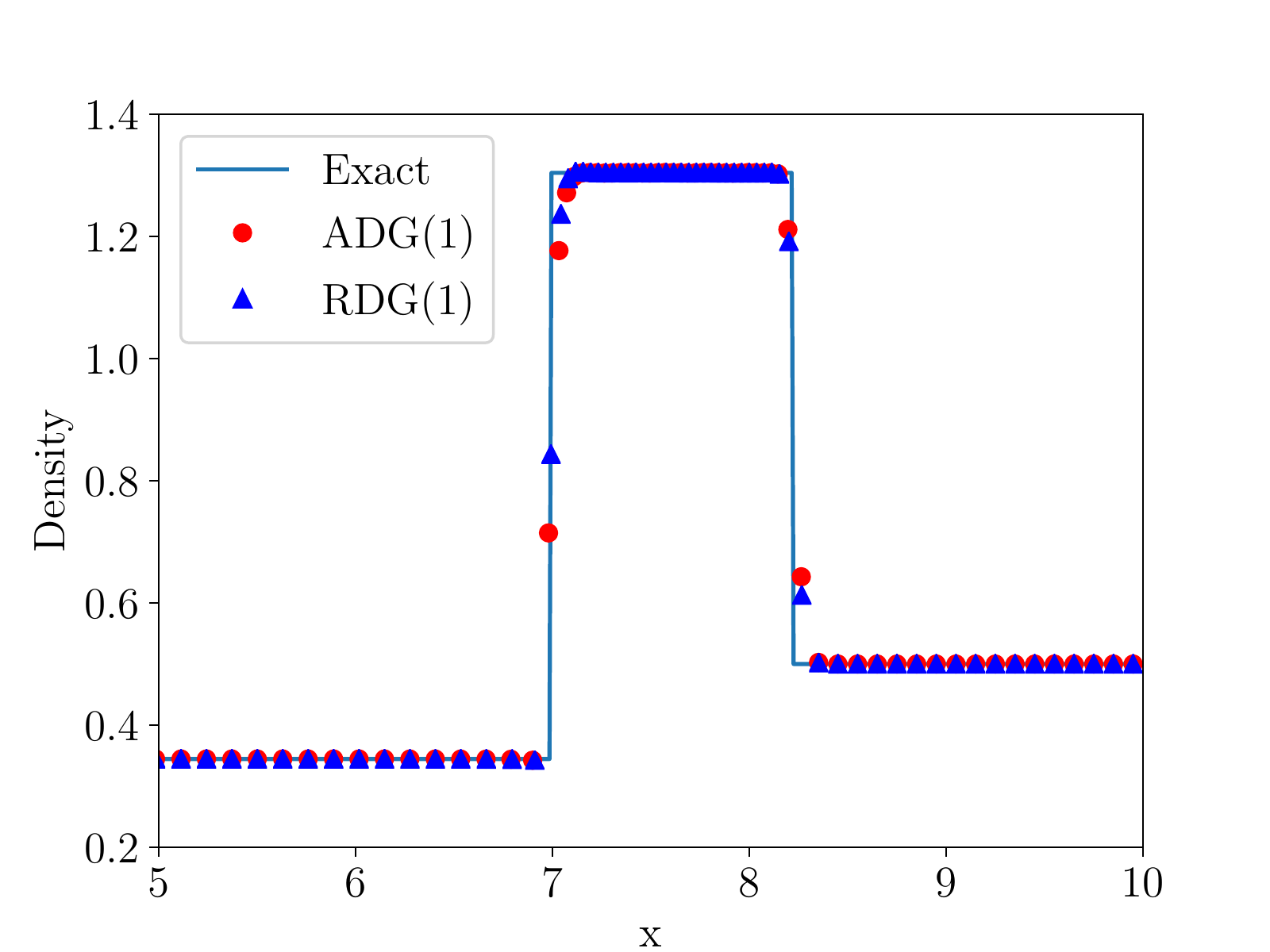}
\caption{Lax problem using HLLC flux, 100 cells and TVD limiter. ADG : Average
Velocity, RDG : Linearized Riemann Velocity} \end{figure}
\subsection{Shu-Osher problem} 
\begin{figure} \begin{center} \begin{tabular}{cc}
  \includegraphics[width=0.48\textwidth]{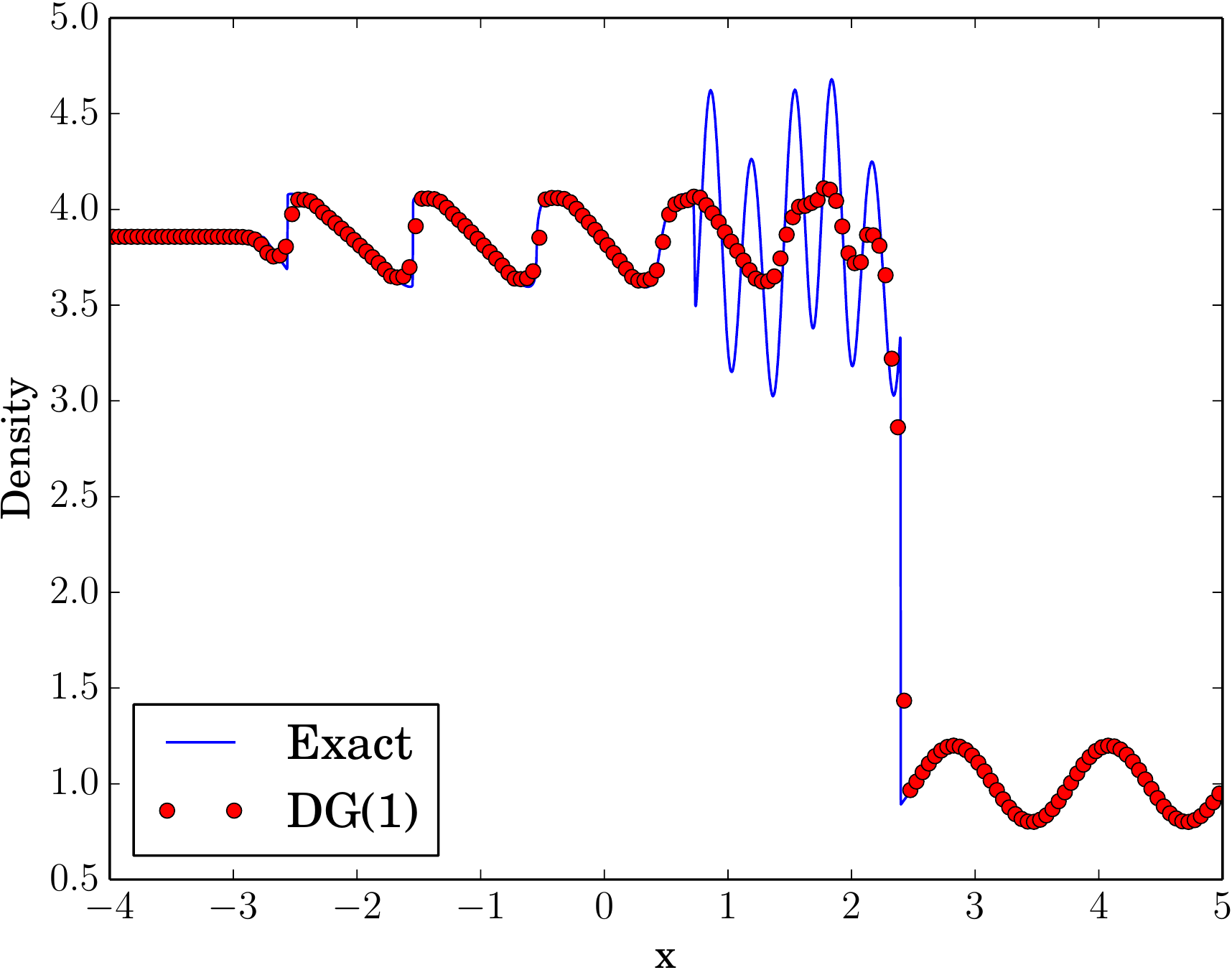} &
  \includegraphics[width=0.48\textwidth]{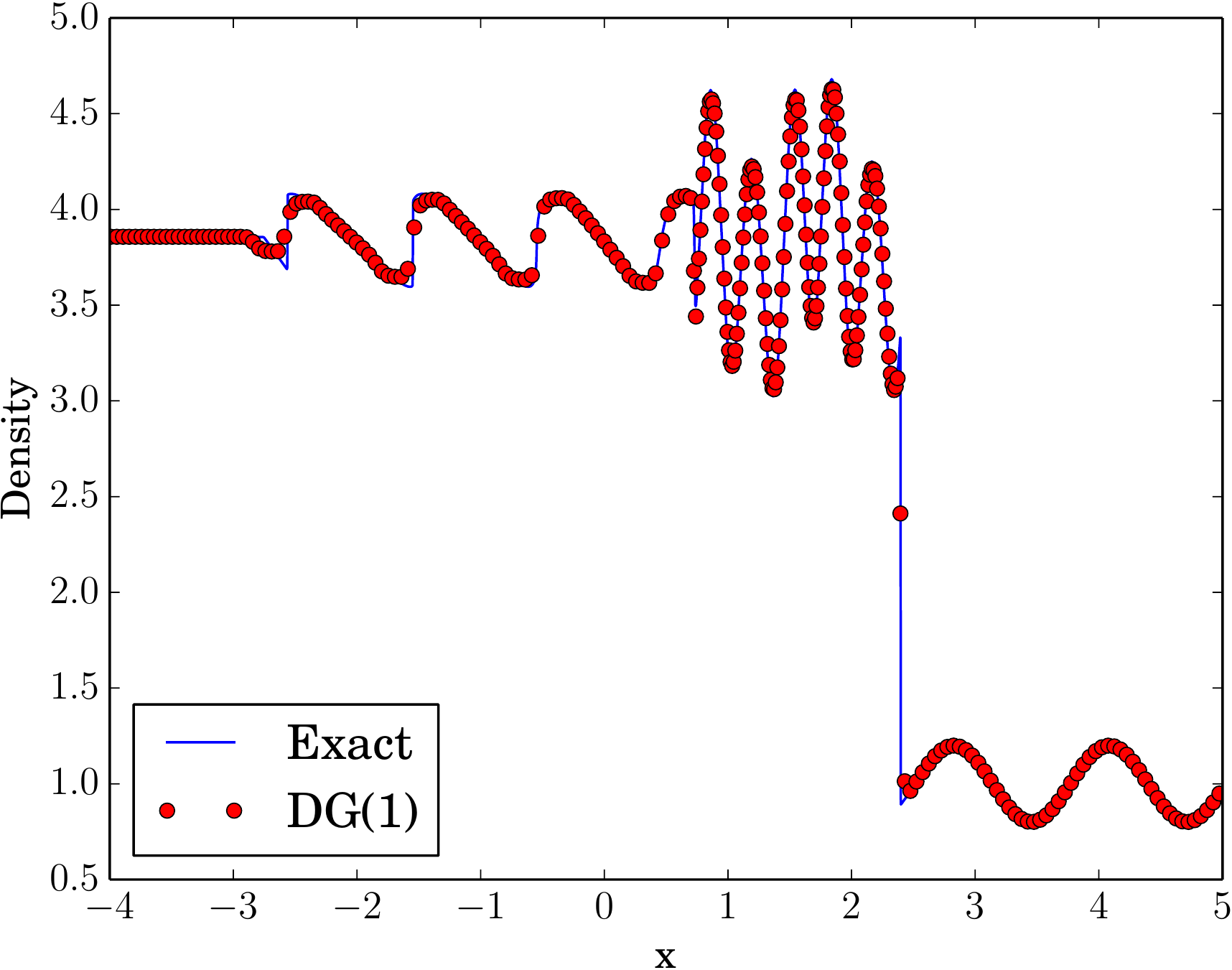} \\ (a)   &
  (b)   \\ \includegraphics[width=0.48\textwidth]{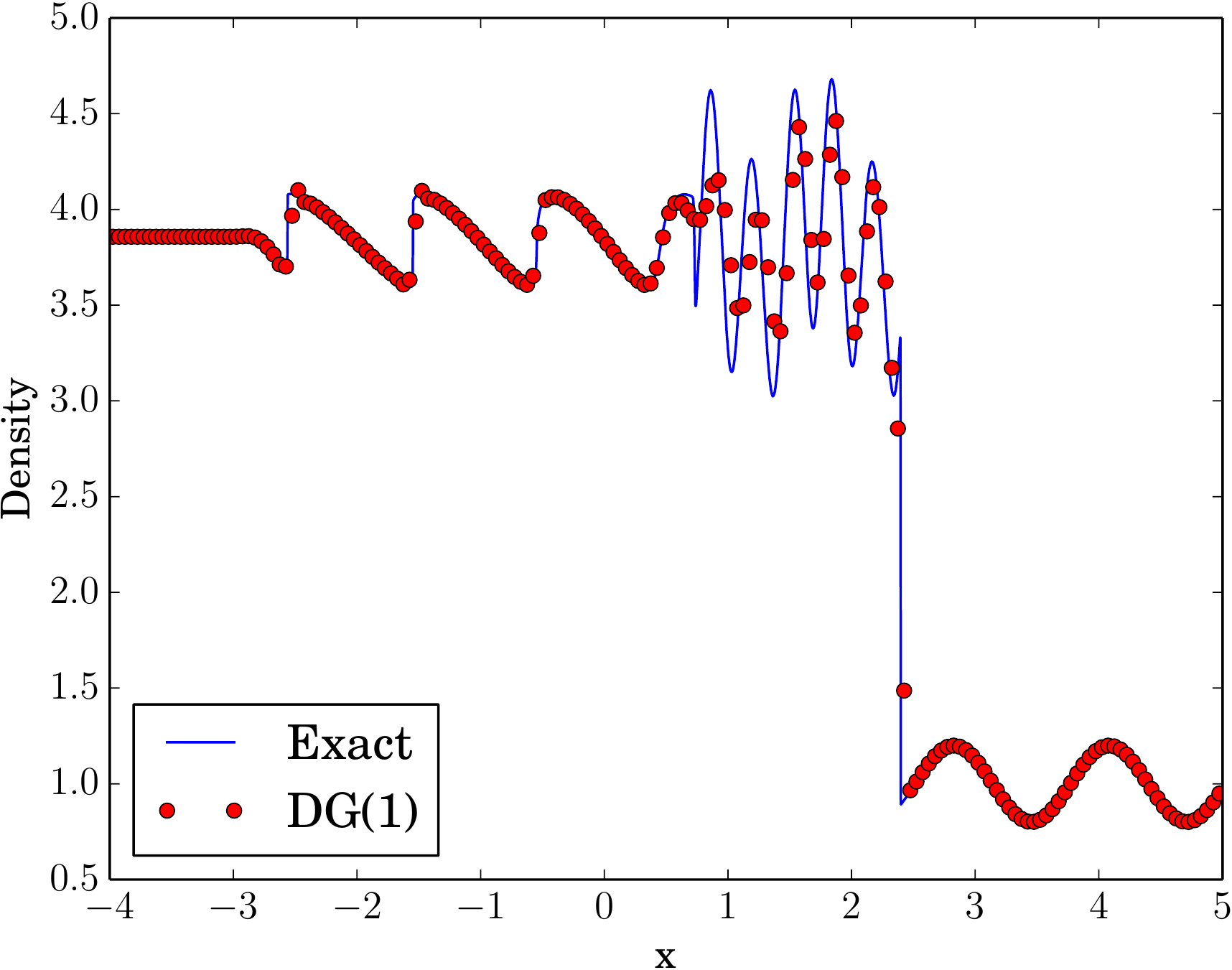}
                                                                               &
\includegraphics[width=0.48\textwidth]{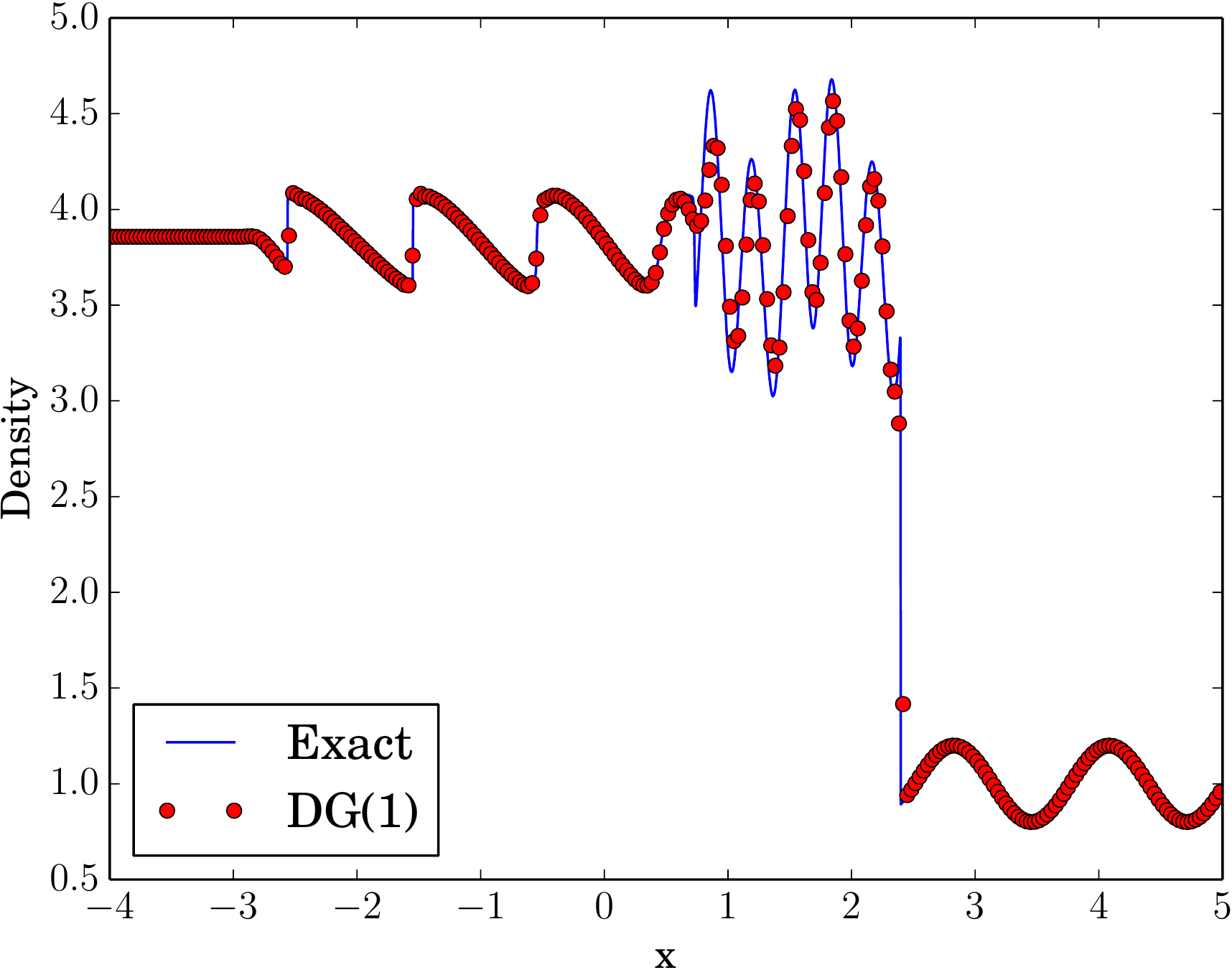} \\ (c) &
(d) \end{tabular} \end{center} \caption{Shu-Osher problem using Roe flux: (a)
static mesh, 200 cells, $M=0$ (b) moving mesh, 200 cells, $M=0$ (c) static mesh,
200 cells, $M=100$ (d) static mesh, 300 cells, $M=100$} \label{fig:so1}
\end{figure}

The initial condition is given
by~\cite{Shu1988439} \[ (\rho,v,p) = \begin{cases} (3.857143, 2.629369,
10.333333) &  \textrm{if } x < -4 \\ (1+0.2\sin(5x), 0.0, 1.0) &  \textrm{if } x
> -4 \end{cases} \] which involves a smooth sinusoidal density wave which
interacts with a shock. The domain is $[-5,+5]$ and the solution is computed
up to a final time of $T=1.8$. The solutions are shown in
figure~(\ref{fig:so1}a)-(\ref{fig:so1}b) on static and moving meshes using 200
cells and TVD limiter. The moving mesh scheme is considerably more accurate in
resolving the sinusoidal wave structure that arises after interaction with the
shock. In figure~(\ref{fig:so1}c) we compute the solution on static mesh with
TVB limiter and the parameter $M=100$ in equation~(\ref{eq:tvblim}). In this
case the solutions on static mesh are more accurate compared to the case of TVD
limiter but still not as good as the moving mesh results. The moving mesh result
has more than 200 cells in the interval $[-5,+5]$ at the final time since cells
enter the domain from the left side. Hence in figure~(\ref{fig:so1}d), we show
the static mesh results with 300 cells and using TVB limiter. The results are
further improved for the static mesh case but still not as accurate as the
moving mesh case. The choice of parameters in the TVB limiter is very critical
but we do not have a rigorous algorithm to choose a good value for this. Hence
it is still advantageous to use the moving mesh scheme which gives improved
solutions even with TVD limiter.

\begin{figure} \begin{center}
\includegraphics[width=0.8\textwidth]{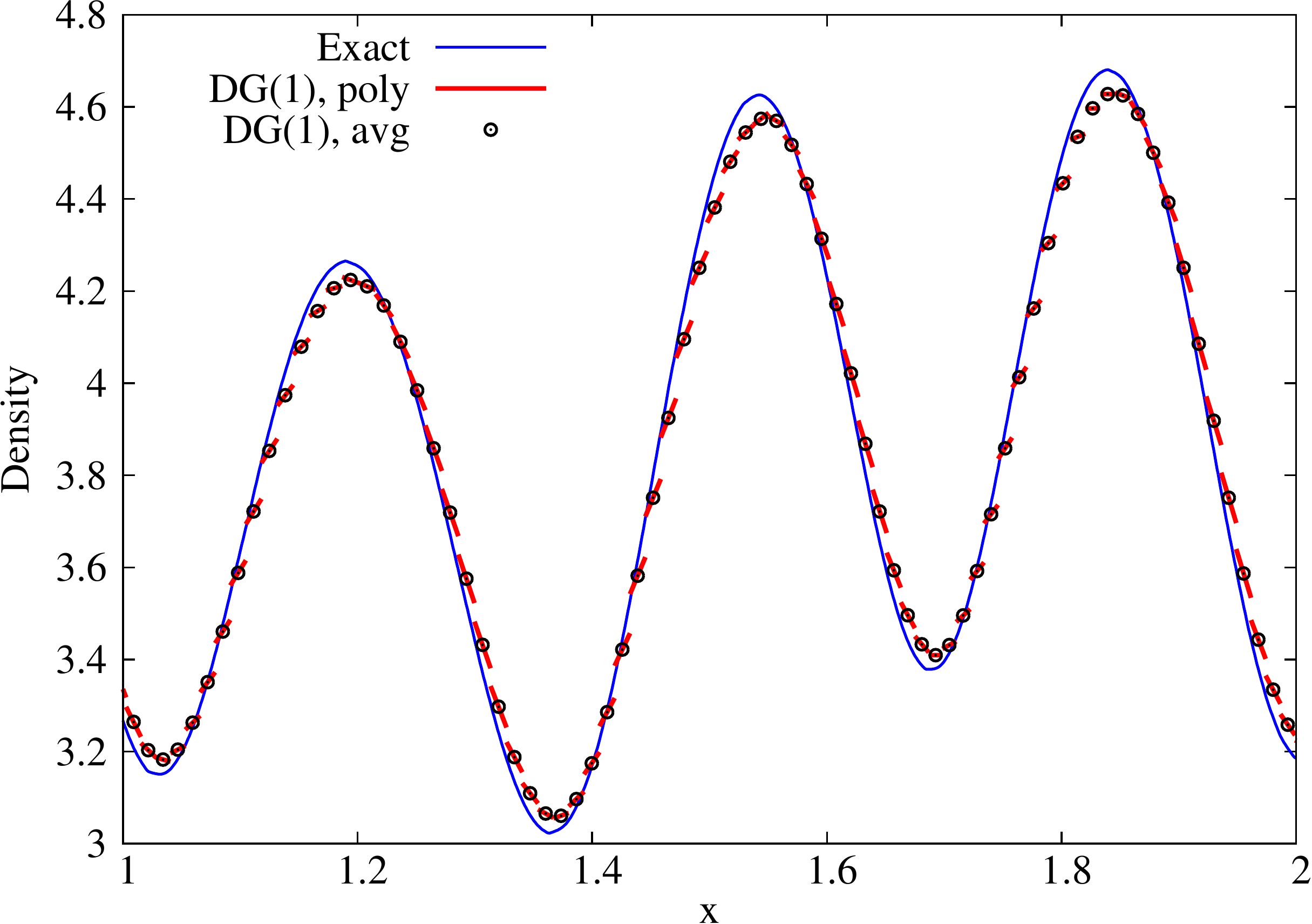}
\end{center} \caption{Shu-Osher problem using Roe flux on moving mesh}
\label{fig:so2} \end{figure}

\begin{figure} \begin{center}
  \includegraphics[width=0.8\textwidth]{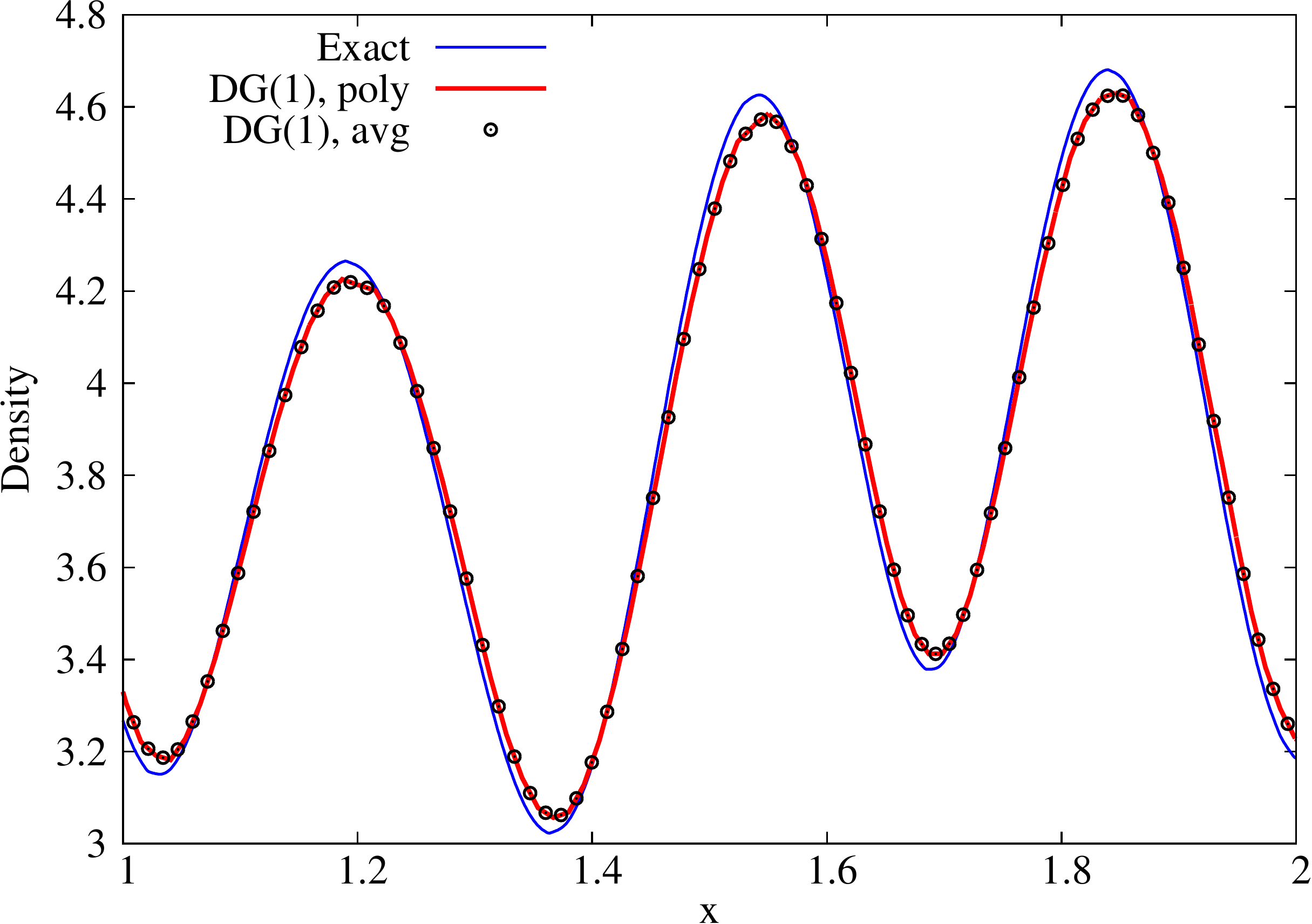}
 \end{center} \caption{Shu-Osher problem using modified Roe flux on
moving mesh} \label{fig:so3} \end{figure}

The above results show the ALE method is very accurate in terms of the cell
averages. In figure~(\ref{fig:so2}), we show a zoomed view of density and
pressure, where we also plot the linear polynomial solution. The slope of the
solution is not accurately predicted with the Roe scheme and there are spurious
contact discontinuities as the pressure and velocity are nearly continuous. This
behaviour is observed with all contact preserving fluxes like Roe, HLLC and
HLL-CPS but not with the Rusanov flux. Due to the almost Lagrangian character of
the scheme, the eigenvalue corresponding to the contact wave, $\lambda_2 = v-w$,
is nearly zero, which leads to loss of dissipation in the corresponding
characteristic field. If a spurious contact wave is generated during the violent
dynamics, then this wave will be preserved by the scheme leading to wrong
solutions. We modify the Roe scheme by preventing this eigenvalue from becoming
too small or zero, which is similar to the approach used for the entropy fix.
The eigenvalue $|\lambda_2|$ used in the dissipative part of the Roe flux is
determined from \[ |\lambda_2| = \begin{cases} |v-w| & \textrm{if } |v-w| >
\delta = \alpha c \\ \half(\delta + |v-w|^2/\delta)         & \textrm{otherwise}
\end{cases} \] With this modification and using $\alpha=0.1$, the solution on
moving mesh is shown in figure~(\ref{fig:so3}) and we do not observe the
spurious contact discontinuities which arise with the standard Roe flux, while
at the same time, the solution accuracy compares favourably with the previous
results that did not use the eigenvalue fix.  

\begin{figure}
  \begin{center} \begin{tabular}{cc}
    \includegraphics[width=0.48\textwidth]{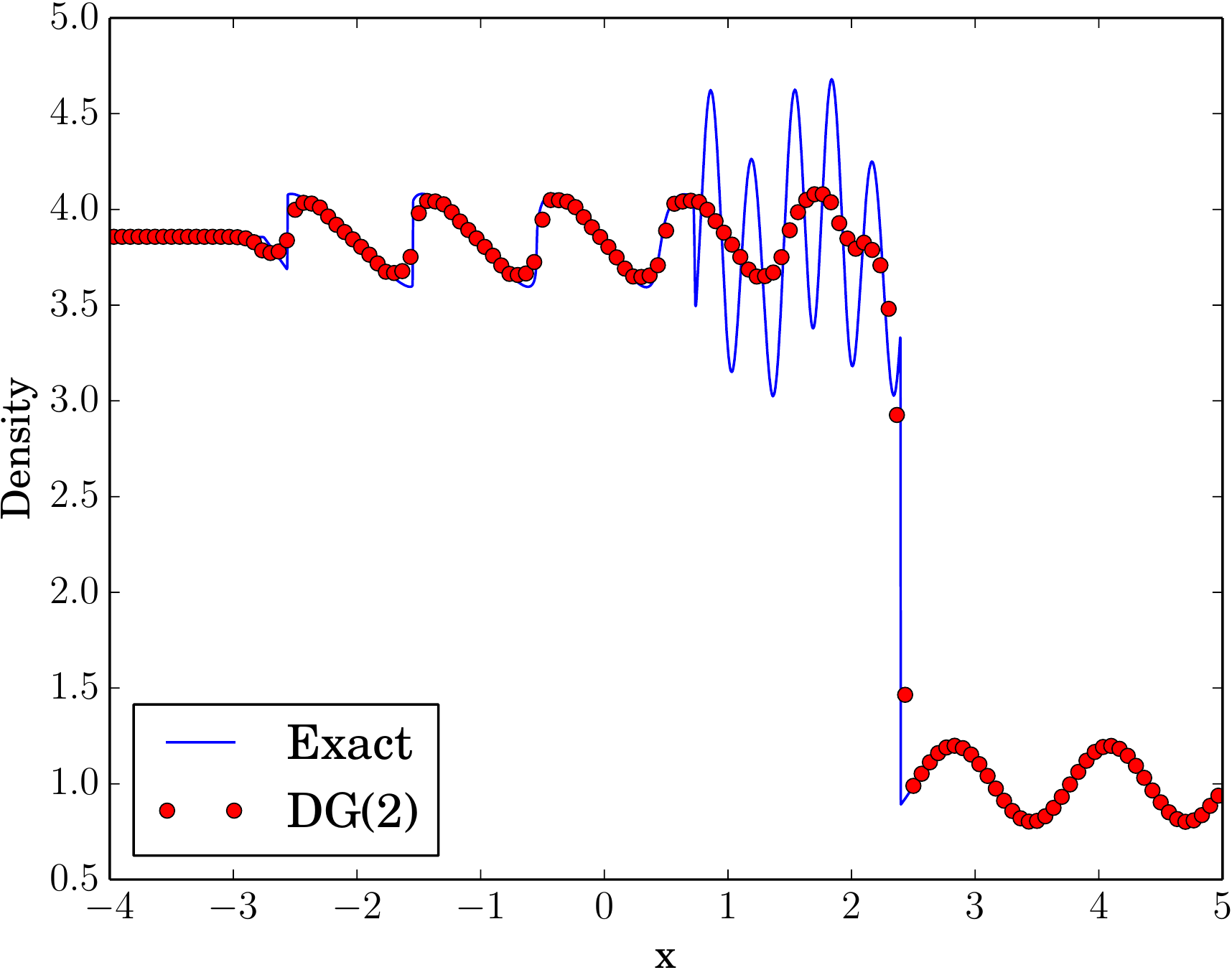}
    &
    \includegraphics[width=0.48\textwidth]{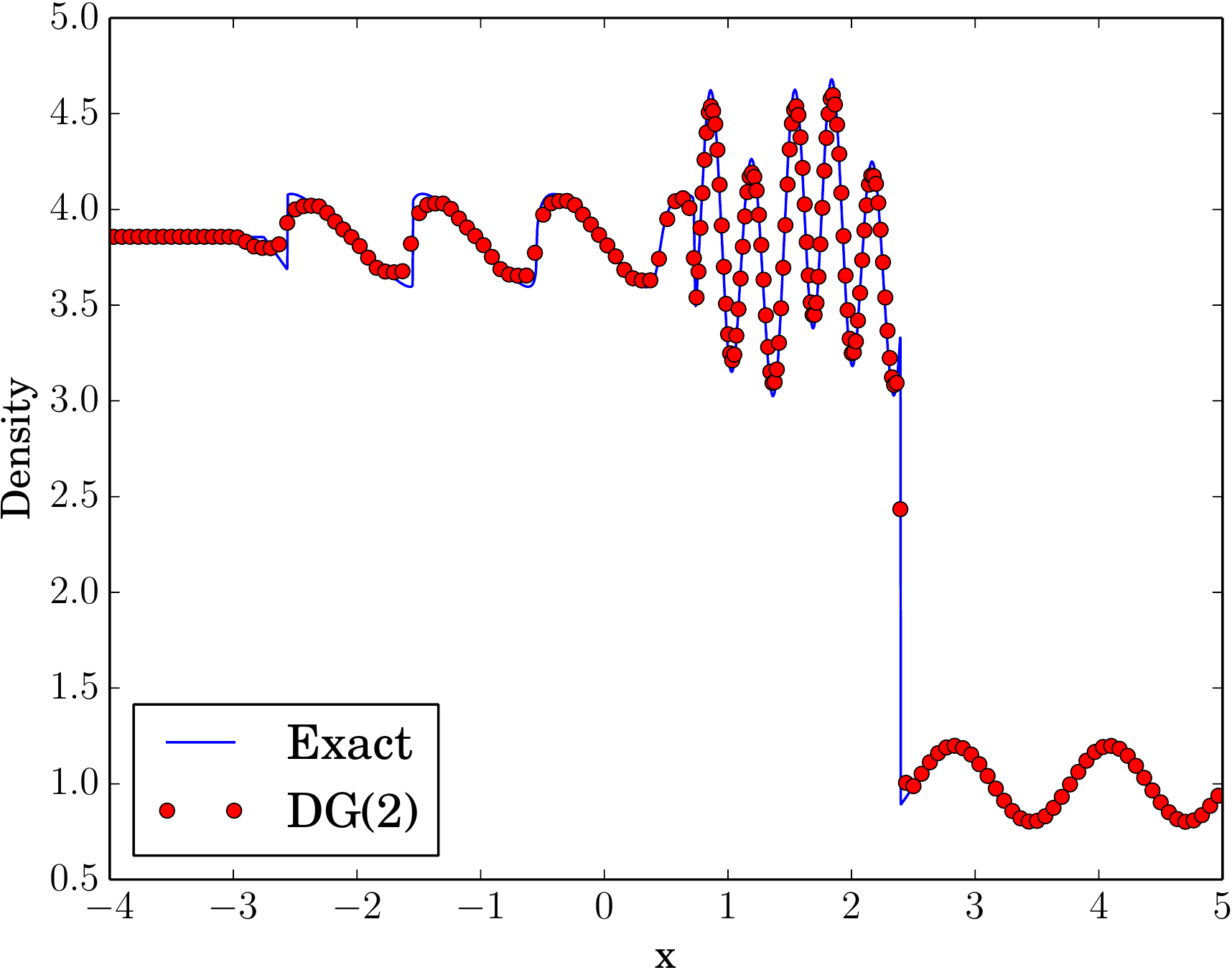}
    \\ (a) & (b) \end{tabular} \end{center} \caption{Shu-Osher problem using
  modified Roe flux, TVD limiter, quadratic polynomials and 150 cells. (a)
  static mesh, (b) moving mesh} \label{fig:so4a} \end{figure}

\begin{figure} \begin{center} \begin{tabular}{cc}
\includegraphics[width=0.48\textwidth]{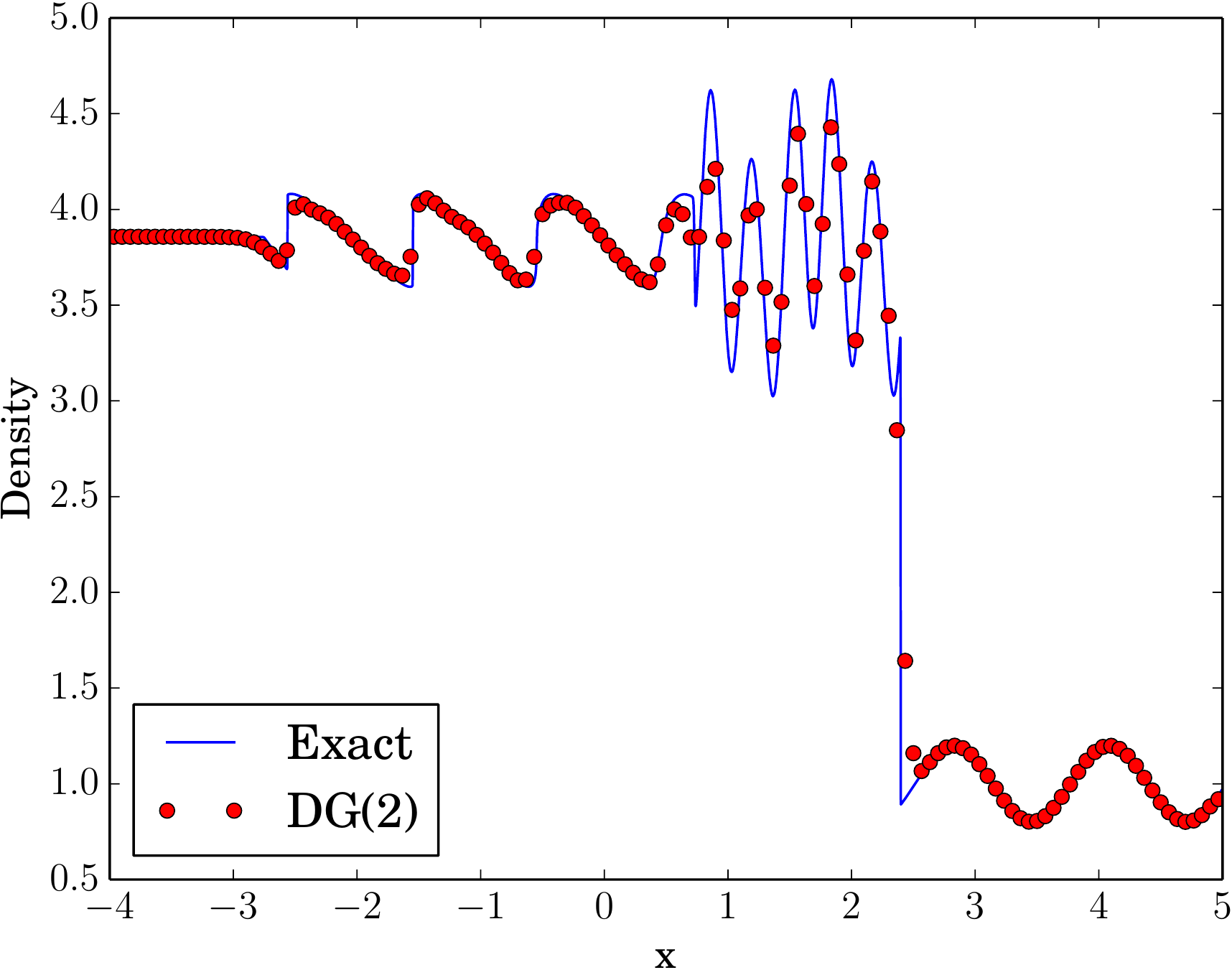} &
\includegraphics[width=0.48\textwidth]{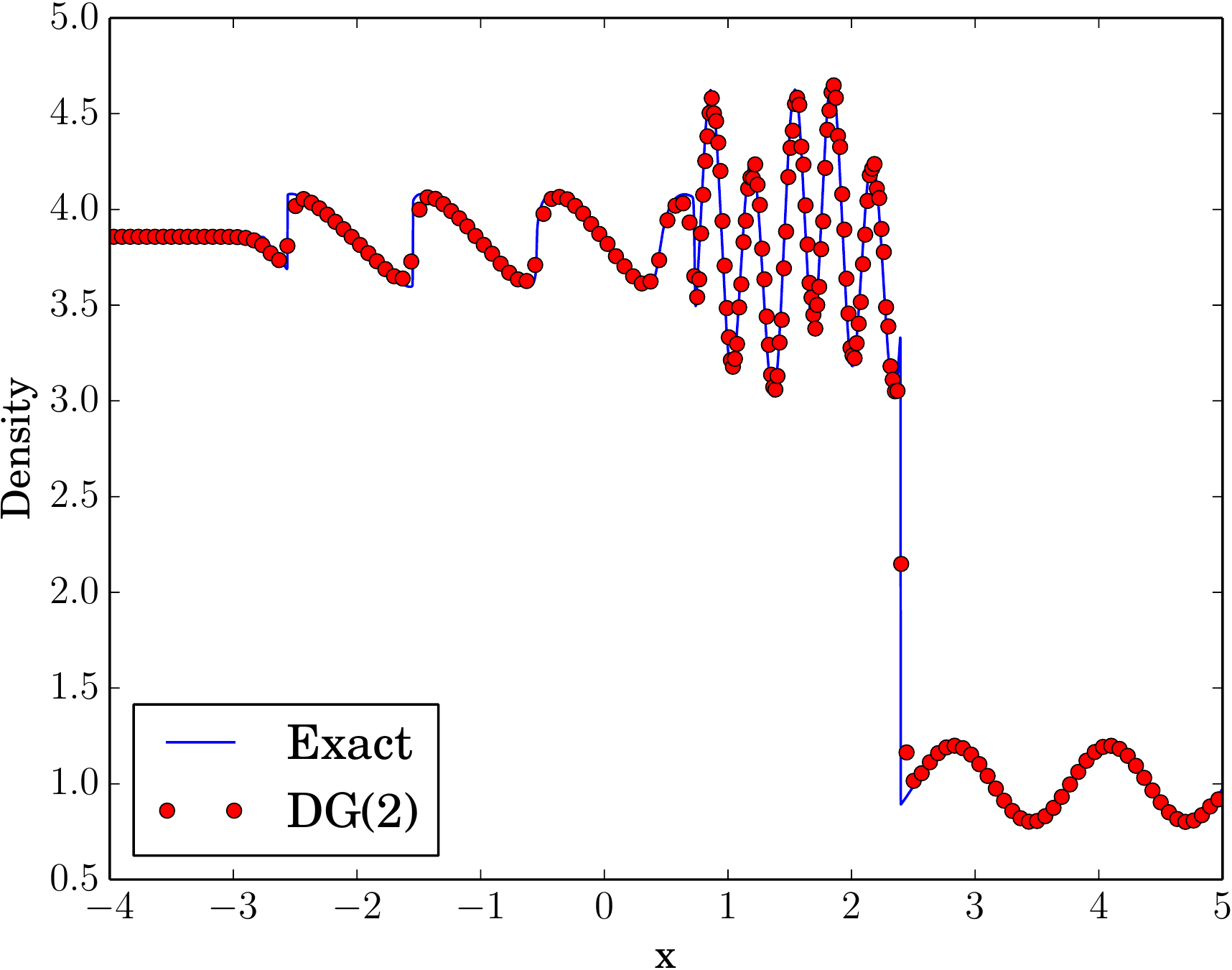} \\ (a)
                                                                          & (b)
                                                                          \end{tabular}
    \end{center} \caption{Shu-Osher problem using modified Roe flux, WENO
    limiter, quadratic polynomials and 150 cells. (a) static mesh, (b) moving
  mesh} \label{fig:so4b} \end{figure}

\begin{figure} \begin{center} \begin{tabular}{cc}
\includegraphics[width=0.48\textwidth]{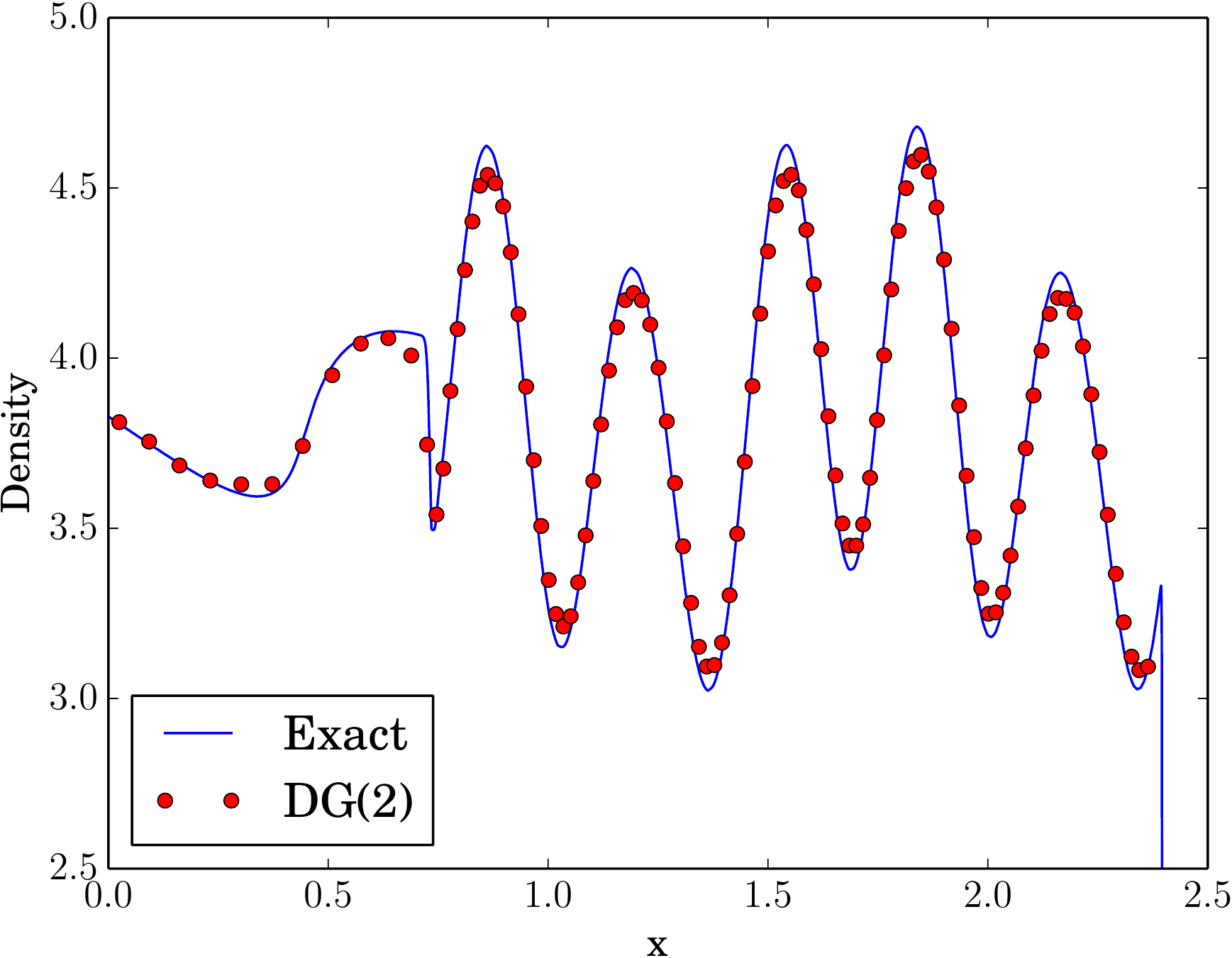}
& \includegraphics[width=0.48\textwidth]{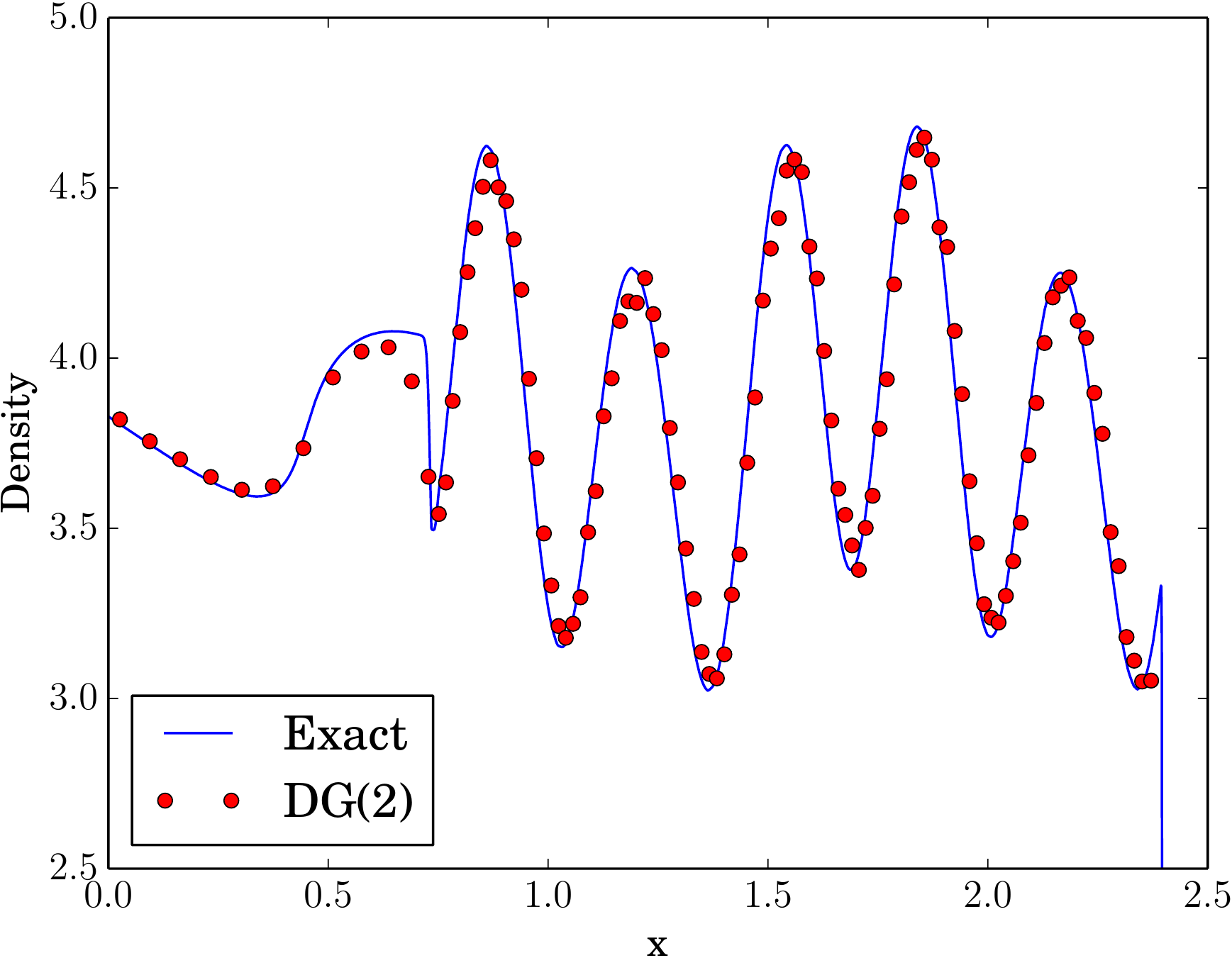}
\\ (a) & (b) \end{tabular} \end{center} \caption{Shu-Osher problem using
modified Roe flux, moving mesh, quadratic polynomials and 150 cells. (a) TVD
limiter, (b) WENO limiter} \label{fig:so5} \end{figure}

We next compute the solutions using quadratic polynomials.
Figure~(\ref{fig:so4a}) shows the results obtained with the TVD limiter which
shows the dramatically better accuracy that is achieved on moving mesh compared
to static mesh. In figure~(\ref{fig:so4b}) we perform the same computation with
a WENO limiter taken from~\cite{Zhong:2013:SWE:2397205.2397446}. The static mesh
results are now improved over the case of TVD limiter but still not as good as
the moving mesh results in terms of capturing the extrema. In
figure~(\ref{fig:so5}) we show a zoomed view of the results on moving mesh with
TVD and WENO limiter. We see that the TVD limiter is also able to capture all
the features and is almost comparable to the WENO limiter.  

\subsection{Titarev-Toro problem }
\begin{figure} \centering
\includegraphics[width=0.8\textwidth]{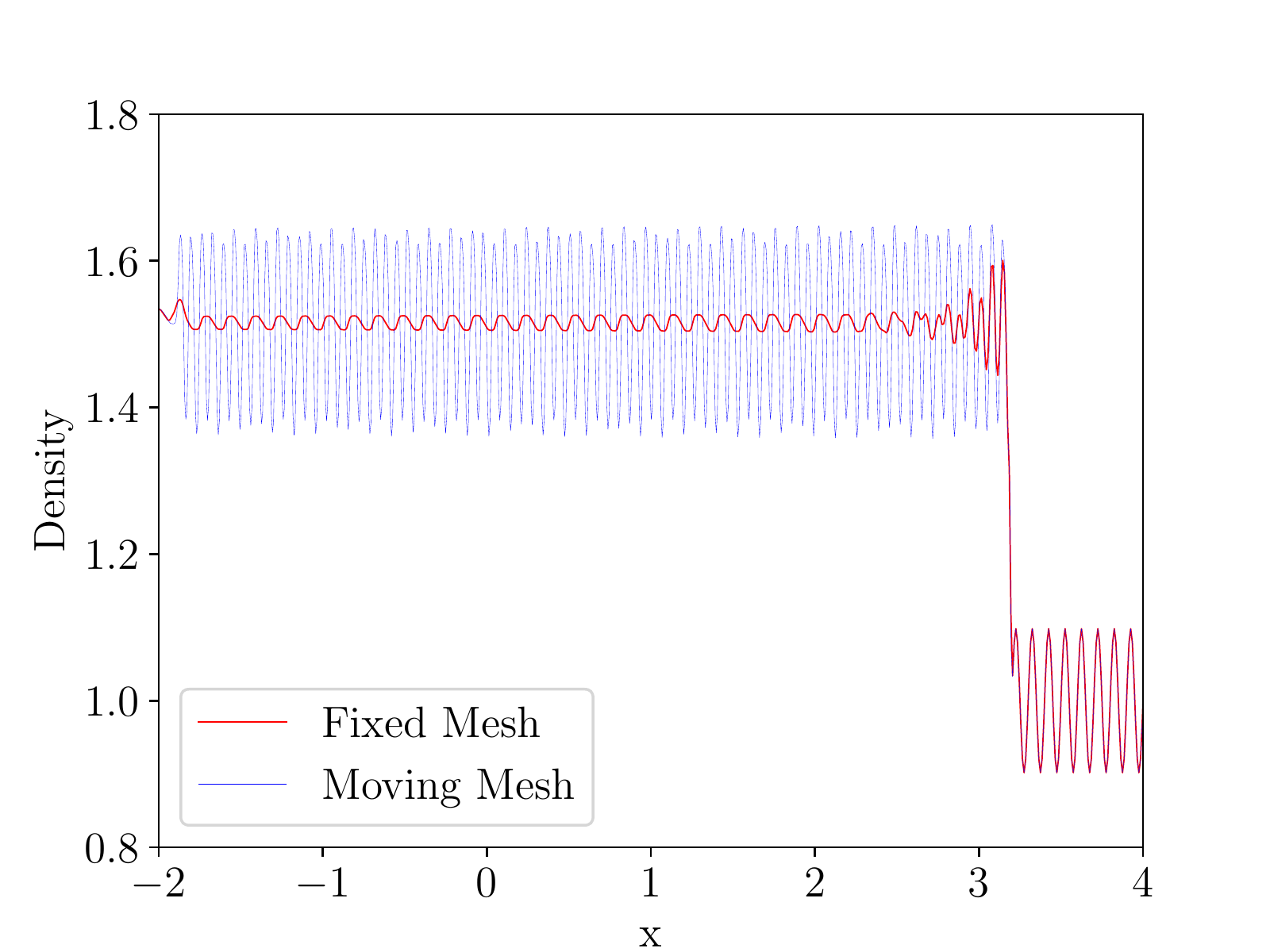}
\caption{Titarev Problem with HLLC flux, 1000 cells and TVD limiter}
\label{fig:ttoro1}
\end{figure}

\begin{figure} \centering
\includegraphics[width=0.8\textwidth]{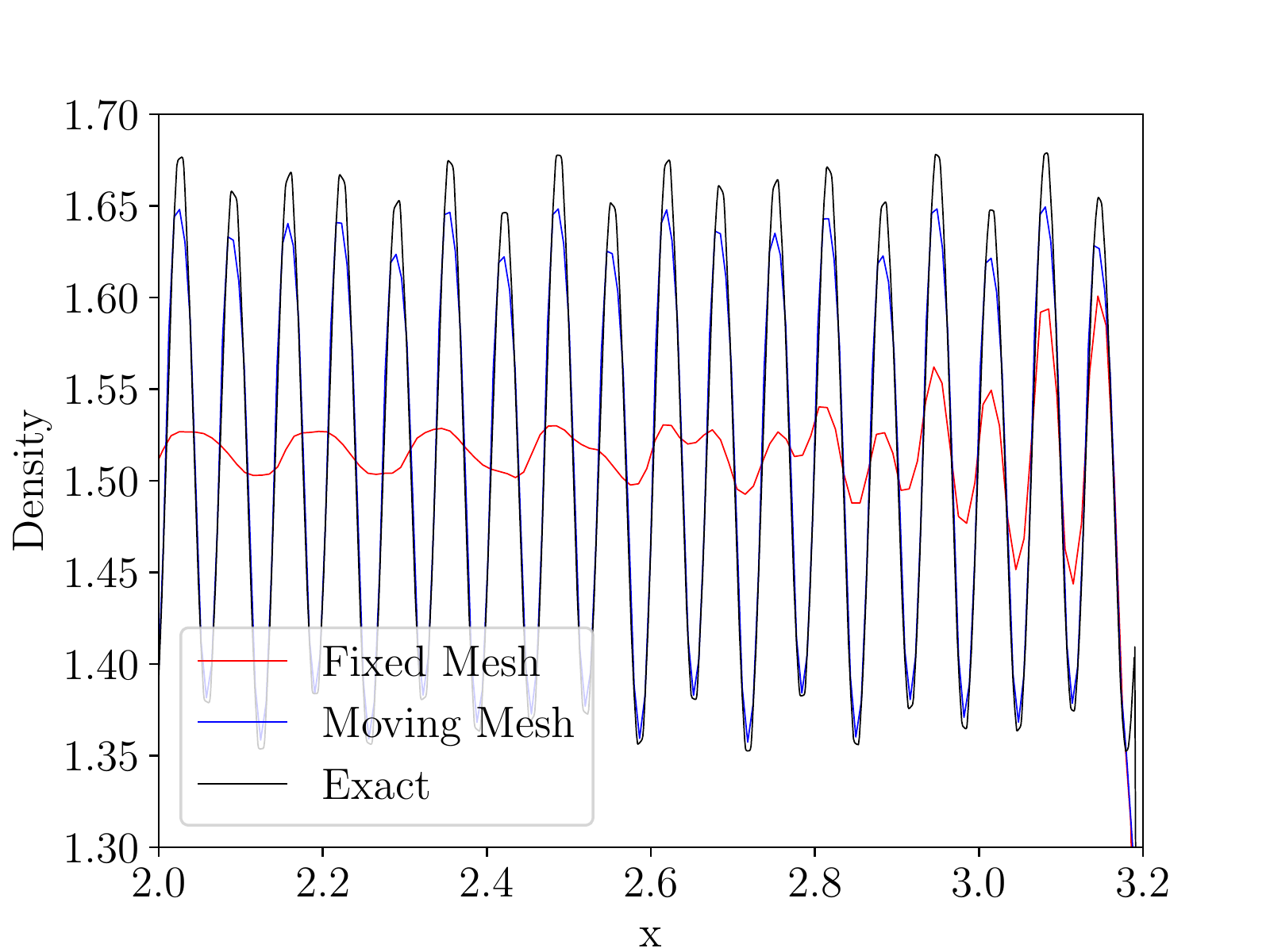}
\caption{Titarev Problem with HLLC flux, 1000 cells and TVD limiter. (Zoomed
Version)}
\label{fig:ttoro2}
\end{figure}

Titarev-Toro problem is an extension of the Shu-Osher problem \cite{Titarev2014}
to test a severely oscillatory wave interacting with a shock wave. It aims to test the
ability of higher-order methods to capture the extremely high frequency waves. The initial condition is given by \begin{align} (\rho, v, p) = \begin{cases}
(1.515695, 0.523346, 1.805), & -5 < x \leq -4.5 \\ (1 + 0.1 \sin(20 \pi x), 0,
1), & -4.5 < x \leq 5 \end{cases} \end{align}
The computation is carried out on a mesh of 1000 cells with the final time
$T=5$ and the density at this final time is shown in Figures~(\ref{fig:ttoro1}),~(\ref{fig:ttoro2}). The fixed mesh is not able to resolve the high frequency oscillations due to dissipation in the fluxes and the TVD limiter, but the ALE scheme gives an excellent resolution of these high frequency oscillations. Note that the ALE scheme also uses the same TVD limiter but it is still able to resolve the solution to a very degree of accuracy. This result again demonstrates the superior accuracy that can be achieved by using a nearly Lagrangian ALE scheme in problems involving interaction of shocks and smooth flow structures.

\subsection{123 problem} 
\begin{figure} \begin{center}
  \begin{tabular}{cc}
    \includegraphics[width=0.48\textwidth]{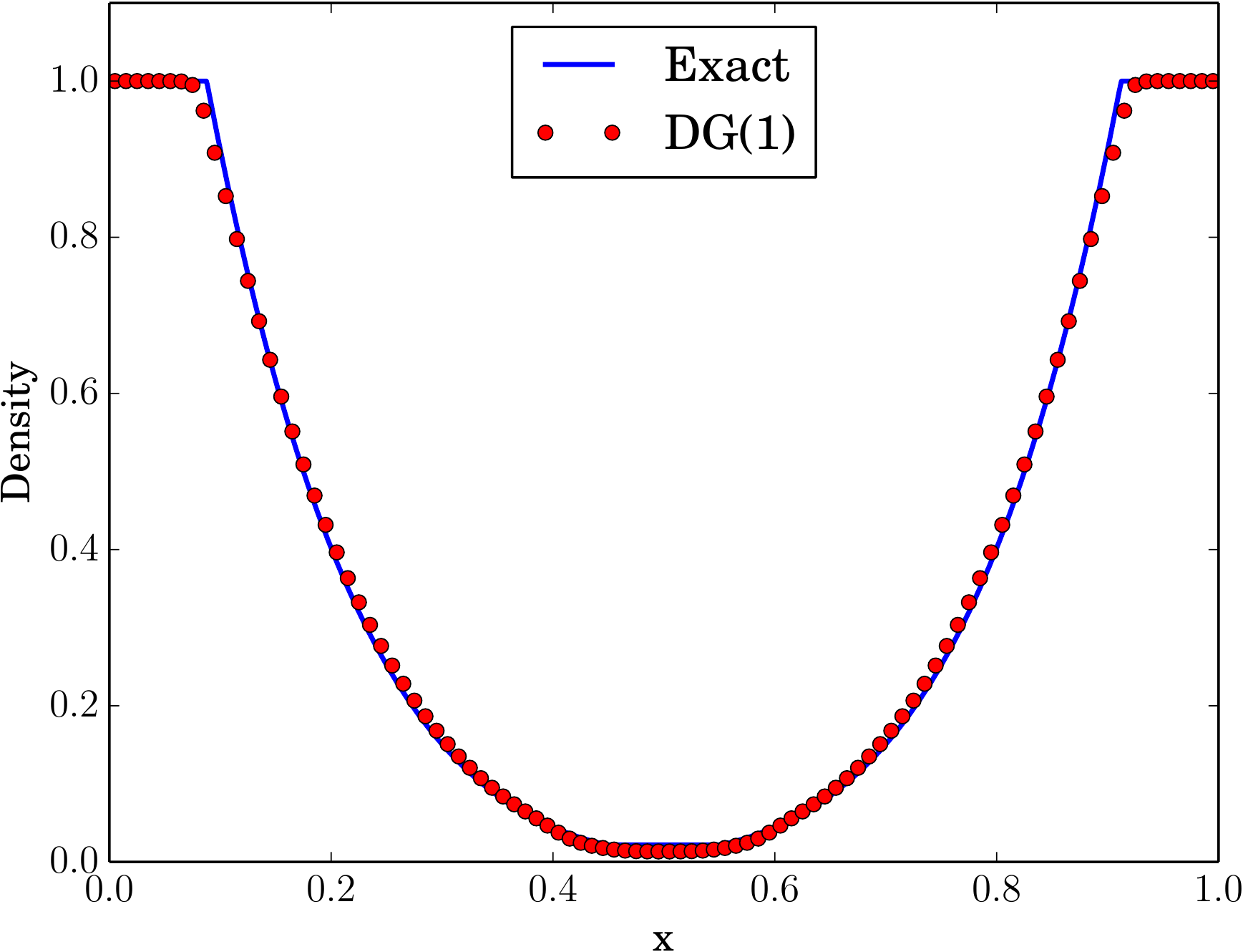} &
    \includegraphics[width=0.48\textwidth]{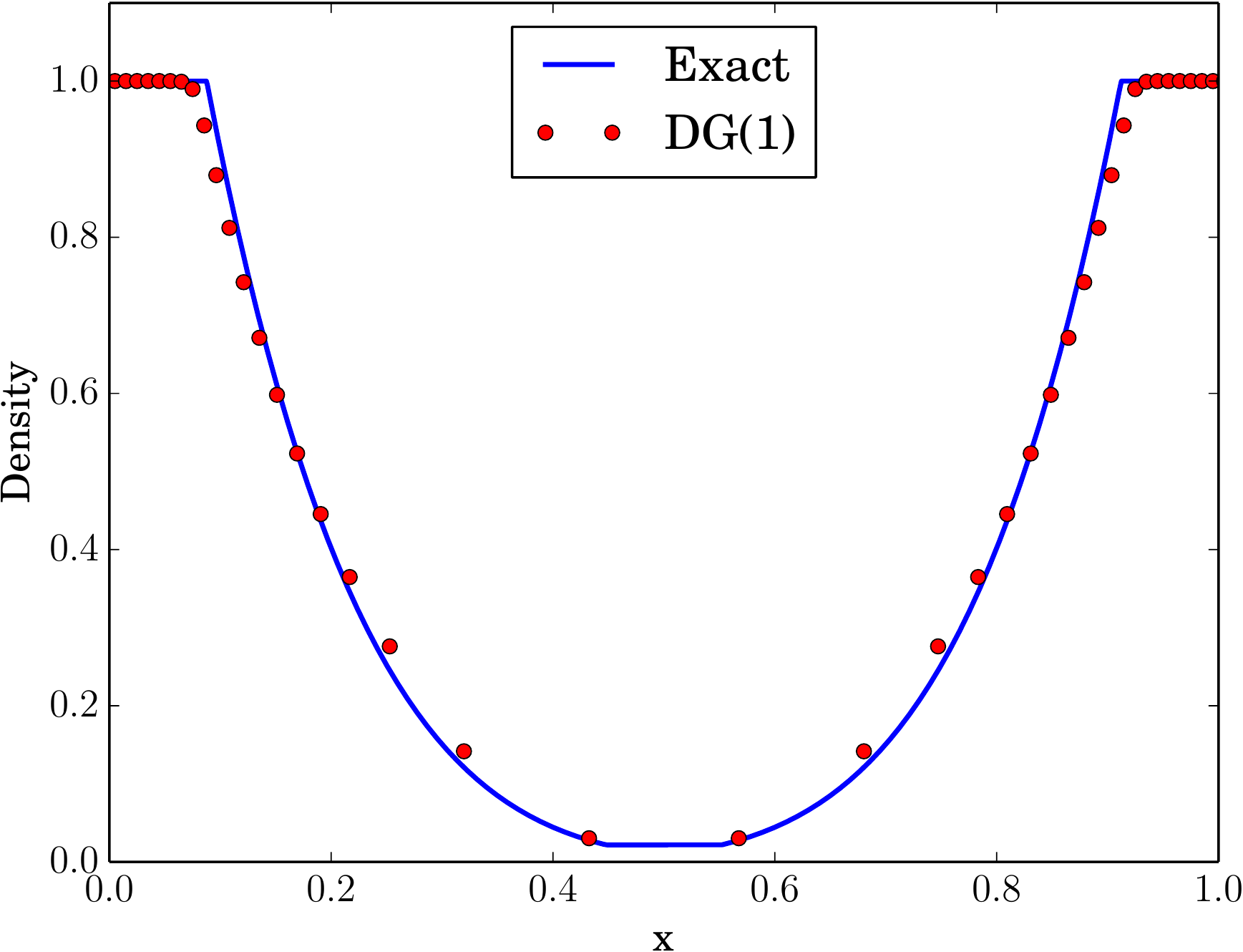}  \\ (a) &
    (b) \end{tabular} \end{center} \caption{123 problem using HLLC flux and
  100 cells: (a) static mesh, (b) moving mesh} \label{fig:lowd1} \end{figure}

\begin{figure} \begin{center} \begin{tabular}{cc}
\includegraphics[width=0.48\textwidth]{lowd_mov_rho} &
\includegraphics[width=0.48\textwidth]{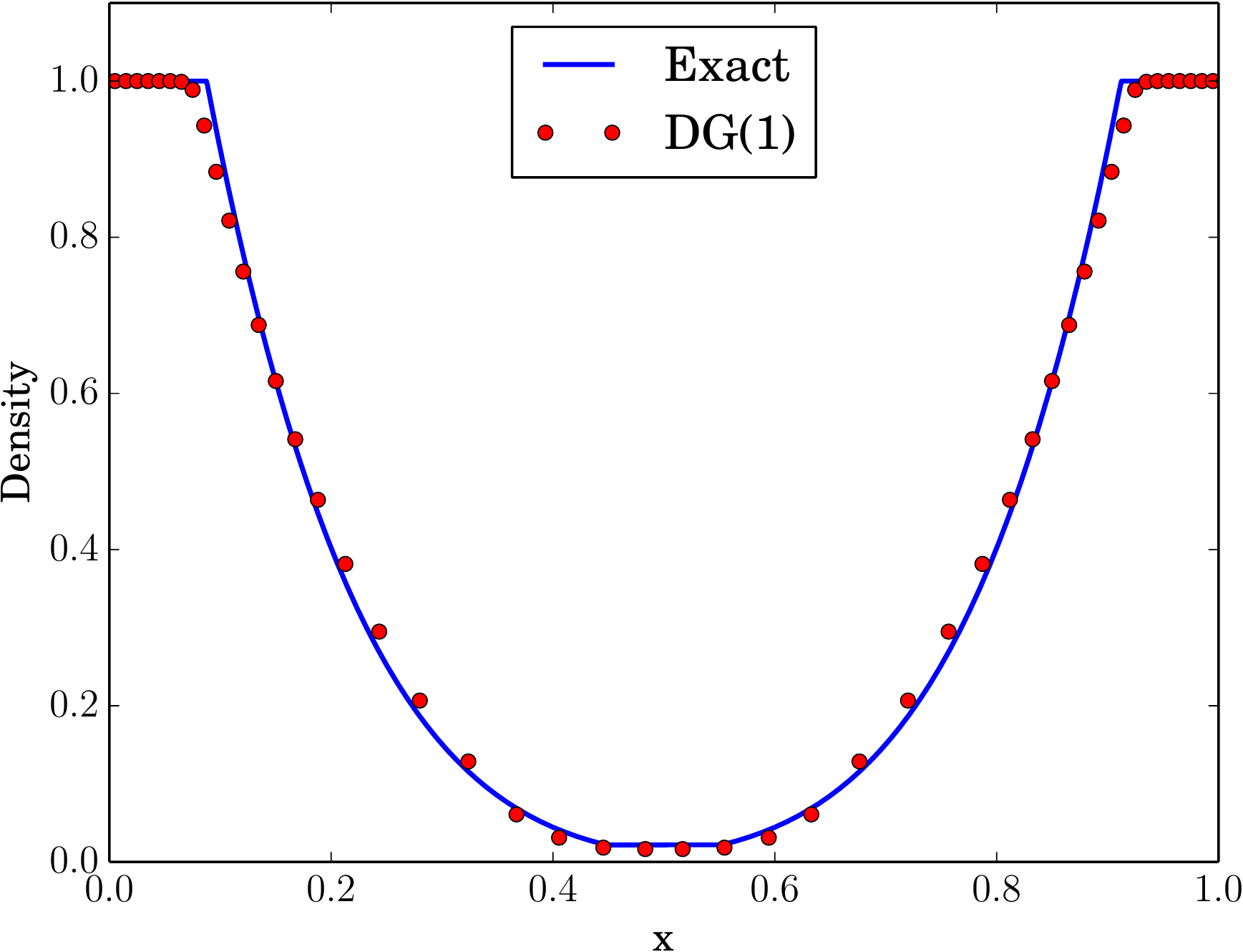}  \\ (a) &
(b) \end{tabular} \caption{123 problem using HLLC flux and grid refinement: (a)
static mesh, (b) moving mesh with mesh adaptation ($h_{max}=0.05$) leading to
108 cells at final time} \label{fig:lowd2} \end{center} \end{figure}

The initial condition is given by~\cite{torobook} \[
(\rho,v,p) = \begin{cases} (1.0, -2.0, 0.4) & x < 0.5 \\ (1.0, +2.0, 0.4) & x >
0.5 \end{cases} \] The computational domain is $[0,1]$ and the final time is
$T=0.15$. The density using 100 cells is shown in figure~(\ref{fig:lowd1}) with
static and moving meshes. The mesh motion does not significantly improve the
solution compared to the static mesh case since the solution is smooth. On the
contrary, the mesh becomes rather coarse in the expansion region, though the
solution is still well resolved. However, severe expansion may lead to very
coarse meshes which may be undesirable. To prevent very coarse cells, we switch
on the mesh refinement algorithm as described before and use the upper bound on
the mesh size as $\hmax=0.05$. The resulting solution is shown in
figure~(\ref{fig:lowd2}) where the number of cells has increased to 108 at the
time shown. The central expansion region is now resolved by more uniformly sized
cells compared to the case of no grid refinement.  

\begin{figure} \centering
\includegraphics[width=0.6\textwidth]{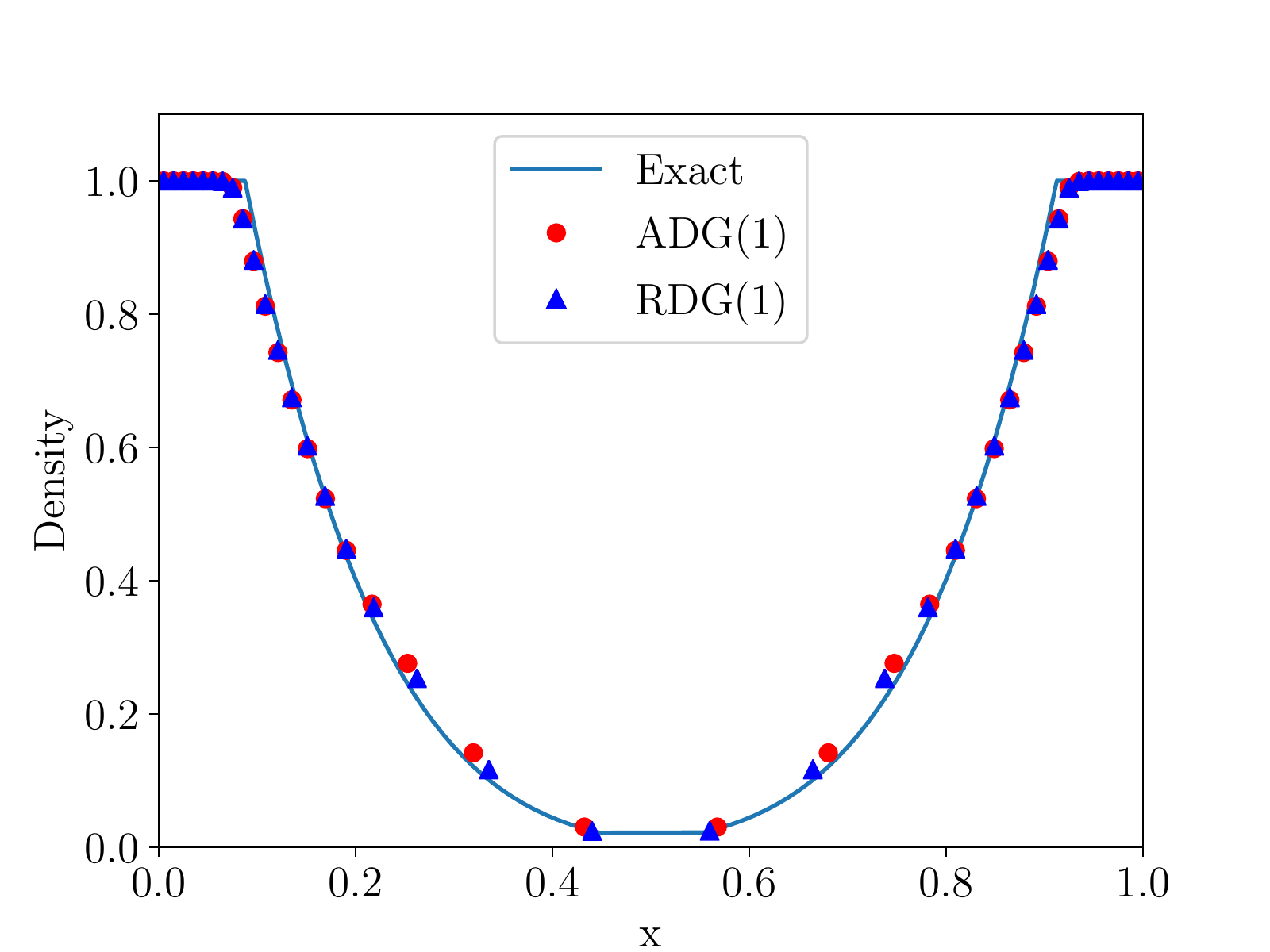}
\caption{123 problem using HLLC flux, 100 cells and TVD limiter. ADG : Average
Velocity, RDG : Linearized Riemann Velocity} \end{figure}

\subsection{Blast problem} 
\begin{figure} \begin{center} \begin{tabular}{cc}
  \includegraphics[width=0.48\textwidth]{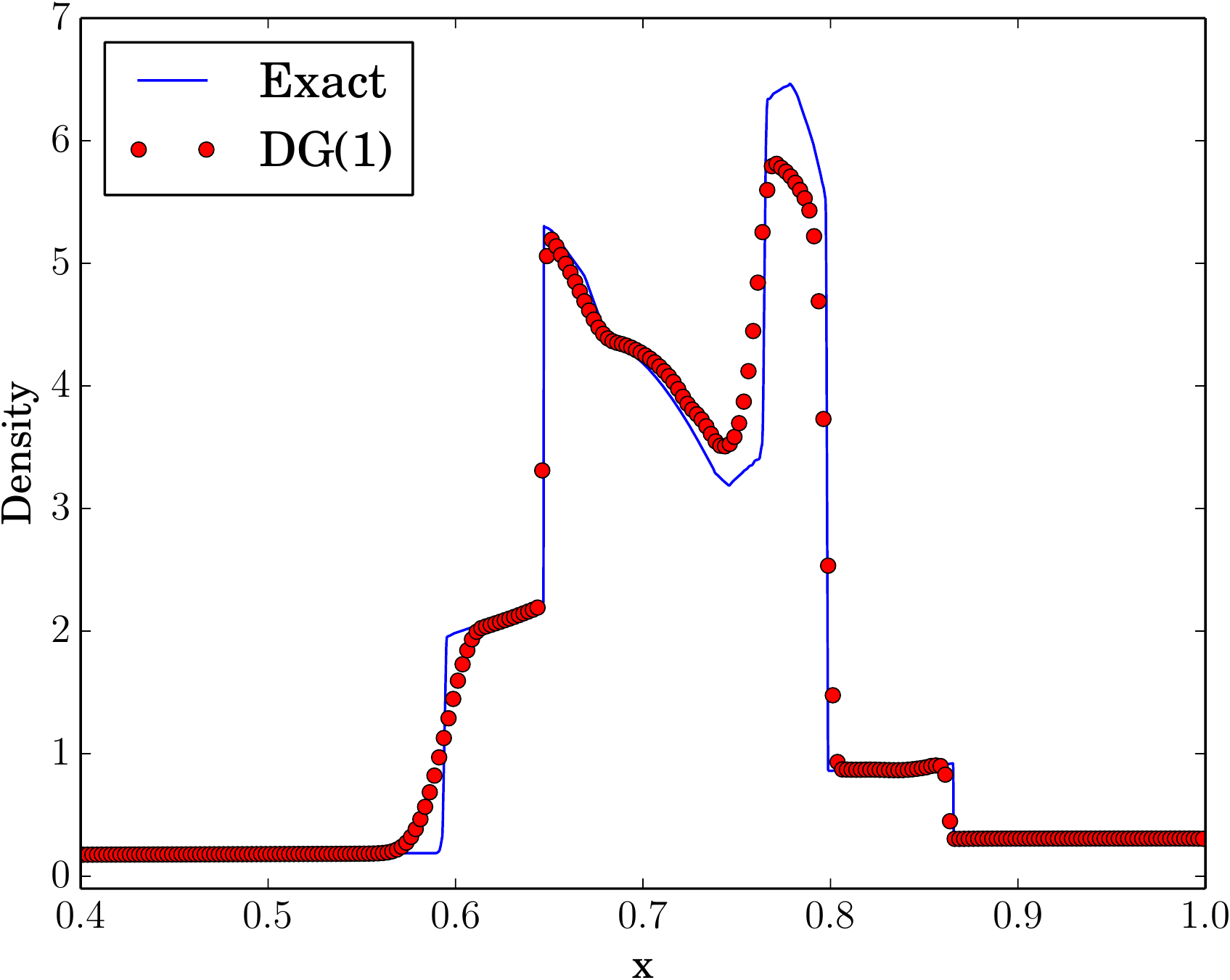}
  & \includegraphics[width=0.48\textwidth]{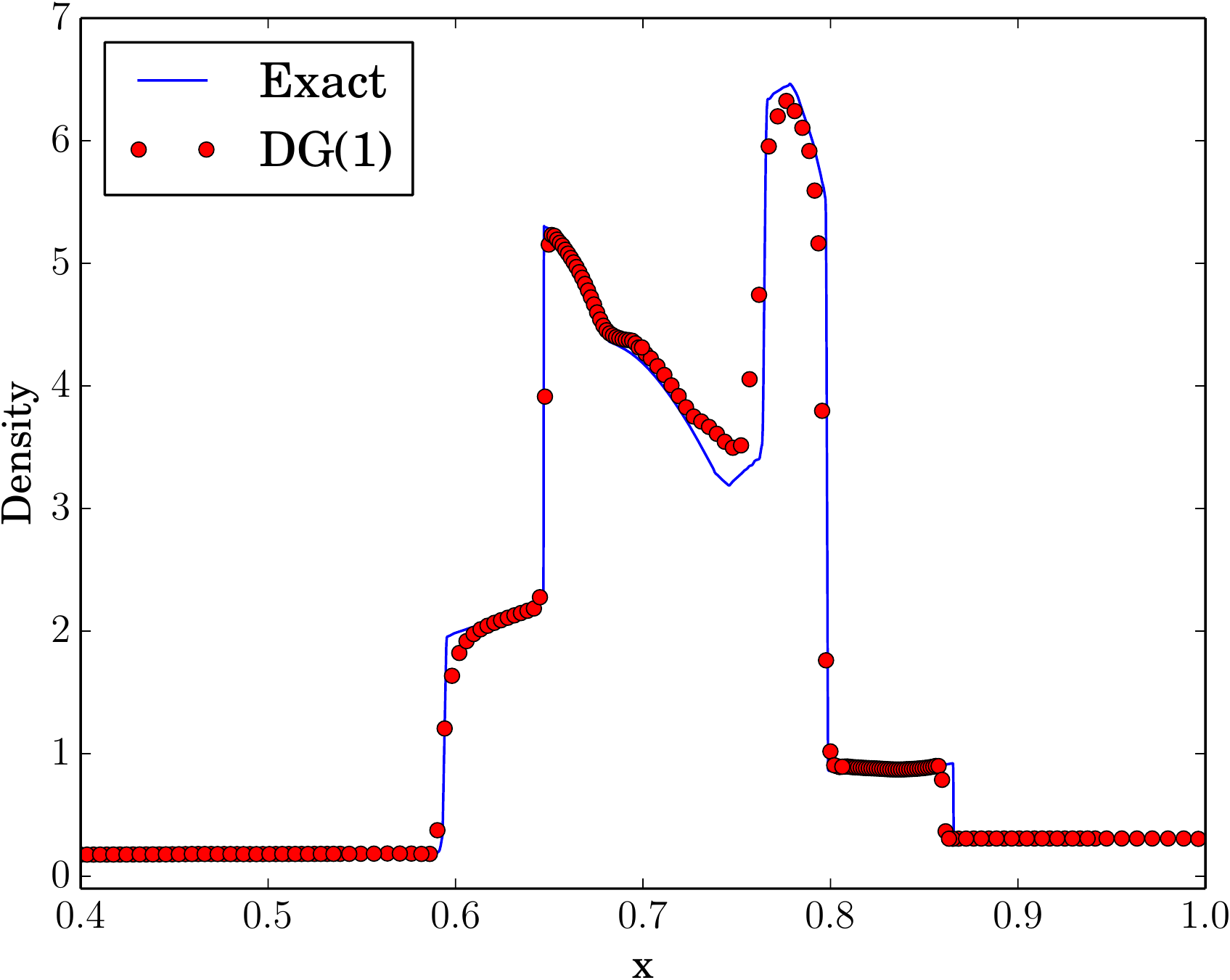}
  \\ (a) & (b) \end{tabular} \caption{Blast problem using HLLC flux and 400
cells. (a) static mesh, (b) moving mesh with adaptation ($h_{min}=0.001$)
leading to 303 cells at final time.} \label{fig:blast} \end{center} \end{figure}

\begin{figure}
  \begin{center} \begin{tabular}{cc}
    \includegraphics[width=0.48\textwidth]{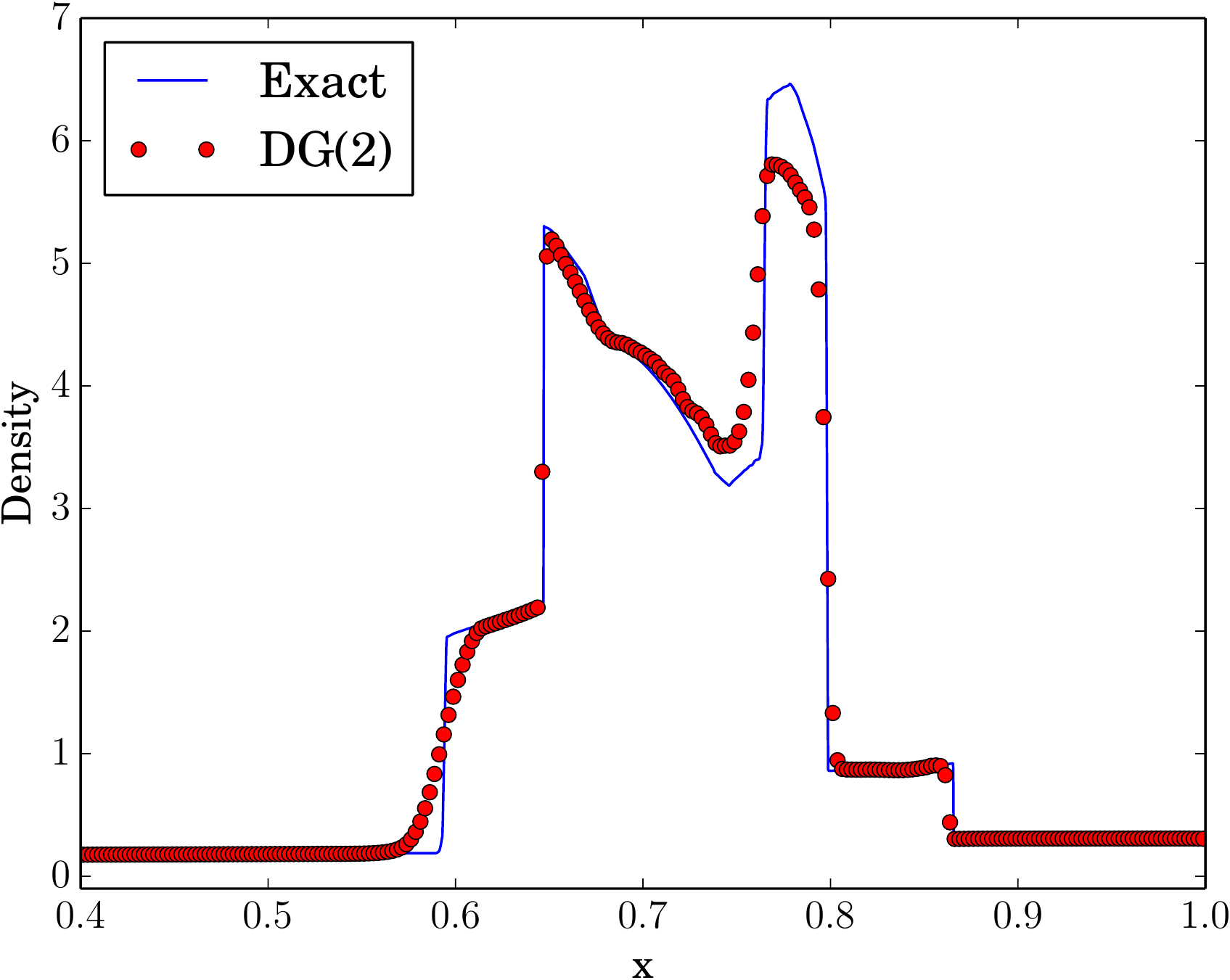}
    &
    \includegraphics[width=0.48\textwidth]{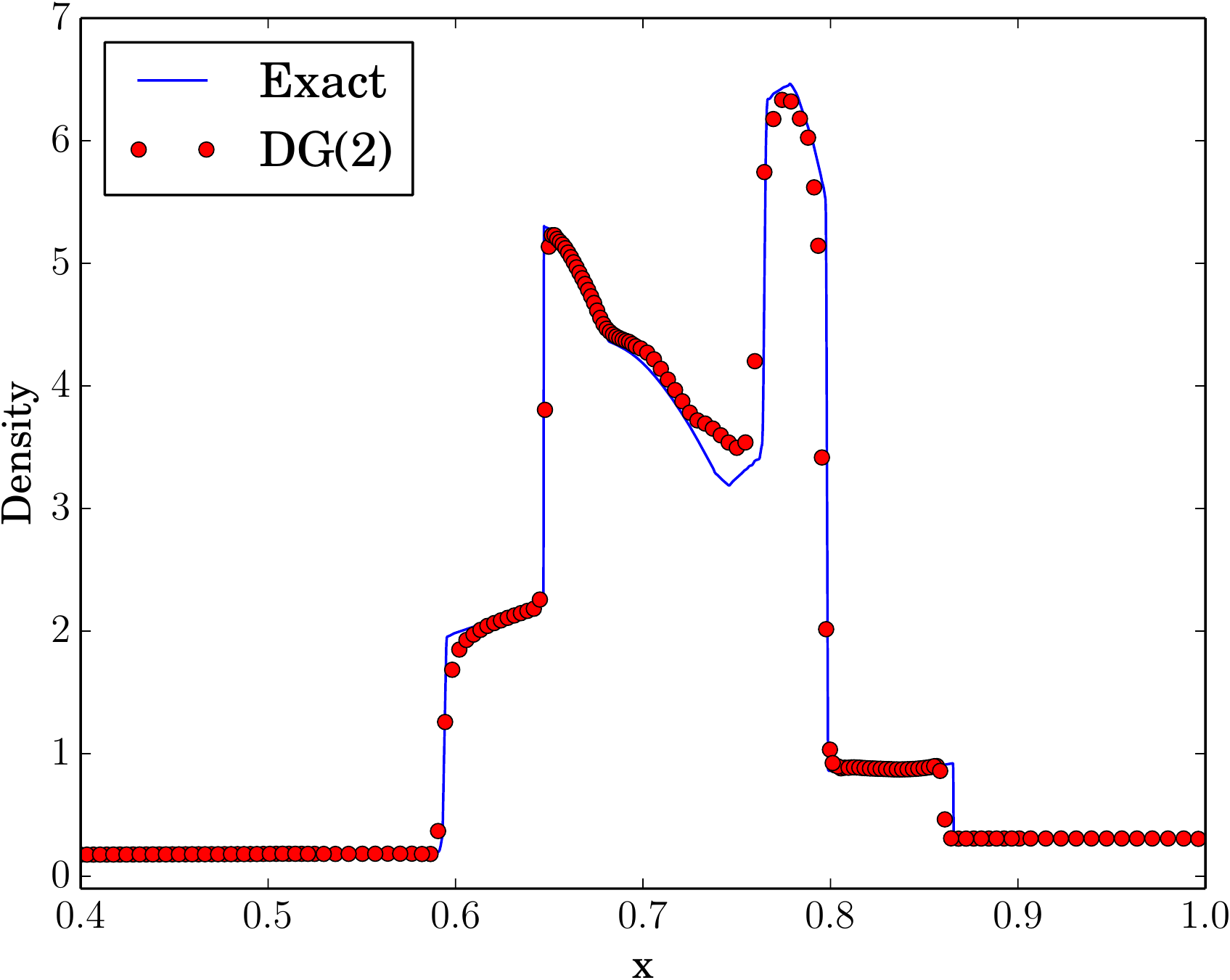}
    \\ (a) & (b) \end{tabular} \caption{Blast problem using HLLC flux, quadratic
  polynomials and 400 cells. (a) static mesh, (b) moving mesh with adaptation
  ($h_{min}=0.001$) leading to 293 cells at final time.} \label{fig:blast2}
    \end{center} \end{figure}

The initial condition is given
by~\cite{Woodward1984115} \[ (\rho,v,p) = \begin{cases} (1.0, 0.0, 1000.0) & x <
0.1 \\ (1.0, 0.0, 0.01) & 0.1 < x < 0.9 \\ (1.0, 0.0, 100.0) & x > 0.9
\end{cases} \] with a domain of $[0,1]$ and the final time is $T=0.038$. A
reflective boundary condition is used at $x=0$ and $x=1$. A mesh of 400 cells is
used for this simulation and in case of moving mesh, we perform grid adaptation
with $h_{min}=0.001$ since some cells become very small during the collision of
the two shocks. The positivity preserving limiter
of~\cite{Zhang:2010:PHO:1864819.1865066} is applied together with TVD limiter
and HLLC flux. The static mesh results shown in figure~(\ref{fig:blast}a)
indicate too much numerical viscosity in the contact wave around $x=0.6$. This
wave is more accurately resolved in the moving mesh scheme as seen in
figure~(\ref{fig:blast}b) which is an advantage due to the ALE scheme and is a very good indicator of the scheme accuracy as this is a very challenging feature to compute accurately.
We next compute the same problem using quadratic polynomials with all other
parameters being as before. The solutions are shown in figure~(\ref{fig:blast2})
and indicate that the Lagrangian moving mesh scheme is more accurate in
resolving the contact discontinuity. The higher polynomial degree does not show
any major improvement in the solution compared to the linear case, which could be a
consequence of the strong shock interactions present in this problem, see
figure~(4.11-4.12) in~\cite{Zhong:2013:SWE:2397205.2397446} and figure~(3.7) in
\cite{Zhu2016110} in comparison to current results.  

\begin{figure} \centering
\includegraphics[width=0.8\textwidth]{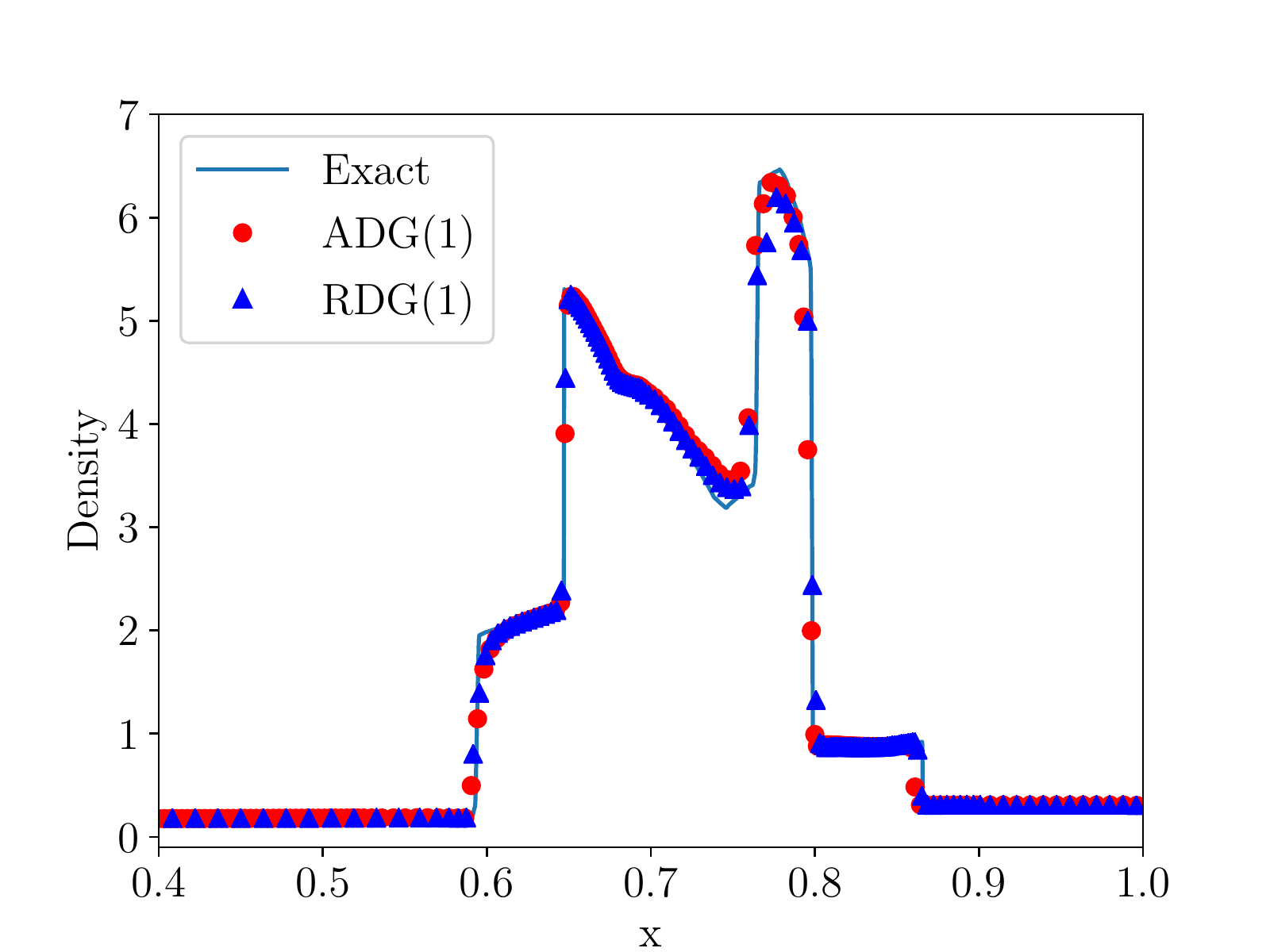}
\caption{Blast problem using HLLC flux, 100 cells and TVD limiter. ADG : Average
Velocity, RDG : Linearized Riemann Velocity} \end{figure}

\subsection{Le Blanc shock tube test case}
\begin{figure} \centering
\includegraphics[width=\textwidth]{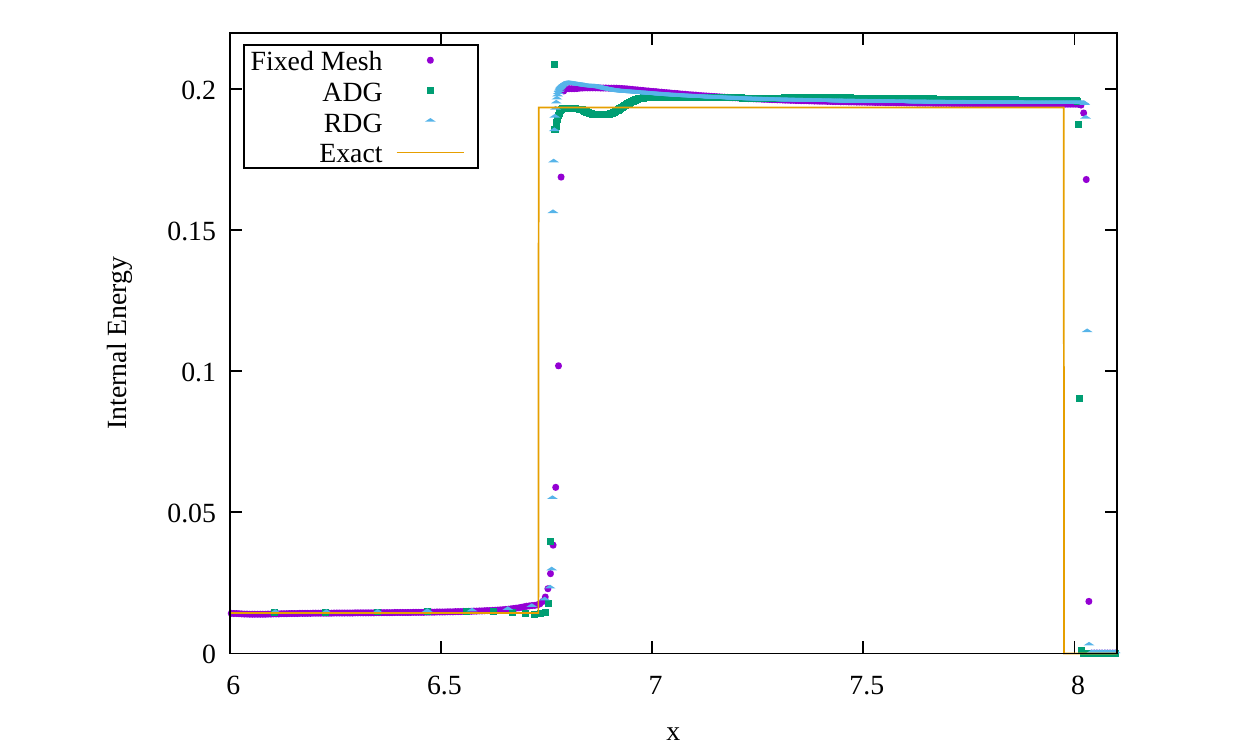}
\caption{Internal energy for Le Blanc Shock Tube with Rusanov flux, 1400 cells and TVD limiter, ADG
: Average Velocity, RDG : Linearized Riemann Velocity}
\label{fig: leblanc shock tube1}
\end{figure}

\begin{figure} \centering
\includegraphics[width=\textwidth]{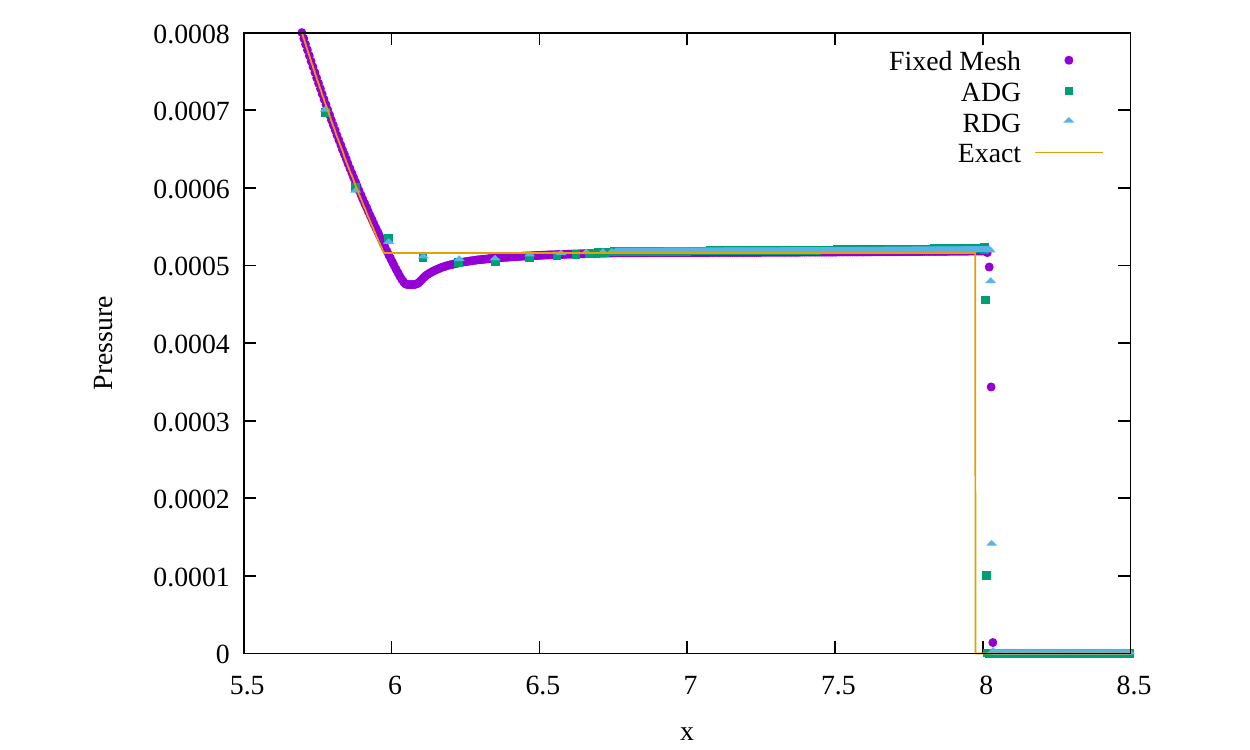}
\caption{Pressure for Le Blanc Shock Tube with Rusanov flux, 1400 cells and TVD limiter, ADG
  : Average Velocity, RDG : Linearized Riemann Velocity}
\label{fig: leblanc shock tube2}
\end{figure}

The Le Blanc shock tube test case is an extreme shock tube problem where the initial
discontinuity separates a region of high energy and density from one of low
energy and density. This is a much more severe test than the Sod problem and hence more challenging for numerical schemes. The
computational domain is $ 0 \leq x \leq 9$ and is filled with an ideal gas with
$\gamma = 5/3$. The gas is initially at rest and we perform the simulation up to a time of $T=6$ units.  The initial discontinuity is at $x=3$ and the initial condition is given by
\begin{align}
(\rho, v, p) = \begin{cases}
  (1.0, 0.0, 0.1) &\text{if } x < 3 \\
  (0.001, 0.0, 10^{-7}) &\text{if } x > 3
        \end{cases}
\end{align}
Note that both the density and pressure have a very large jump in the initial condition. The solution that develops from this initial condition consists of a
rarefaction wave moving to the left and a contact discontinuity and a strong shock
moving to the right. In Figure~(\ref{fig: leblanc shock tube1}), we show the
comparison of the internal energy profile at final time between a fixed mesh solution and moving mesh solutions with two
different mesh velocities as described before. Most methods tend to generate a very large spike in the internal energy in the contact region, e.g., compare with Figure~(11) in \cite{Shashkov2005}, while the present ALE method here is able to give a better profile. We plot the pressure profile in Figure~(\ref{fig: leblanc shock tube2})
 which shows that the ALE scheme is able to better represent the region around the contact wave as compared to fixed mesh method.

\color{red}

\subsection{Two Dimensional Isentropic Vortex Test Case}
\begin{figure}
\begin{center}
\includegraphics[width=0.48\textwidth]{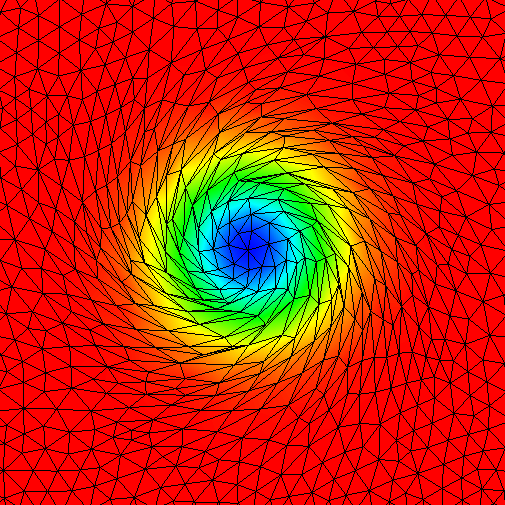}
\caption{Isentropic Vortex in 2-D: Skewed Mesh without Remeshing $t=2.660534$}
\label{fig:badisentropicvortex}
\end{center}
\end{figure}

\begin{figure}
\begin{center}
\begin{tabular}{cc}
\includegraphics[width=0.48\textwidth]{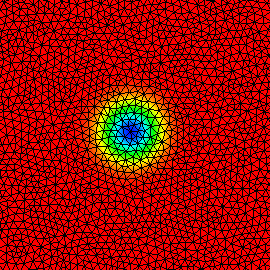}
&\includegraphics[width=0.48\textwidth]{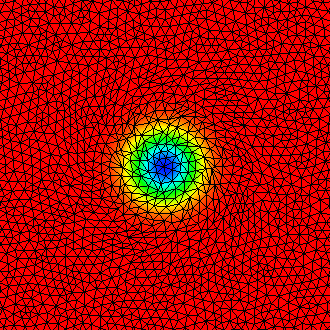}
\\
(a) & (b)\\
\includegraphics[width=0.48\textwidth]{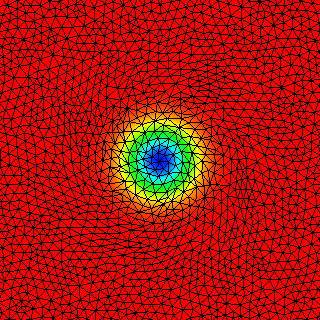}
& \includegraphics[width=0.48\textwidth]{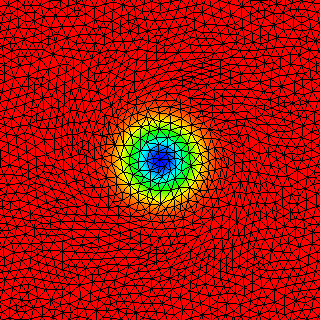}
\\
(c) & (d)
\end{tabular}
\caption{Isentropic vortex in 2-D: Mesh and pressure solution at various times (a)~$t=0$
(b)~$t=6$
(c)~$t=12$
(d)~$t=20$
}
\label{fig:isentropicvortex}
\end{center}
\end{figure}

The extension to two dimensions involves two aspects that need to be addressed. The first issue is how to handle the grid motion and the second is how to formulate the ALE-DG scheme. The second part is a natural generalization of the DG scheme we have described for the 1-D case in this paper, except that we have to construct basis functions on triangles and perform some numerical quadrature. The first part involving grid movement is more complicated and we present only very preliminary results in this section to demonstrate that the idea  has merit. We now consider the two dimensional Euler equations written as
\begin{align}
\frac{\partial u}{\partial t}
+ \frac{\partial f(u)}{\partial x}
+ \frac{\partial g(u)}{\partial y}&= 0
\end{align}
where
\begin{align}
  u &=
  \begin{bmatrix}
    \rho \\
    \rho u \\
    \rho v \\
    E
  \end{bmatrix}
 &
  f(u) &=
  \begin{bmatrix}
    \rho u \\
    p + \rho u^2 \\
    \rho u v \\
    (E+p) u
  \end{bmatrix}
 &
  g(u) &=
  \begin{bmatrix}
    \rho v \\
    \rho u v \\
    p + \rho v^2 \\
    (E+p) v
  \end{bmatrix}
\end{align}
\begin{align}
  p = (\gamma -1) \left[ E - \frac{1}{2} \rho (u^2 + v^2)\right]
\end{align}
The test case we consider involves an isentropic vortex that is advecting with constant velocity and is a smooth solution for which error norms can be calculated. The test is carried out on a square domain $[-10,10] \times [-10,10]$ with periodic boundary conditions. The initial conditions  is an isentropic vortex
\begin{align}
  T &= 1 - \frac{(\gamma-1)\beta ^2}{8 \gamma \pi^2} e^{1-r^2} \\
  \rho &= T^{\frac{1}{\gamma -1}} \\
  u &= u_{\infty} - \frac{\beta}{2\pi} ye^{\frac{1-r^2}{2}} \\
  v &= v_{\infty} - \frac{\beta}{2\pi} ye^{\frac{1-r^2}{2}} \\
  p &= \rho^\gamma
\end{align}
with $u_\infty=1, v_\infty=0, \gamma=1.4, \beta=10$. As the solution evolves in
time, the mesh becomes quite deformed because the vortex is continually shearing
the mesh, which can lead to degenerate meshes, as shown in
figure~\ref{fig:badisentropicvortex}. We
avoid the occurence of badly shaped triangles by using a combination of face
swapping and mesh velocity smoothing algorithms. The mesh modification is a very
local procedure and does not require global remeshing which is a costly process.
With these techniques, we are able to maintain a good mesh quality even after
the vortex has rotated 4 times around its center as shown in
figures~\ref{fig:isentropicvortex}. As the vortex is translating, we plot the
solution in a window centered at the vortex center. We can see that the method
maintains its high order of accuracy from the convergence rates of the error
shown in table~\ref{tab:isentropicvortex}; using linear basis functions yields
second order convergence while quadratic basis functions lead to third order
convergence.

\begin{table}
\begin{center}
\begin{tabular}{|c|c|c|c|c|} \hline
  \multirow{2}{*}{$N$}  & \multicolumn{2}{|c|}{$k=1$} &
  \multicolumn{2}{|c|}{$k=2$}
  \\ \cline{2-5} &
  Error & Rate & Error & Rate \\
  \hline
  50x50     &2.230e-03   &       &1.762E-04  &        \\
  100x100   &5.987E-04   &1.945  &2.305E-05  &2.934   \\
  200x200   &1.498E-04   &1.998  &2.973E-06  &2.955   \\
  400x400   &3.786E-05   &1.984  &3.762E-07  &2.982   \\
  800x800   &9.617E-06   &1.977  &3.474E-08  &2.991   \\
  \hline
\end{tabular}

\caption{Isentropic Vortex in 2D: Order of accuracy study on two dimensional static mesh}
\label{tab:isentropicvortexfixed}
\end{center}
\end{table}

\begin{table}
\begin{center}
\begin{tabular}{|c|c|c|c|c|} \hline
  \multirow{2}{*}{$N$}  & \multicolumn{2}{|c|}{$k=1$} &
  \multicolumn{2}{|c|}{$k=2$}
  \\ \cline{2-5} &
  Error & Rate & Error & Rate \\
  \hline
  50x50     &2.230e-03   &       &1.762E-04  &        \\
  100x100   &5.987E-04   &1.945  &2.305E-05  &2.934   \\
  200x200   &1.498E-04   &1.998  &2.973E-06  &2.955   \\
  400x400   &3.786E-05   &1.984  &3.762E-07  &2.982   \\
  800x800   &9.617E-06   &1.977  &3.474E-08  &2.991   \\
  \hline
\end{tabular}

\caption{Isentropic Vortex in 2D: Order of accuracy study on two dimensional moving mesh}
\label{tab:isentropicvortex}
\end{center}
\end{table}

\color{black}

\section{Summary and conclusions} We have developed an explicit DG scheme on
moving meshes using ALE framework and space-time expansion of the solutions
within each cell. The near Lagrangian nature of the mesh motion dramatically
reduces the numerical dissipation especially for contact waves. Even moving
contact waves can be exactly computed with a numerical flux that is exact for
stationary contact waves. The scheme is shown to yield superior results even in
the presence of large boost velocity of the coordinate system indicating its
Galilean invariance property. The standard Roe flux does not suffer from entropy violation when applied in the current nearly Lagrangian framework. However, in some problems with strong shocks, spurious contact waves can appear and we propose to fix the dissipation in Roe-type schemes that eliminates this issue. The method yields accurate solutions even in
combination with standard TVD limiters, where fixed grid methods perform poorly.
The mesh motion provides automatic grid adaptation near shocks but may lead to
very coarse cells inside expansion waves. A grid adaptation strategy is
developed to handle the problem of very small or very large cells. The presence
of the DG polynomials makes it easy to transfer the solution during grid
adaptation without loss of accuracy. The proposed methodology is general enough
to be applicable to other systems of conservation laws modelling fluid flows. The
basic idea can be extended to multi-dimensions but additional considerations are
required to maintain good mesh quality under fluid deformations. The preliminary
results shown for the isentropic vortex is very promising for the 2-D case.
\section*{Acknowledgements} This work was started when Jayesh Badwaik was a project assistant at TIFR-CAM, Bangalore. The first two authors gratefully acknowledge the financial support received from the Airbus Foundation Chair on Mathematics of Complex Systems established in TIFR-CAM, Bangalore, for carrying out this work. Christian Klingenberg acknowledges the support of the {\em Priority Program 1648: Software for Exascale Computing} by the German Science Foundation.
On behalf of all authors, the corresponding author states that there is no conflict of interest.
\appendix
\section{Numerical flux} The ALE scheme requires a numerical flux
$\nfl(\con_l,\con_r,w)$ which is usually based on some approximate Riemann
solver.  The numerical flux function is assumed to be consistent in the sense
that \[ \nfl(\con, \con, w) = \mfl(\con,w), \qquad \forall \ \con \in \re^3, w
  \in \re \] Since the ALE versions of the numerical fluxes are not so well
  known, here we list the formulae used in the present
  work.
  \subsection{Rusanov flux} \label{sec:rusanov} The Rusanov flux is a variant of
  the Lax-Friedrich flux and is given by \[ \nfl(\con_l,\con_r,w) = \half[
  \mfl(\con_l,w) + \mfl(\con_r,w)] - \half \lambda_{lr} (\con_r - \con_l) \]
  where $\lambda_{lr} = \lambda(\con_l,\con_r,w)$ \[ \lambda(\con_l,\con_r,w) =
  \max\{ |v_l - w|+c_l, |v_r - w|+c_r \} \] which is an estimate of the largest
  wave speed in the Riemann problem. Since the mesh velocity is close to the
  fluid velocity, the value of $\lambda$ is close to the local sound speed. Thus
  the numerical dissipation is independent of the velocity scale.
\subsection{Roe flux} The Roe scheme~\cite{Roe1981357} is based on a local
linearization of the conservation law and then exactly solving the Riemann
problem for the linear approximation. The flux can be written as \[
  \nfl(\con_l,\con_r,w) = \half[ \mfl(\con_l,w) + \mfl(\con_r,w)] - \half |A_w|
  (\con_r - \con_l) \] where the Roe average matrix $A_w = A_w(\con_l,\con_r)$
  satisfies \[ \mfl(\con_r,w) - \mfl(\con_l,w) = A_w (\con_r - \con_l) \] and we
  define $|A_w| = R |\Lambda -wI| R^{-1}$. This matrix is evaluated at the Roe
  average state $\con(\bar{\pv})$, $\bar{\pv}=\half(\pv_l + \pv_r)$, where
  $\pv=\sqrt{\rho}[1, \ v, \ H]^\top$ is the parameter vector introduced by Roe.
\subsection{HLLC flux} This is based on a three wave approximate Riemann solver
and the particular ALE version we use can also be found in~\cite{Luo2004304}.
Define the relative velocity $\rv = v - w$; then the numerical flux is given by
\[ \nfl(\con_l,\con_r,w) = \begin{cases} \mfl(\con_l,w) & S_l > 0 \\
    \mfl^*(\con_l^*,w) & S_l \le 0 < S_M \\ \mfl^*(\con_r^*,w) & S_M \le 0 \le
    S_r \\ \mfl(\con_r,w) & S_r < 0 \end{cases} \] where the intermediate states
    are given by \[ \con_\alpha^* = \frac{1}{S_\alpha - S_M} \begin{bmatrix}
      (S_\alpha - \rv_\alpha) \rho_\alpha \\ (S_\alpha - \rv_\alpha)(\rho
    v)_\alpha + p^* - p_\alpha \\ (S_\alpha - \rv_\alpha) E_\alpha - p_\alpha
  \rv_\alpha + p^* S_M \end{bmatrix}, \qquad \alpha=l,r \] and \[ \mfl^*(\con,w)
  = S_M \con + \begin{bmatrix} 0 \\ p^* \\ (S_M + w) p^* \end{bmatrix} \] where
  \[ p^* = \rho_l (\rv_l - S_l) (\rv_l - S_M) + p_l = \rho_r (\rv_r - S_r)(\rv_r
    - S_M) + p_r \] which gives $S_M$ as \[ S_M = \frac{ \rho_r \rv_r (S_r -
    \rv_r) - \rho_l \rv_l (S_l - \rv_l) + p_l - p_r}{\rho_r (S_r - \rv_r) -
  \rho_l (S_l - \rv_l)} \] The signal velocities are defined as \[ S_l = \min\{
\rv_l - c_l, \hat{v}-w-\hat{c}\}, \qquad S_r = \max\{ \rv_r + c_r,
\hat{v}-w+\hat{c}\} \] where $\hat{v}$, $\hat{c}$ are Roe's average velocity and
speed of sound.

\section{Continuous expansion Runge-Kutta (CERK) schemes} We use a Runge-Kutta
scheme to compute the predicted solution used to compute all the integrals in
the DG scheme. In this section, we list down the CERK scheme for the following
ODE \[ \dd{u}{t} = f(u,t) \] Given the solution $u^n$ at time $t_n$, the CERK
scheme gives a polynomial solution in the time interval $[t_n, t_{n+1})$ of the
form \[ u(t_n + \theta h) = u^n + h \sum_{s=1}^{n_s} b_s(\theta) k_s, \qquad
  \theta \in [0,1] \] where $n_s$ is the number of stages and $h$ denotes the
  time step.

\subsection{Second order (CERK2)} The number of stages is $n_s = 2$ and \[
b_1(\theta) = \theta - \theta^2/2, \qquad b_2(\theta) = \theta^2/2 \] and
\begin{eqnarray*} k_1 &=& f(u^n, t_n) \\ k_2 &=& f(u^n + h k_1, t_n + h)
\end{eqnarray*}

\subsection{Third order (CERK3)} The number of stages is $n_s = 4$ and
  \begin{eqnarray*} k_1 &=& f(u^n, t_n) \\ k_2 &=& f(u^n + (12/23) h k_1, t_n +
    12 h /23) \\ k_3 &=& f(u^n + h ((-68/375) k_1 + (368/375) k_2), t_n + 4h/5)
    \\ k_4 &=& f(u^n + h ((31/144) k_1 + (529/1152) k_2 + (125/384) k_3), t_n +
    h) \end{eqnarray*} and \begin{eqnarray*} b_1(\theta) &=& (41/72) \theta^3 -
    (65/48) \theta^2 + \theta \\ b_2(\theta) &=& -(529/576) \theta^3 + (529/384)
    \theta^2 \\ b_3(\theta) &=& -(125/192) \theta^3 + (125/128) \theta^2 \\
    b_4(\theta) &=& \theta^3 - \theta^2 \end{eqnarray*}

\bibliography{bibdesk}

\end{document}